\documentclass{article}

\title{Adding a constant and an axiom to a doctrine}
\author{Francesca Guffanti\footnote{University of Luxembourg, Luxembourg. Email address: francesca.guffanti@uni.lu}\,\,\footnote{Work conducted at Università degli Studi di Milano, Italy.}}
\date{}
\usepackage[utf8]{inputenc}
\usepackage{csquotes}
\usepackage[italian,british]{babel}
\usepackage[margin=89.5pt]{geometry}
\usepackage[pdfa]{hyperref}

\hypersetup{
    pdftitle={},
    pdfauthor={},
    pdfmenubar=false,
    pdffitwindow=true,
    pdfstartview=FitH,
    colorlinks=true,
    linkcolor=black,
    citecolor=black,
    urlcolor=cyan
}

\usepackage{todonotes}
\usepackage{setspace}
\usepackage{mathtools}
\usepackage{amssymb}
\usepackage{bbold}
\usepackage{amsmath,amsthm, graphicx,xspace}
\usepackage{mathrsfs}
\usepackage{latexsym}
\usepackage{rotating}
\usepackage{multicol}
\usepackage{verbatim}
\usepackage{faktor}
\usepackage{epigraph}
\usepackage{quiver}
\usepackage{esvect}
\usepackage{cancel}
\usepackage{oubraces}
\usepackage{bbm}
\usepackage[capitalize]{cleveref}
\usepackage{enumitem}
\setlist[itemize]{parsep=0pt}
\setlist[enumerate]{parsep=0pt}
\usepackage{tikz-cd}
\usepackage{tikz}
\usetikzlibrary{arrows}

\numberwithin{equation}{section}
\theoremstyle{definition}
\newtheorem{theorem}{Theorem}[section]
\newtheorem{proposition}[theorem]{Proposition}

\newtheorem{corollary}[theorem]{Corollary}
\newtheorem{definition}[theorem]{Definition}
\newtheorem{notation}[theorem]{Notation}
\newtheorem{example}[theorem]{Example}
\newtheorem{remark}[theorem]{Remark}
\mathcode`\:="603A 
\mathchardef\colon="303A 

\def\pws{\mathop{\mathscr{P}\kern-1.1ex_{\ast}}}
\DeclareMathOperator{\OPRid}{id}
\newcommand{\id}[1]{\ensuremath{\OPRid_{#1}}}
\newcommand{\ct}[1]{\ensuremath{\mathbb{#1}}\xspace}
\newcommand{\CC}{\ct{C}}
\newcommand{\Pos}{\mathbf{Pos}}
\newcommand{\Cat}{\mathbf{Cat}}
\newcommand{\Hom}{\mathop{\mathrm{Hom}}}
\newcommand{\QDot}{\mathbf{IdxPos}}
\newcommand{\Dott}{\mathbf{Dct}}
\DeclareMathOperator{\Cmd}{Cmd}
\DeclareMathOperator{\inc}{Inc}
\newcommand{\ple}[1]{\langle#1\rangle}
\newcommand{\op}{^{\mathrm{op}}}
\newcommand{\blank}{\mathrm{-}}
\newcommand{\bblank}{\mathrm{=}}
\newcommand{\pr}[1]{\ensuremath{\mathrm{pr}_{#1}}}
\newcommand{\tmn}{\mathbf{t}}

\newcommand{\PD}{\mathbf{PD}}

\newcommand{\Sub}{\mathop{\mathrm{Sub}}}

\newcommand{\Set}{\textnormal{Set}}

\newcommand{\des}{\mathscr{D}es}
\newcommand{\dom}{\mathop{\mathrm{dom}}}

\newcommand{\LT}{\mathtt{LT}}
\newcommand{\ctx}{\mathrm{\mathbb{C}tx}}

\newcommand{\ter}{\mathsf{1}}
\newcommand{\ini}{\mathsf{0}}
\newcommand{\pwo}[1]{\mathsf{P}(#1)}
\begin{document}
\newlength{\myindent} 
\setlength{\myindent}{\parindent}
\parindent 0em 
\maketitle
\begin{abstract}
We study the meaning of “adding a constant to a language” for any doctrine, and “adding an axiom to a theory” for a primary doctrine, by showing how these are actually two instances of the same construction. We prove their universal properties, and how these constructions are compatible with additional structure on the doctrine. Existence of Kleisli object for comonads in the 2-category of indexed poset is proved in order to build these constructions.
\end{abstract}
\setcounter{tocdepth}{1}
\tableofcontents

\section{Introduction}
Given a theory $\mathcal{T}$ in a first-order language $\mathcal{L}$, consider for each set of variables the set of well-formed formulae written with at most those variables. That set can be ordered by provable consequence in the theory $\mathcal{T}$. By that we mean that the formula $\alpha$ is less than or equal to the formula $\beta$ if the consequence $\alpha\vdash_\mathcal{T}\beta$ holds. The logical operations of conjunction, disjunction, implication, negation, true and false give this set the structure of a Boolean algebra. And the assignment $\vec x\mapsto\mathrm{WFF}(\vec x)$ of the Boolean algebra of formulae to a list (i.e.\ a context) of distinct variables can be extended to a functor $\mathrm{WFF}:\ctx\op\to\mathbf{BA}$ from the opposite category of contexts and terms to the category of Boolean algebras and homomorphisms.

This can be considered the motivating example at the basis of the notion of hyperdoctrine which was introduced by Lawvere in 1969 in a series of seminal papers \cite{adjfound,lawdiag,lawequality}. It is a categorical tool that allows the analysis of both syntax and semantics of logical theories through the same mathematical structure. One of the main intuitions of Lawvere was to recognize that quantifiers in logic are instances of adjunctions between the posets of formulae. 

A doctrine is possibly the basic fabric of Lawvere's hyperdoctrine: just a functor $P:\CC\op\to\Pos$, from a category $\CC$ with finite products into the category $\Pos$ of partially ordered sets and monotone functions.

In this framework, a possible direction is to understand what classical results in logic can be interpreted in the language of doctrines, and how to do so. In first-order model theory it is not uncommon for proofs to rely on adding constant symbols to a language and axioms to a theory. For example, the proof of Henkin's Theorem \cite{henkin}---\emph{every consistent first-order theory has a model}--- is done by at first extending the language with new constant symbols, and then adding new axioms of the form $\exists x\varphi(x)\to\varphi(c)$ to the theory. Henkin's Theorem is now the standard way to prove one of the most crucial results in first-order logic: Gödel’s Completeness Theorem \cite{godel}---\emph{if a closed formula is valid in every first-order structure, then there is a first-order proof of it}. A similar process can be used in the proof of the Upward Löwenheim-Skolem-Tarski Theorem \cite{modeltheory}---\emph{a theory which has an infinite model has a model in every cardinality $\kappa\geq\max(|\mathcal{L}|, \aleph_0)$}: after having extended the language $\mathcal{L}$ with a suitable amount of constant symbols, axioms of the form $\lnot(c=d)$ are added to the theory. Many other examples arise in algebra: the theory of pointed sets is obtained from the theory of sets by adding one constant symbol, the theory of abelian groups is obtained from the theory of groups by adding the commutativity axiom, the theory of monoids is obtained from the theory of semigroups by adding a constant $e$ and the axiom $\forall x((x\cdot e = x)\land(e\cdot x = x))$, and so on.

Thus, translating the concepts of adding constants and axioms to a theory is a tool which could be relevant for future applications, with the aim to interpret results of first-order logic in the doctrinal setting. However, before understanding what it means to add \emph{any} amount of constants and axioms to a doctrine, it is reasonable to start this investigation in the finitary case, which is the focus of this work.

In more detail, we extend to doctrines the construction of adding a constant to a language. 
And we also extend the construction of forcing a new axiom for primary doctrines, which are doctrines where all orders are inf-semilattices and reindexing preserves them---essentially, what it amounts to the ability to interpret conjunctions of formulae.
Actually, we do both constructions in one step, using a Kleisli object for a convenient comonad on the original doctrine seen as an indexed poset. The existence of Kleisli objects for comonads in indexed posets is proven in Proposition \ref{prop:kleisli}.
Let $P:\CC\op\to\Pos$ be a primary doctrine, let $X$ be a fixed object, and let $\varphi\in P(X)$ be an element. In keeping a reasonable parallel with the logical intuition, think of the primary doctrine $P$ as the syntactic consequences of a theory---not just formulae of a language, but that is already a good intuition---; think of an object $X$ as a list of fresh variables, and think of $\varphi$ as a formula in the fresh variables. We construct a primary doctrine $P_{(X,\varphi)}$ and a homomorphism of primary doctrines $P\to P_{(X,\varphi)}$. Staying with the parallel, the doctrine $P_{(X,\varphi)}$ acts like the extension of the theory with new constant symbols and with the new axiom $\varphi$ evaluated in those constants. The homomorphism of doctrines acts like a translation of the original theory into the new one.

The construction has a universal property: any other morphism of primary doctrines $P\to R$ such that the interpretation of $\varphi$ evaluated in some constant in $R$ is true factors through $P\to P_{(X,\varphi)}$, essentially in a unique way. This is, in broad terms, the statement of Theorem \ref{thm:univ_prop}. Moreover, the result is extended to 2-arrows in Theorem \ref{thm:prec}. We also show throughout \cref{sect:pres_prop} that the construction $P\to P_{(X,\varphi)}$ preserves many additional structures and properties that the original doctrine may already enjoy.

The construction includes the two constructions we discussed at the
beginning: adding no axiom has a structural parallel in adding $\top$
to the axioms, while adding no constant corresponds to performing the
construction picking the terminal object for $X$.
We can clearly decide to add just a constant of sort $X$,
without adding any axiom to the theory, and obtain a homomorphism $P\to
P_X$.
Similarly, we can decide to add just an axiom $\varphi$ to a
theory, without adding any constant symbol, and obtain the homomorphism
$P\to P_\varphi$.

\section{Preliminaries}
In this section, we lay the groundwork for the paper by introducing the language of doctrines and establishing their key properties. We then compute the Kleisli objects for comonads in the 2-category of indexed posets, which will be the basis for our main construction. Most of the notions and results concerning category theory used in this thesis are standard, and we refer to any textbook, for instance \cite{borhandbk,elephant,CatWorMat}.

\subsection{Doctrines}
In this subsection, we define the 2-category of doctrines and show some relevant examples. Then we will gradually add more structure in order to be able to interpret symbols of first-order logic---such as connectives and quantifiers---in the context of doctrines. 
\begin{definition}
Let $\ct{C}$ be a category with finite products and let $\Pos$ be the category of partially-ordered sets and monotone functions. A \emph{doctrine} is a functor $P:\ct{C}\op\to\Pos$. The category $\CC$ is called \emph{base category of $P$}, each poset $P(X)$ for an object $X\in\CC$ is called \emph{fiber}, the function $P(f)$ for an arrow $f$ in $\CC$ is called \emph{reindexing}.
\end{definition}
By viewing doctrines as a broad generalization of doctrines of well-formed formulae, we can interpret the objects of category $\CC$ as lists of variables, the arrows as terms, the fibers as sets of formulae, and reindexing as substitutions, providing an intuitive understanding of the structure.
\begin{example}\label{ex:doctr}We propose the following examples.
\begin{enumerate}[label=(\alph*)]
\item The functor $\mathscr{P}:\Set\op\to\Pos$, sending each set in the poset of its subsets, ordered by inclusion, and each function $f:A\to B$ to the inverse image $f^{-1}:\mathscr{P}(B)\to\mathscr{P}(A)$ is a doctrine.
\item For a given category $\CC$ with finite limits, the functor $\Sub_\CC:\CC\op\to\Pos$ sending each object to the poset of its subobjects in $\CC$ and each arrow $f:A\to B$ to the pullback function $f^*:\Sub_\CC(B)\to\Sub_\CC(A)$, is a doctrine.
\item \label{i:fbf} For a given theory $\mathcal{T}$ on a one-sorted first-order language $\mathcal{L}$, define the category $\ctx_\mathcal{L}$ of contexts: an object is a finite list of distinct variables and an arrow between two lists $\vec x=(x_1,\dots, x_n)$ and $\vec y=(y_1,\dots, y_m)$ is
\begin{equation*}(t_1(\vec x),\dots,t_m(\vec x)):(x_1,\dots, x_n) \to(y_1,\dots, y_m)\end{equation*}
an $m$-tuple of terms in the context $\vec x$. The empty list $()$ is the terminal object in $\ctx_\mathcal{L}$, the product of two lists $\vec x$ and $\vec y$ in $\ctx_\mathcal{L}$ is given by any list whose length is the sum of the length of $\vec x$ and $\vec y$---if the variables in the two lists are all distinct, their product can be written as the juxtaposition $\langle\vec x;\vec y\rangle=(x_1\dots,x_n,y_1,\dots,y_m)$. The functor $\LT^\mathcal{L}_{\mathcal{T}}:\ctx_\mathcal{L}\op\to\Pos$ sends each list of variables $\vec{x}$ to the Lindenbaum-Tarski algebra of well-formed formulae with free variables in $\vec{x}$ modulo the theory $\mathcal{T}$, meaning that $\LT^\mathcal{L}_{\mathcal{T}}(\vec{x})$ is the poset reflection of well-formed formulae written with free variables in $\vec{x}$, possibly dummy, ordered by provable consequence in $\mathcal{T}$. Moreover, $\LT^\mathcal{L}_{\mathcal{T}}:\ctx_\mathcal{L}\op\to\Pos$ sends an arrow $\vec{t}(\vec{x}):\vec{x}\to\vec{y}$ into the substitution $[\vec{t}(\vec{x})/\vec{y}]$, that maps the equivalence class of a formula $\alpha(\vec y)$ to the equivalence class of the formula $\alpha(\vec{t}(\vec{x})/\vec{y})$. We refer to any standard textbook about first-order logic for definitions of concepts including language, variables, theory, terms, substitution, formulae, see for instance \cite{modth}.
\item For a given category $\ct{D}$ with finite products and weak pullbacks, the functor of weak subobjects $\mathbf{\Psi}_\ct{D}:\ct{D}\op\to\Pos$ sending each object $A$ to the poset reflection of the comma category $\ct{D}/A$ is a doctrine: for each arrow $f:A\to B$, $\mathbf{\Psi}_\ct{D}(f)$ sends the equivalence class of an arrow $\alpha:\dom\alpha\to B$ to the equivalence class of the projection $\pi_1$ of a chosen weak pullback of $\alpha$ along $f$---see Example 2.9 in \cite{quotcomplfoun} for more details.
\[\begin{tikzcd}
	W & \dom\alpha \\
	A & B
	\arrow["{\pi_1}"', from=1-1, to=2-1]
	\arrow["{\pi_2}", from=1-1, to=1-2]
	\arrow["\alpha", from=1-2, to=2-2]
	\arrow["f"', from=2-1, to=2-2]
\end{tikzcd}\]
\end{enumerate}
\end{example}
\begin{definition}
A \emph{doctrine homomorphism}---1-cell or 1-arrow---between $P:\ct{C}\op\to\Pos$ and $R:\ct{D}\op\to\Pos$ is a pair $(F,\mathfrak f)$ where $F:\ct{C}\to\ct{D}$ is a functor that preserves finite products and $\mathfrak{f}:P\xrightarrow{\cdot}R\circ F\op$ is a natural transformation. Sometimes a morphism between $P$ and $R$ will be called a model of $P$ in $R$.
A \emph{2-cell} between $(F,\mathfrak f)$ and $(G,\mathfrak g)$ from $P$ to $R$ is a natural transformation $\theta:F\xrightarrow{\cdot} G$ such that $\mathfrak{f}_A(\alpha)\leq R(\theta_A)(\mathfrak{g}_A(\alpha))$ for any object $A$ in $\ct{C}$ and $\alpha\in P(A)$. Doctrine, doctrine morphisms with 2-cells defined here form a $2$-category, that will be denoted $\Dott$.
\end{definition}
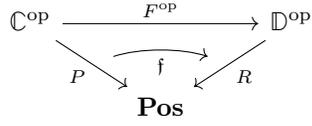
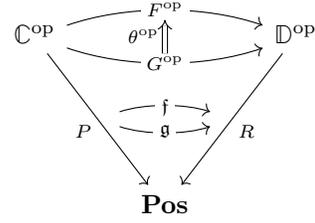
\[\begin{tikzcd}
	{\ct{C}\op} && {\ct{D}\op} && {\ct{C}\op} && {\ct{D}\op} \\
	& \Pos \\
	&&&&& \Pos
	\arrow["F\op", from=1-1, to=1-3]
	\arrow[""{name=0, anchor=center, inner sep=0}, "P"', from=1-1, to=2-2]
	\arrow[""{name=1, anchor=center, inner sep=0}, "R", from=1-3, to=2-2]
	\arrow[""{name=2, anchor=center, inner sep=0}, "F\op"{description}, curve={height=-12pt}, from=1-5, to=1-7]
	\arrow[""{name=3, anchor=center, inner sep=0}, "P"', from=1-5, to=3-6]
	\arrow[""{name=4, anchor=center, inner sep=0}, "R", from=1-7, to=3-6]
	\arrow[""{name=5, anchor=center, inner sep=0}, "G\op"{description}, curve={height=12pt}, from=1-5, to=1-7]
	\arrow["{\mathfrak{f}}"', curve={height=-6pt}, shorten <=8pt, shorten >=8pt, from=0, to=1]
	\arrow["{\mathfrak f}"{description}, curve={height=-6pt}, shorten <=8pt, shorten >=8pt, from=3, to=4]
	\arrow["{\mathfrak g}"{description}, curve={height=6pt}, shorten <=8pt, shorten >=8pt, from=3, to=4]
	\arrow["\theta\op", shorten <=3pt, shorten >=3pt, Rightarrow, from=5, to=2]
\end{tikzcd}\]
By definition of doctrine, the fibers are simply posets. However, we can define specific doctrines by imposing additional structure on these posets or by requiring the existence of adjoints to certain reindexing. To work in a setting that interprets the conjunction of formulae and the true constant, primary doctrines are necessary.
\begin{definition}
A \emph{primary doctrine} $P:\ct{C}\op\to\Pos$ is a doctrine such that for each object $A$ in $\ct{C}$, the poset $P(A)$ has finite meets, and the related operations $\land:P\times P\xrightarrow{\cdot} P$ and $\top:\mathbf{1}\xrightarrow{\cdot} P$ yield natural transformations.
\end{definition}
\begin{example}
Examples seen in \ref{ex:doctr} are primary doctrines:
\begin{enumerate}[label=(\alph*)]
\item For any set $A$, intersection of two subsets is their meet, $A$ is the top element.
\item For any object $A$ in $\CC$, the pullback of a subobject along another defines their meet.
\[\begin{tikzcd}
	{\dom(\alpha\land\beta)} & \dom\alpha \\
	\dom\beta & A
	\arrow["{\pi_1}"', tail, from=1-1, to=2-1]
	\arrow["{\pi_2}", tail, from=1-1, to=1-2]
	\arrow["\alpha", tail, from=1-2, to=2-2]
	\arrow["\beta"', tail, from=2-1, to=2-2]
	\arrow["\alpha\land\beta"{description}, tail, from=1-1, to=2-2]
\end{tikzcd}\]
The arrow $\id{A}$ is the top element.
\item For any list $\vec x$, the conjunction of two formulae is their binary meet, the true constant $\top$ is the top element.
\item For any object $A$ in $\ct{D}$, a choice of a weak pullback of a representative of a weak subobject along another defines their meet
\[\begin{tikzcd}
	{\dom(\alpha\land\beta)} & \dom\alpha \\
	\dom\beta & A
	\arrow["{\pi_1}"', from=1-1, to=2-1]
	\arrow["{\pi_2}", from=1-1, to=1-2]
	\arrow["\alpha", from=1-2, to=2-2]
	\arrow["\beta"', from=2-1, to=2-2]
	\arrow["\alpha\land\beta"{description}, from=1-1, to=2-2]
\end{tikzcd},\]
the class of $\id{A}$ is the top element.
\end{enumerate}
\end{example}
In order to interpret equality, we recall the notion of elementary doctrines. The definition we provide is not the classical one given by Lawvere in \cite{lawequality}, but it is the equivalent version that can be found in Proposition 2.5 of \cite{EmPaRo}:
\begin{definition}\label{def:elem}
A primary doctrine $P:\ct{C}\op\to\Pos$ is \emph{elementary} if for any object $A$ in $\ct{C}$ there exists an element $\delta_A\in P(A\times A)$ such that:
\begin{enumerate}
\item $\top_A\leq P(\Delta_A)(\delta_A)$;
\item $P(A)=\des_{\delta_A}:=\{\alpha\in P(A)\mid P(\pr1)(\alpha)\land\delta_A\leq P(\pr2)(\alpha)\}$;
\item $\delta_A\boxtimes\delta_B\leq\delta_{A\times B}$, where $\delta_A\boxtimes\delta_B=P(\ple{\pr1,\pr3})(\delta_A)\land P(\ple{\pr2,\pr4})(\delta_B)$.
\end{enumerate}
In 2., $\pr1$ and $\pr2$ are the projections from $A\times A$ in $A$; in 3., the projections are from $A\times B\times A\times B$. The element $\delta_A$ will be called \emph{fibered equality} on $A$.
\end{definition}
\begin{example}
Examples seen in Example \ref{ex:doctr} are elementary doctrines:
\begin{enumerate}[label=(\alph*)]
\item For any set $A$, the subset $\Delta_A\subseteq A\times A$ is the fibered equality on $A$.
\item For any object $A$ in $\CC$, the map $\Delta_A:A\rightarrowtail A\times A$ is the fibered equality on $A$---see in \cite{ElemQuotCompl} the Example 2.4.a.
\item Provided that the language $\mathcal{L}$ has equality, for any list $\vec x$, the formula $\big(x_1=x_1'\land\dots\land x_n=x_n'\big)$ in $\LT^\mathcal{L}_\mathcal{T}(\vec x;\vec x')$ is the fibered equality on $\vec{x}$.
\item For any object $A$ in $\ct{D}$, the equivalence class of the map $\Delta_A:A\rightarrowtail A\times A$ is the fibered equality on $A$.
\end{enumerate}
\end{example}
We now generalize the existential and universal quantifiers, which are defined as adjoint to some reindexing.
\begin{definition}\label{def:exist}
A primary doctrine $P:\ct{C}\op\to\Pos$ is \emph{existential} if for any pair of objects $B, C$ of $\ct{C}$, the map $P(\pr1):P(C)\to P(C\times B)$ has a left adjoint
\begin{equation*}\exists^B_C:P(C\times B)\to P(C),\end{equation*}
satisfying:
\begin{itemize}
\item[-] Beck-Chevalley condition with respect to pullback diagrams of the form:
\[\begin{tikzcd}[ampersand replacement=\&]
	{C\times B} \& C \\
	{C'\times B} \& {C'}
	\arrow["{f\times \id{B}}"', from=1-1, to=2-1]
	\arrow["f", from=1-2, to=2-2]
	\arrow["\pr1"', from=2-1, to=2-2]
	\arrow["\lrcorner"{anchor=center, pos=0.125}, draw=none, from=1-1, to=2-2]
	\arrow["\pr1", from=1-1, to=1-2]
\end{tikzcd}\]
that is, $\exists^B_CP(f\times\id{B})=P(f)\exists^B_{C'}$.
\item[-] Frobenius reciprocity, that is, for any $\alpha\in P(C\times B)$ and $\beta\in P(C)$ it holds $\exists^B_C(\alpha\land P(\pr1)(\beta))=\exists^B_C(\alpha)\land \beta$.
\end{itemize}
\end{definition}
\begin{definition}
A doctrine $P:\ct{C}\op\to\Pos$ is \emph{universal} if for any pair of objects $B, C$ of $\ct{C}$, the map $P(\pr1):P(C)\to P(C\times B)$ has a right adjoint $\forall^B_C:P(C\times B)\to P(C)$, satisfying Beck-Chevalley condition with respect to pullback diagrams of the same form as Definition \ref{def:exist}, that is $P(f)\forall^B_{C'}=\forall^B_CP(f\times\id{B})$.
\end{definition}
\begin{definition}
A doctrine $P:\ct{C}\op\to\Pos$:
\begin{itemize}
\item is \emph{implicational} if for any object $A$, the poset $P(A)$ is cartesian closed, and the related operations $\land:P\times P\xrightarrow{\cdot} P$, $\top:\mathbf{1}\xrightarrow{\cdot} P$, $\to:P\op\times P\xrightarrow{\cdot} P$ yield natural transformations--in particular it is a primary doctrine;
\item \emph{has bottom element} if for any object $A$, the poset $P(A)$ has a bottom element, and the related operation, $\bot:\mathbf{1}\xrightarrow{\cdot} P$ yields a natural transformation;
\item \emph{is bounded} if for any object $A$, the poset $P(A)$ has a top and a bottom element, and the related operation, $\top:\mathbf{1}\xrightarrow{\cdot} P$ and $\bot:\mathbf{1}\xrightarrow{\cdot} P$ yield natural transformations;
\item \emph{has finite joins} if for any object $A$, the poset $P(A)$ has finite joins, and the related operations $\lor:P\times P\xrightarrow{\cdot} P$, $\bot:\mathbf{1}\xrightarrow{\cdot} P$ yield natural transformations;
\item is \emph{Heyting} if for any object $A$, the poset $P(A)$ is an Heyting algebra, and the related operations $\land:P\times P\xrightarrow{\cdot} P$, $\top:\mathbf{1}\xrightarrow{\cdot} P$, $\to:P\op\times P\xrightarrow{\cdot} P$, $\lor:P\times P\xrightarrow{\cdot} P$, $\bot:\mathbf{1}\xrightarrow{\cdot} P$ yield natural transformations;
\item is \emph{Boolean} if it is Heyting and the operation $\lnot(\blank):=(\blank)\to\bot:P\op\xrightarrow{\cdot}P$ is an isomorphism. 
\end{itemize}
\end{definition}
\begin{example}
The doctrine $\LT^\mathcal{L}_{\mathcal{T}}:\ctx_\mathcal{L}\op\to\Pos$ is Boolean elementary existential universal: in addition to the structure mentioned in the examples above, the implication of two formulae gives the implicational structure, the disjunction of two formulae is their join, the false is the bottom element, existential and universal quantifier define the left and the right adjoint to the inclusions of formulae $\LT^\mathcal{L}_{\mathcal{T}}(\vec x)\subseteq\LT^\mathcal{L}_{\mathcal{T}}(\vec x;\vec y)$ for any pair $\vec x=(x_1,\dots, x_n)$, $\vec y=(y_1,\dots, y_m)$:
\[\begin{tikzcd}
	{\LT^\mathcal{L}_{\mathcal{T}}(\vec x;\vec y)} & {\LT^\mathcal{L}_{\mathcal{T}}(\vec x)}
	\arrow[""{name=0, anchor=center, inner sep=0}, from=1-2, to=1-1, hook']
	\arrow[""{name=1, anchor=center, inner sep=0}, "{\forall y_1\dots\forall y_m}"', shift right=4, from=1-1, to=1-2]
	\arrow[""{name=2, anchor=center, inner sep=0}, "{\exists y_1\dots\exists y_m}", shift left=4, from=1-1, to=1-2]
	\arrow["\dashv"{anchor=center, rotate=-90}, draw=none, from=0, to=1]
	\arrow["\dashv"{anchor=center, rotate=-90}, draw=none, from=2, to=0]
\end{tikzcd}.\]
\end{example}
\begin{definition}
Any homomorphism $(F,\mathfrak f):P\to R$ from $P:\ct{C}\op\to\Pos$ to $R:\ct{D}\op\to\Pos$ is called respectively \emph{primary, elementary, existential, universal, implicational, bounded, Heyting, Boolean} if both $P$ and $R$ are, and $\mathfrak{f}$ preserves the said structure. 
For example, an elementary homomorphism is such that for any object $A$ in $\ct{C}$, and any $\alpha,\alpha'\in P(A)$:
\begin{equation*}\mathfrak{f}_A(\alpha\land_A\alpha')=\mathfrak{f}_A(\alpha)\land_{FA}\mathfrak{f}(\alpha');\qquad\qquad\mathfrak{f}_A(\top_A)=\top_{FA};\qquad\qquad\mathfrak{f}_{A\times A}(\delta_A)=\delta_{FA};\end{equation*}
while a universal homomorphism is such that for any pair of objects $B,C$ in $\ct{C}$, and any element $\alpha\in P(C\times B)$:
\begin{equation*}\mathfrak{f}_C\forall^B_C(\alpha)=\forall^{FB}_{FC}\mathfrak{f}_{C\times B}(\alpha).\end{equation*}
\end{definition}
\begin{notation}
We will write $\PD$ for the 2-full 2-subcategory of $\Dott$ of primary doctrines and primary homomorphisms.
\end{notation}
\begin{example}
For a given category $\CC$ with finite limits, the inclusion of $\Sub_\CC(A)$ into the poset reflection of $\CC/A$ yields a natural transformation $\Sub_\CC\to\mathbf{\Psi}_\CC$; pairing it with the identity on the base category $\CC$, this defines an elementary homomorphism.
\end{example}
\subsection{Kleisli constructions in the 2-category of indexed posets}\label{sub:kleisli}
The remaining part of the preliminaries section is devoted to showing the existence of Kleisli objects for comonads in the 2-category of indexed posets---that are essentially doctrines without the assumption of having a base category with finite products. We will use the Kleisi construction to define a new doctrine $P_{(X,\varphi)}$ from a starting doctrine $P$ in Section \ref{sect:kleisli} and to prove its universal property in Section \ref{sect:univ_prop}. Before we delve into the details, we will provide a brief overview of the relevant definitions and concepts. In the 2-category $\QDot$ of indexed posets the cells are defined as follows:
\begin{itemize}
\item a \emph{0-cell} is a functor $P:\CC\op\to\Pos$;
\item a \emph{1-cell} between $P:\CC^{\text{op}}\to\Pos$ and $R:\ct{D}^{\text{op}}\to\Pos$ is a pair $(F,\mathfrak{f})$ where $F:\CC\to\ct{D}$ is a functor and $\mathfrak{f}:P\xrightarrow{\cdot}R\circ F^{\text{op}}$ is a natural transformation;
\item a \emph{2-cell} between $(F,\mathfrak{f}),(G,\mathfrak{g}):P\to R$ is a natural transformation $\theta:F\xrightarrow{\cdot} G$ such that $\mathfrak{f}_A(\alpha)\leq R(\theta_A)(\mathfrak{g}_A(\alpha))$ for any object $A$ in $\CC$ and $\alpha\in P(A)$.
\end{itemize}
\begin{remark}
The 2-category $\Dott$ is a 2-full 2-subcategory of $\QDot$.
\end{remark}
In the following, definitions of comonads, Eilenberg--Moore and Kleisli objects in a 2-category are taken from \cite{MonCom} (see also \cite{street}). 

A \emph{comonad} in the 2-category of indexed posets is a list $(P:\CC\op\to\Pos,(K,\mathfrak{k}),\gamma, \varepsilon)$ where $P$ is a indexed poset, $(K,\mathfrak{k})$ is a 1-arrow, $\gamma$ and $\varepsilon$ are 2-arrows, and $(K,\gamma,\varepsilon)$ is a comonad in $\CC$.
In particular, the following diagrams commute:
\begin{center}
\begin{tikzcd}
K\arrow[r,"\gamma"]\arrow[d,"\gamma"]&K^2\arrow[d,"K(\gamma)"]&&K\arrow[dl,"\gamma"']\arrow[dr,"\gamma"]\arrow[d,equal]\\
K^2\arrow[r,"\gamma_K"]&K^3&K^2\arrow[r,"K(\varepsilon)"']&K&K^2\arrow[l,"\varepsilon_K"]
\end{tikzcd}.
\end{center}
Moreover, since $\gamma$ and $\varepsilon$ are 2-arrows, the following inequalities hold:
\begin{equation*}\mathfrak{k}_A(\alpha)\leq P(\gamma_A)\mathfrak{k}_{KA}\mathfrak{k}_A(\alpha);\qquad\qquad\mathfrak{k}_A(\alpha)\leq P(\varepsilon_A)(\alpha).\end{equation*}
Before studying the Kleisli objects, we briefly recall Eilenberg--Moore objects for comonads.
Define the 2-category $\Cmd(\QDot)$.
\begin{itemize}
\item a \emph{0-cell} is a comonad $(P:\CC\op\to\Pos,(K,\mathfrak{k}),\gamma, \varepsilon)$;
\item a \emph{1-cell} from the comonad $(P:\CC\op\to\Pos,(K,\mathfrak{k}),\gamma, \varepsilon)$ to $(P':{\CC'}\op\to\Pos,(K',\mathfrak{k}'),\gamma', \varepsilon')$ is a lax morphism of comonads, i.e.\ a pair $((F,\mathfrak{f}),\mathfrak{j})$ where the first entry $(F,\mathfrak{f}):P\to P'$ is a 1-arrow in $\QDot$ and the second one $\mathfrak{j}:(F,\mathfrak{f})(K,\mathfrak{k})\to(K',\mathfrak{k}')(F,\mathfrak{f})$ is a 2-arrow, i.e.\ $\mathfrak{j}:FK\xrightarrow{\cdot} K'F$ such that $\mathfrak{f}_{KA}\mathfrak{k}_A(\alpha)\leq P'(\mathfrak{j}_A)\mathfrak{k}'_{FA}\mathfrak{f}_A(\alpha)$, satisfying the coherence diagrams below;
\end{itemize}
\begin{center}
\begin{tikzcd}
FK\arrow[r,"\mathfrak{j}"]\arrow[d,"F(\gamma)"]&K'F\arrow[rd,"\gamma'_F", bend left]&&FK\arrow[r,"\mathfrak{j}"]\arrow[rd,"F(\varepsilon)",bend right]&K'F\arrow[d,"\varepsilon'_F"]\\
FK^2\arrow[r,"\mathfrak{j}_K"]&K'FK\arrow[r,"K'(\mathfrak{j})"]&K'^2F&&F
\end{tikzcd}
\end{center}
\begin{itemize}
\item a \emph{2-cell} between $((F,\mathfrak{f}),\mathfrak{j})$ and $((G,\mathfrak{g}),\mathfrak{h})$ is a 2-arrow $\eta:(F,\mathfrak{f})\to(G,\mathfrak{g})$ in $\QDot$, i.e.\ $\eta:F\xrightarrow{\cdot} G$ such that $\mathfrak{f}_A(\alpha)\leq P'(\eta_A)\mathfrak{g}_A(\alpha)$, satisfying the coherence diagram below.
\end{itemize}
\begin{center}
\begin{tikzcd}
FK\arrow[r,"\mathfrak{j}"]\arrow[d,"\eta_K"]&K'F\arrow[d,"K'(\eta)"]\\
GK\arrow[r,"\mathfrak{h}"]&K'G
\end{tikzcd}
\end{center}
\begin{definition}
A 2-category $\chi$ has \emph{Eilenberg--Moore object for comonads} if the 2-functor $\inc:\chi\to\Cmd(\chi)$, which associates to every object the identity comonad, has a right 2-adjoint
\begin{equation*}(\blank)\text{-}\mathrm{Coalg}:\Cmd(\chi)\to\chi.\end{equation*}
\end{definition}
\begin{proposition}\label{prop:eil_moo}
The 2-category $\QDot$ has Eilenberg--Moore object.\end{proposition}
The proof of this result can be found in \cite{DocModCom}. Here we just display how the right 2-adjoint $(\blank)\text{-}\mathrm{Coalg}:\Cmd(\QDot)\to\QDot$ is computed on 0-cells, since it will be useful for defining Kleisli objects: for a comonad $(P:\CC\op\to\Pos,(K,\mathfrak{k}),\gamma, \varepsilon)$, define the functor $P^K:{\CC^K}\op\to\Pos$ as follows. Let $\CC^K$ be the category of coalgebras for the comonad $(\CC,K,\gamma, \varepsilon)$ in $\Cat$: its objects are pairs $(A,c)$ where $A$ is an object of $\CC$ and $c:A\to KA$ is a $\CC$-arrow such that
\begin{center}
\begin{tikzcd}

A\arrow[r,"c"]\arrow[d,"c"]&KA\arrow[d,"K(c)"]&A\arrow[r,"c"]\arrow[dr,"\id{A}", bend right]&KA\arrow[d,"\varepsilon_A"]\\
KA\arrow[r,"\gamma_A"]&K^2A&&A

\end{tikzcd},
\end{center}

while an arrow between $(A,c)$ and $(B,d)$ is a $\CC$-arrow $f:A\to B$ such that
\begin{center}
\begin{tikzcd}

A\arrow[r,"f"]\arrow[d,"c"]&B\arrow[d,"d"]\\
KA\arrow[r,"K(f)"]&KB

\end{tikzcd}.
\end{center}

Let $P^K(A,c)$ be $\{\alpha\in P(A)\mid\alpha\leq P(c)\mathfrak{k}_A(\alpha)\}$ and $P^K(f)$ be the restriction of $P(f)$.
We now define the 2-category $\Cmd^*(\QDot):=\Cmd(\QDot\op)\op$.
\begin{itemize}
\item a \emph{0-cell} is a comonad $(P:\CC\op\to\Pos,(K,\mathfrak{k}),\gamma, \varepsilon)$;
\item a \emph{1-cell} from the comonad $(P:\CC\op\to\Pos,(K,\mathfrak{k}),\gamma, \varepsilon)$ to $(P':{\CC'}\op\to\Pos,(K',\mathfrak{k}'),\gamma', \varepsilon')$ is an oplax morphism of comonads, i.e.\ a pair $((F,\mathfrak{f}),\mathfrak{j})$ where the first entry $(F,\mathfrak{f}):P\to P'$ is a 1-arrow in $\QDot$ and the second one $\mathfrak{j}:(K',\mathfrak{k}')(F,\mathfrak{f})\to(F,\mathfrak{f})(K,\mathfrak{k})$ is a 2-arrow, i.e.\ $\mathfrak{j}:K'F\xrightarrow{\cdot} FK$ such that $\mathfrak{k}'_{FA}\mathfrak{f}_A(\alpha)\leq P'(\mathfrak{j}_A)\mathfrak{f}_{KA}\mathfrak{k}_A(\alpha)$, satisfying the coherence diagrams below;
\end{itemize}
\begin{center}
\begin{tikzcd}

K'F\arrow[r,"\mathfrak{j}"]\arrow[d,"\gamma'_F"]&FK\arrow[rd,"F(\gamma)", bend left]&&K'F\arrow[r,"\mathfrak{j}"]\arrow[rd,"\varepsilon'_F",bend right]&FK\arrow[d,"F(\varepsilon)"]\\
K'^2F\arrow[r,"K'(\mathfrak{j})"]&K'FK\arrow[r,"\mathfrak{j}_K"]&FK^2&&F

\end{tikzcd}
\end{center}
\begin{itemize}
\item a \emph{2-cell} between $((F,\mathfrak{f}),\mathfrak{j})$ and $((G,\mathfrak{g}),\mathfrak{h})$ is a 2-arrow $\eta:(F,\mathfrak{f})\to(G,\mathfrak{g})$ in $\QDot$, i.e.\ $\eta:F\xrightarrow{\cdot} G$ such that $\mathfrak{f}_A(\alpha)\leq P'(\eta_A)\mathfrak{g}_A(\alpha)$, satisfying the coherence diagram below.
\end{itemize}
\begin{center}
\begin{tikzcd}

K'F\arrow[r,"\mathfrak{j}"]\arrow[d,"K'(\eta)"]&FK\arrow[d,"\eta_K"]\\
K'G\arrow[r,"\mathfrak{h}"]&GK

\end{tikzcd}
\end{center}
\begin{definition}
A 2-category $\chi$ has \emph{Kleisli object for comonads} if the 2-functor associating to every object the identity comonad $\inc:\chi\to\Cmd^*(\chi)$, has a left 2-adjoint
\begin{equation*}(\blank)\text{-}\mathrm{coKl}:\Cmd^*(\chi)\to\chi.\end{equation*}
\end{definition}

We will devote the remaining of the section to prove Proposition \ref{prop:kleisli} below. The reader may skip the proof and proceed directly to Section \ref{sect:comonad}, coming back here to look up the definition of the left 2-adjoint on 0-cells and of the universal arrow when needed.
\begin{proposition}\label{prop:kleisli}
The 2-category $\QDot$ has Kleisli object.
\begin{proof} In order to prove the statement, we shall explicitly construct the left 2-adjoint
\begin{equation*}(\blank)\text{-}\mathrm{coKl}:\Cmd^*(\QDot)\to\QDot.\end{equation*}
We obviously begin with\newline
{\bf 0-cells:}
Fix a comonad $(P:\CC\op\to\Pos,(K,\mathfrak{k}),\gamma, \varepsilon)$, and consider the functor $P_K:{\CC_K}\op\to\Pos$. Let $\CC_K$ be the category of free coalgebras for the comonad $(\CC,K,\gamma, \varepsilon)$ in $\Cat$: it is the full subcategory of $\CC^K$ whose objects are pairs $(KA,\gamma_A)$ where $A$ is an object of $\CC$. 
Let $P_K(KA,\gamma_A)$ be $P^K(KA,\gamma_A)=\{\alpha\in P(KA)\mid\alpha\leq P(\gamma_A)\mathfrak{k}_{KA}(\alpha)\}$ and $P_K(f)$ be $ P^K(f)=P(f)$. This restriction is well defined because $P^K$ is.
\begin{remark}\label{rmk:squiggly}
The category $\CC_K$ is isomorphic to the category $\overline{\CC_K}$ whose objects are the same as $\CC$, and a $\overline{\CC_K}$-arrow $A\rightsquigarrow B$ is a $\CC$-arrow $KA\to B$; composition between $g:A\rightsquigarrow B$ and $h:B\rightsquigarrow C$ is computed as
\begin{equation*}KA\xrightarrow{\gamma_A}K^2A\xrightarrow{K(g)}KB\xrightarrow{h}C;\end{equation*}
the identity of $A$ is given by $\varepsilon_A$.

The functor ${\CC_K}\to \overline{\CC_K}$ sends $f:(KA,\gamma_A)\to(KB,\gamma_B)$ to $\varepsilon_Bf:A\rightsquigarrow B$.
The inverse $\overline{\CC_K}\to\CC_K$ sends $g:A\rightsquigarrow B$ to $K(g)\gamma_A$ (for more details see, for example, \cite[Proposition 4.1.6]{borvol2} for the proof of the dual statement concerning the Kleisli category of a monad).
\end{remark}

We now resume the proof of Proposition \ref{prop:kleisli}.\newline
{\bf 1-cells:}
Consider $((F,\mathfrak{f}),\mathfrak{j}):(P,(K,\mathfrak{k}),\gamma, \varepsilon)\to(P':{\CC'}\op\to\Pos,(K',\mathfrak{k}'),\gamma', \varepsilon')$ in $\Cmd^*(\QDot)$. The corresponding 1-cell in $\QDot$ will be $(F',\mathfrak{f}')$ from $P_K$ to $P'_{K'}$.
\begin{center}
\begin{tikzcd}

{\CC_K}^{\text{op}}\arrow[rr,"{F'}\op"] \arrow[dr,"P_K"' ,""{name=L}]&&{\CC'_{K'}}^{\text{op}}\arrow[dl,"{P'}_{K'}" ,""'{name=R}]\\
&\Pos\arrow[rightarrow,"\mathfrak{f}'","\cdot"', from=L, to=R, bend left=10]

\end{tikzcd}
\end{center}

Define $F'(KA,\gamma_A):=(K'FA,\gamma'_{FA})$, which is by definition a free $K'$-coalgebra.

Then, take a morphism between free $K$-coalgebras $f:(KA,\gamma_A)\to(KB,\gamma_B)$, and let $F'(f)$ be $K'(F(\varepsilon_Bf)\mathfrak{j}_A)\gamma'_{FA}$.
\begin{equation*}K'FA\xrightarrow{\mathfrak{j}_A}FKA\xrightarrow{F(f)}FKB\xrightarrow{F(\varepsilon_B)}FB\end{equation*}
\begin{equation*}F'(f):= \bigg(K'FA\xrightarrow{\gamma'_{FA}}K'^2FA\xrightarrow{K'(F(\varepsilon_Bf)\mathfrak{j}_A)}K'FB\bigg)\end{equation*}
This is a morphism of $K'$-coalgebras if and only if the following diagram commutes:
\begin{center}
\begin{tikzcd}

K'FA\arrow[r,"F'(f)"]\arrow[d,"\gamma'_{FA}"]&K'FB\arrow[d,"\gamma'_{FB}"]\\
K'^2FA\arrow[r,"K'F'(f)"]&K'^2FB

\end{tikzcd},
\end{center}
however
\begin{equation*}\gamma'_{FB}K'(F(\varepsilon_Bf)\mathfrak{j}_A)\gamma'_{FA}=K'^2(F(\varepsilon_Bf)\mathfrak{j}_A)\gamma'_{K'FA}\gamma'_{FA}=K'^2(F(\varepsilon_Bf)\mathfrak{j}_A)K'(\gamma'_{FA})\gamma'_{FA},\end{equation*}
using naturality of $\gamma'$ and its comultiplication property.

Consider that $\mathfrak{f}':P_K\xrightarrow{\cdot}{P'}_{K'}{F'}\op$, i.e.\ for any $(KA,\gamma_A)$ we look for
\begin{equation*}\mathfrak{f}'_{(KA,\gamma_A)}:{P_K(KA,\gamma_A)}\to {P'_{K'}(K'FA,\gamma'_{FA})},\end{equation*}
where ${P_K(KA,\gamma_A)}{\subseteq PKA}$ and ${P'_{K'}(K'FA,\gamma'_{FA})}{\subseteq P'K'FA}$.
So define $\mathfrak{f}'_{(KA,\gamma_A)}$ to be the restriction of the following composition:
\begin{equation*}PKA\xrightarrow{\mathfrak{f}_{KA}}P'FKA\xrightarrow{P'(\mathfrak{j}_A)}P'K'FA\xrightarrow{\mathfrak{k}'_{K'FA}}P'K'^2FA\xrightarrow{P'(\gamma'_{FA})}P'K'FA.\end{equation*}
To prove that the restriction is well defined, take $\alpha\in P_K(KA,\gamma_A)$, i.e.\ $\alpha\in PKA$ such that $\alpha\leq P(\gamma_A)\mathfrak{k}_{KA}(\alpha)$, we want to check that
\begin{equation*}\mathfrak{f}'_{(KA,\gamma_A)}(\alpha)\leq P'(\gamma'_{FA})\mathfrak{k}'_{K'FA}(\mathfrak{f}'_{(KA,\gamma_A)}(\alpha)).\end{equation*}
This can be computed using that $\gamma':(K',\mathfrak{k}')\to(K',\mathfrak{k}')^2$ is a 2-arrow, comultiplication property of $\gamma'$ and naturality of $\mathfrak{k}'$.

Then we check naturality of $\mathfrak{f}'$ proving that for any $f:(KA,\gamma_A)\to(KB,\gamma_B)$, the following diagram commutes:
\begin{center}
\begin{tikzcd}

P_K(KB,\gamma_B)\arrow[r,"\mathfrak{f}'_{(KB,\gamma_B)}"]\arrow[d,"P(f)"]&P'_{K'}(K'FB,\gamma'_{FB})\arrow[d,"P'(K'(F(\varepsilon_Bf)\mathfrak{j}_A)\gamma'_{FA})"]\\
P_K(KA,\gamma_A)\arrow[r,"\mathfrak{f}'_{(KA,\gamma_A)}"]&P'_{K'}(K'FA,\gamma'_{FA})

\end{tikzcd}.
\end{center}
Observe that we can decompose the diagram above as follows:
\begin{center}
\begin{tikzcd}

PKB\arrow[r,"\mathfrak{f}_{KB}"]\arrow[d,"P(f)"name=X]&P'FKB\arrow[r,"\mathfrak{k}'_{FKB}"]\arrow[d,"P'F(f)",""name=Y]&P'K'FKB\arrow[r,"P'K'(\mathfrak{j}_B)"]\arrow[d,"P'K'F(f)"name=Z]&P'K'^2FB\arrow[r,"P'(\gamma'_{FB})"]&P'K'FB\arrow[d,"P'(K'(F(\varepsilon_Bf)\mathfrak{j}_A)\gamma'_{FA})"]\\

PKA\arrow[r,"\mathfrak{f}_{KA}"]&P'FKA\arrow[r,"\mathfrak{k}'_{FKA}"]&P'K'FKA\arrow[r,"P'K'(\mathfrak{j}_A)"]&P'K'^2FA\arrow[r,"P'(\gamma'_{FA})"]&P'K'FA

\end{tikzcd}.
\end{center}

The first two squares commute because they are naturality squares of $\mathfrak{f}$ and $\mathfrak{k}$ respectively. To prove commutativity of the third square, it is enough to prove that the following square commutes:
\begin{center}
\begin{tikzcd}
K'FA\arrow[r,"\gamma'_{FA}"]\arrow[d,"K'(F(\varepsilon_Bf)\mathfrak{j}_A)\gamma'_{FA}"]&K'^2FA\arrow[r,"K'(\mathfrak{j}_A)"]&K'FKA\arrow[d,"K'F(f)"]\\
K'FB\arrow[r,"\gamma'_{FB}"]&K'^2FB\arrow[r,"K'(\mathfrak{j}_B)"]&K'FKB
\end{tikzcd}
\end{center}
and this can be shown by doing some computation using naturality of $\gamma'$, naturality of $\mathfrak{j}$, $\gamma'$ comultiplication property, $\mathfrak{j}$ coherence. In particular, we use that $K(\varepsilon_Bf)\gamma_A=K(\varepsilon_B)\gamma_Bf=f$ because of the definition of morphism between coalgebras and a property of $\varepsilon$.\newline
{\bf 2-cells:} Take a 2-cell $\eta:((F,\mathfrak{f}),\mathfrak{j})\to((G,\mathfrak{g}),\mathfrak{h})$, and look for $\eta':F'\xrightarrow{\cdot}G'$ such that
\begin{equation*}\mathfrak{f}'_{(KA,\gamma_A)}(\alpha)\leq P'_{K'}(\eta'_{(KA,\gamma_A)})(\mathfrak{g}'_{(KA,\gamma_A)}(\alpha)).\end{equation*}
Define $\eta'_{(KA,\gamma_A)}:(K'FA,\gamma'_{FA})\to(K'GA,\gamma'_{GA})$, as $\eta'_{(KA,\gamma_A)}:= K'(\eta_A)$.

Naturality diagram of $\gamma'$ applied to $\eta_A$ proves that $\eta'$ is a morphism of coalgebras.

To prove naturality of $\eta'$ we have to check that for any $f:(KA,\gamma_A)\to(KB,\gamma_B)$ the following diagram commutes:
\begin{center}
\begin{tikzcd}

(K'FA,\gamma'_{FA})\arrow[r,"K'(\eta_A)"]\arrow[d,"K'(F(\varepsilon_Bf)\mathfrak{j}_A)\gamma'_{FA}"]&(K'GA,\gamma'_{GA})\arrow[d,"K'(G(\varepsilon_Bf)\mathfrak{h}_A)\gamma'_{GA}"]\\

(K'FB,\gamma'_{FB})\arrow[r,"K'(\eta_B)"]&(K'GB,\gamma'_{GB})

\end{tikzcd}
\end{center}
but this follows from naturality of $\gamma'$, coherence diagram of $\eta$ and its naturality.

Finally, $\eta'$ is indeed a 2-arrow, i.e.\ $\mathfrak{f}'_{(KA,\gamma_A)}(\alpha)\leq P'_{K'}(\eta'_{(KA,\gamma_A)})(\mathfrak{g}'_{(KA,\gamma_A)}(\alpha))$, using that $\eta$ is a 2-arrow, its coherence, naturality of $\mathfrak{k}$ and naturality of $\gamma'$.\newline
{\bf Universal property:}
In order to prove that $(\blank)\text{-coKl}:\Cmd^*(\QDot)\to\QDot$ is indeed a left adjoint, we have to find for each comonad $(P,(K,\mathfrak{k}),\gamma, \varepsilon)$ a universal arrow
\begin{equation*}(P,(K,\mathfrak{k}),\gamma, \varepsilon)\longrightarrow(P_K,(\id{},\id{}),\id{}, \id{})\end{equation*}

i.e.\ a 1-arrow $((F_K,\mathfrak{k}'),\mathfrak{j}')$ such that, for any indexed poset $R:\ct{D}\op\to\Pos$ and any arrow
\begin{equation*}((F,\mathfrak{f}),\mathfrak{j}):(P,(K,\mathfrak{k}),\gamma, \varepsilon)\to\text{Inc}(R),\end{equation*}
there exists a unique morphism $\overline{((F,\mathfrak{f}),\mathfrak{j})}$ between the indexed posets $P_K$ and $R$ such that $((F,\mathfrak{f}),\mathfrak{j})=(\overline{((F,\mathfrak{f}),\mathfrak{j})},\id{})\circ((F_K,\mathfrak{k}'),\mathfrak{j}')$.
\begin{equation}\label{eq:uni_prop}
\begin{tikzcd}
(P,(K,\mathfrak{k}),\gamma, \varepsilon)\arrow[dr, bend right=15pt,"{((F,\mathfrak{f}),\mathfrak{j})}"']\arrow[rr, "{((F_K,\mathfrak{k}'),\mathfrak{j}')}"]&&(P_K,(\id{},\id{}),\id{}, \id{})\arrow[dl,dashed, bend left=15pt]\\
&(R,(\id{},\id{}),\id{}, \id{})
\end{tikzcd}
\end{equation}
Define $(F_K,\mathfrak{k}'):P\to P_K$, where $F_K:\CC\to\CC_K$ is the cofree functor
\begin{equation*}\bigg(f:A\to B\bigg)\longmapsto\bigg(K(f):(KA,\gamma_A)\to(KB,\gamma_B)\bigg)\end{equation*}
and the natural transformation $\mathfrak{k}':P\xrightarrow{\cdot}P_K{F_K}\op$ is computed as $\mathfrak{k}$: this is well defined since, recalling that $\gamma$ is a 2-arrow, we know that $\mathfrak{k}_A(\alpha)\leq P(\gamma_A)\mathfrak{k}_{KA}\mathfrak{k}_A(\alpha)$, i.e.\ $\mathfrak{k}_A(\alpha)\in P_KF_KA$. Naturality of $\mathfrak{k}'$ follows from naturality of $\mathfrak{k}$.

Finally, define the 2-arrow $\mathfrak{j}':(F_K,\mathfrak{k})\to(F_K,\mathfrak{k})(K,\mathfrak{k})$ to be $\mathfrak{j}'_A:= \gamma_A$. This is natural and a 2-arrow because $\gamma$ is.

Now consider $((F,\mathfrak{f}),\mathfrak{j})$; by definition of 1-cells in $\Cmd^*(\QDot)$, we know that $(F,\mathfrak{f})$ is a 1-arrow from $P$ to $R$, and $\mathfrak{j}:F\xrightarrow{\cdot}FK$ is such that $\mathfrak{f}_{A}(\alpha)\leq R(\mathfrak{j}_A)\mathfrak{f}_{KA}\mathfrak{k}_{A}(\alpha)$, and the coherence diagrams become:
\begin{center}
\begin{tikzcd}

F\arrow[r,"\mathfrak{j}"]\arrow[d,"\mathfrak{j}"]&FK\arrow[d,"F(\gamma)"]&F\arrow[r,"\mathfrak{j}"]\arrow[rd,"\id{}",bend right]&FK\arrow[d,"F(\varepsilon)"]\\
FK\arrow[r,"\mathfrak{j}_K"]&FK^2&&F

\end{tikzcd}.
\end{center}

Define the functor $F':\CC_K\to\ct{D}$ to be the one that maps $f:(KA,\gamma_A)\to(KB,\gamma_B)$ to the composition $F(\varepsilon_Bf)\mathfrak{j}_A:FA\to FB$. We check that this is indeed a functor: for any pair of composable arrows $f:KA\to KB$, $g:KB\to KC$ between free coalgebras we prove that $F'(gf)=F'(g)F'(f)$ using naturality of $\mathfrak{j}$, its coherence, definition of morphism between coalgebras and property of the counit. The identity is trivially preserved by $F'$.

To conclude the definition of $\overline{((F,\mathfrak{f}),\mathfrak{j})}$, we have to find a natural transformation $\mathfrak{f}':P_K \xrightarrow{\cdot}R{F'}\op$. Define $\mathfrak{f}'_{(KA,\gamma_A)}:P_K(KA,\gamma_A)\to RFA$ to be the restriction of $R(\mathfrak{j}_A)\mathfrak{f}_{KA}$. To show $\mathfrak{f}'$ is natural we need for any $f:KA\to KB$ between free coalgebras that the following diagram commutes:
\begin{center}
\begin{tikzcd}
P_K(KB,\gamma_B)\arrow[r,"\mathfrak{f}'_{(KB,\gamma_B)}"]\arrow[d,"P(f)"]&RF'(KB,\gamma_B)\arrow[d,"RF'(f)"]\\

P_K(KA,\gamma_A)\arrow[r,"\mathfrak{f}'_{(KA,\gamma_A)}"]&RF'(KA,\gamma_A)
\end{tikzcd}.
\end{center}

To see this, observe that it is enough to prove the commutativity of the second square of
\begin{center}
\begin{tikzcd}

PKB\arrow[r,"\mathfrak{f}_{KB}"]\arrow[d,"P(f)"]&RFKB\arrow[r,"R(\mathfrak{j}_B)"]\arrow[d,"RF(f)"]&RFB\arrow[d,"R(F(\varepsilon_Bf)\mathfrak{j}_A)"]\\

PKA\arrow[r,"\mathfrak{f}_{KA}"]&RFKA\arrow[r,"R(\mathfrak{j}_A)"]&RFA
\end{tikzcd},
\end{center}
since commutativity of the first square follows from naturality of $\mathfrak{f}$, but again it is enough to prove:
\begin{center}
\begin{tikzcd}

FA\arrow[r,"\mathfrak{j}_A"]\arrow[d,"F(\varepsilon_Bf)\mathfrak{j}_A"]&FKA\arrow[d,"F(f)"]\\
FB\arrow[r,"\mathfrak{j}_B"]&FKB

\end{tikzcd},
\end{center}

but this follows from naturality of $\mathfrak{j}$, its coherence, definition of morphism between coalgebras and counit property; so $\mathfrak{f}'$ is indeed a natural transformation.

We now prove that $((F',\mathfrak{f}'),\id{})\circ((F_K,\mathfrak{k}),\gamma)=((F,\mathfrak{f}),\mathfrak{j})$.
\begin{center}
\begin{tikzcd}

\CC\op\arrow[r,"{F_K}\op"]\arrow[dr,"P"',""{name=L},bend right]&{\CC_K}\op\arrow[r,"{F'}\op"]\arrow[d,""'{name=C},"P_K"{name=M}]&\ct{D}\op\arrow[dl, "R",""'{name=R}, bend left]\\
&\Pos

\arrow[rightarrow,"\mathfrak{k}"near end,"\cdot"', from=L, to=C, bend left=10]
\arrow[rightarrow,"\mathfrak{f}'"near start,"\cdot"', from=M, to=R, bend left=10]

\end{tikzcd}
\end{center}

The composition of the functors is indeed $F$:
\begin{multline*}\bigg(f:A\to B\bigg)\longmapsto\bigg(K(f):F_KA\to F_KB\bigg)\\
=\bigg(K(f):(KA,\gamma_A)\to(KB,\gamma_B)\bigg)\longmapsto\bigg(F(\varepsilon_BK(f))\mathfrak{j}_A:FA\to FB\bigg),\end{multline*}
but $F(\varepsilon_BK(f))\mathfrak{j}_A=F(f\varepsilon_A)\mathfrak{j}_A=F(f)$ from naturality of $\varepsilon$ and coherence of $\mathfrak{j}$.

Concerning the composition of the natural transformations, we need to check that $R(\mathfrak{j}_A)\mathfrak{f}_{KA}\mathfrak{k}_A(\alpha)=\mathfrak{f}_A(\alpha)$. The direction $(\geq)$ follows from the definition of $\mathfrak{j}$. To prove the converse, recall that coherence of $\mathfrak{j}$ implies that $F(\varepsilon_A)\mathfrak{j}_A$ is the identity, so
\begin{equation*}\mathfrak{f}_A=R(\mathfrak{j}_A)RF(\varepsilon_A)\mathfrak{f}_A=R(\mathfrak{j}_A)\mathfrak{f}_{KA}P(\varepsilon_A)\end{equation*}
 from naturality of $\mathfrak{f}$. Moreover, since $\varepsilon$ is a 2-arrow, we know that $\mathfrak{k}_A(\alpha)\leq P(\varepsilon_A)(\alpha)$, so that $R(\mathfrak{j}_A)\mathfrak{f}_{KA}\mathfrak{k}_A(\alpha)\leq R(\mathfrak{j}_A)\mathfrak{f}_{KA}P(\varepsilon_A)(\alpha)=\mathfrak{f}_A(\alpha)$, i.e.\ $(\leq)$ holds.

The composition $(\id{}\circ\gamma)_A$ is $F'(\gamma_A)=F(\varepsilon_{KA}\gamma_A)\mathfrak{j}_A=\mathfrak{j}_A$.

Finally, suppose that also $(G,\mathfrak{g})$ is such that $((G,\mathfrak{g}),\id{})\circ((F_K,\mathfrak{k}),\gamma)=((F,\mathfrak{f}),\mathfrak{j})$. Then in particular $GF_K=F$, so that $G=F'$ on objects; moreover, $\id{}\circ\gamma=\mathfrak{j}$ means that $G(\gamma_A)=\mathfrak{j}_A$. Observe also that, given a $\CC$-arrow $g$, $GK(g)=GF_K(g)=F(g)$. We claim that for any morphism $f:KA\to KB$ between free coalgebras, $G(f)=F'(f)$, i.e.\ $G(f)=F(\varepsilon_Bf)\mathfrak{j}_A$.
First of all, coherence of $\mathfrak{j}$ proves that
\begin{equation*}G(f)=F(\varepsilon_B)\mathfrak{j}_BG(f)F(\varepsilon_A)\mathfrak{j}_A.\end{equation*}
However, $\mathfrak{j}_BG(f)F(\varepsilon_A)=F(f)\mathfrak{j}_AF(\varepsilon_A)$, so we obtain $G(f)=F(\varepsilon_B)F(f)\mathfrak{j}_AF(\varepsilon_A)\mathfrak{j}_A=F(\varepsilon_Bf)\mathfrak{j}_A=F'(f)$, i.e.\ the functor $G$ is indeed the functor $F'$.

To conclude, we have to prove that $\mathfrak{g}=\mathfrak{f}'$, i.e.\ $\mathfrak{g}_{(KA,\gamma_A)}=R(\mathfrak{j}_A)\mathfrak{f}_{KA}$. We know that $\mathfrak{g}\mathfrak{k}=\mathfrak{f}$, i.e.\ $\mathfrak{g}_{(KA,\gamma_A)}\mathfrak{k}_A=\mathfrak{f}_A$, where $\mathfrak{g}_{(KA,\gamma_A)}:P_K(KA,\gamma_A)\to RFA$.

Note that
\begin{align*}\mathfrak{f}'_{(KA,\gamma_A)}&=R(\mathfrak{j}_A)\mathfrak{f}_{KA}=R(\mathfrak{j}_A)\mathfrak{g}_{(K^2A,\gamma_{KA})}\mathfrak{k}_{KA}\\
&=RG(\gamma_A)\mathfrak{g}_{(K^2A,\gamma_{KA})}\mathfrak{k}_{KA}=\mathfrak{g}_{(KA,\gamma_A)}P(\gamma_A)\mathfrak{k}_{KA},\end{align*}
because of the property of the composition of $\mathfrak{k}$ and $\mathfrak{g}$ described above and naturality of $\mathfrak{g}$. We only need to prove that the composition $P(\gamma_A)\mathfrak{k}_{KA}$ acts like the identity on $P_K(KA,\gamma_A)$. So take $\alpha\in PKA$ such that $\alpha\leq P(\gamma_A)\mathfrak{k}_{KA}(\alpha)$; we claim that $\alpha= P(\gamma_A)\mathfrak{k}_{KA}(\alpha)$. Clearly $(\leq)$ holds by definition, so we prove the converse. Recall that $\varepsilon$ is a 2-arrow, so $\mathfrak{k}_{KA}(\alpha)\leq P(\varepsilon_{KA})(\alpha)$, and apply $P(\gamma_A)$:
\begin{equation*}P(\gamma_A)\mathfrak{k}_{KA}(\alpha)\leq P(\gamma_A)P(\varepsilon_{KA})(\alpha)=\alpha\end{equation*}
so we proved that $\mathfrak{g}=\mathfrak{f}'$ and $((F_K,\mathfrak{k}),\gamma)$ is a universal arrow.
\newline
{\bf The isomorphism between the \text{Hom}-categories:}
The adjunction we proved above induces a bijection on objects of the \text{Hom}-categories below for any indexed poset $R$ and any comonad $(P,(K,\mathfrak{k}),\gamma, \varepsilon)\in\Cmd^*(\QDot)$. We need to extend it on $2$-arrows and prove it is an isomorphism of categories.

\[\begin{tikzcd}[row sep=small]
	\QDot & {\Cmd^*(\QDot)} \\
	\\
	{\Cmd^*(\QDot)[(P,(K,\mathfrak{k}),\gamma, \varepsilon),\inc(R)]} & {\QDot[P_K,R]} \\
	{\big((F,\mathfrak f),\mathfrak j\big)} & {(F',R(\mathfrak j)\mathfrak{f}_K)} \\
	\\
	{\big((G,\mathfrak g),\mathfrak h\big)} & {(G',R(\mathfrak h)\mathfrak{g}_K)}
	\arrow[""{name=0, anchor=center, inner sep=0}, "\eta", Rightarrow, from=4-1, to=6-1]
	\arrow["\cong"{description}, draw=none, from=3-1, to=3-2]
	\arrow[""{name=1, anchor=center, inner sep=0}, "\eta", Rightarrow, from=4-2, to=6-2]
	\arrow[""{name=2, anchor=center, inner sep=0}, "{(\blank)\text{-coKl}}", shift left=2, from=1-2, to=1-1]
	\arrow[""{name=3, anchor=center, inner sep=0}, "\inc", shift left=2, from=1-1, to=1-2]
	\arrow[shorten <=22pt, shorten >=22pt, maps to, from=0, to=1]
	\arrow["\dashv"{anchor=center, rotate=90}, draw=none, from=2, to=3]
\end{tikzcd}\]
where $F'\big(f:(KA,\gamma_A)\to(KB,\gamma_B)\big)=\big(F(\varepsilon_Bf)\mathfrak{j}_A:FA\to FB\big)$, and similarly for $G'$.

Take $\eta:\big((F,\mathfrak f),\mathfrak j\big)\to\big((G,\mathfrak g),\mathfrak h\big)$, i.e.\:
\begin{enumerate}
\item $\eta:F\xrightarrow{\cdot} G$ is a natural transformation;
\item $\mathfrak{f}_A(\alpha)\leq P(\eta_A)\mathfrak{g}_A(\alpha)$ for any object $A$ in $\CC$ and $\alpha\in PA$;
\item \begin{tikzcd}
	F & KF \\
	G & KG
	\arrow["{\mathfrak j}", from=1-1, to=1-2]
	\arrow["\eta"', from=1-1, to=2-1]
	\arrow["{\mathfrak h}"', from=2-1, to=2-2]
	\arrow["\eta_K", from=1-2, to=2-2]
\end{tikzcd} is commutative.
\end{enumerate}

We prove that $\eta$ in also a $2$-arrow between the correspondent indexed posets, defining $\eta_{(KA,\gamma_A)}=\eta_A:FA\to GA$. It is a natural transformation between the functors $F'$ and $G'$ if for any $f:(KA,\gamma_A)\to(KB,\gamma_B)$, we have $\eta_BF(\varepsilon_Bf)\mathfrak{j}_A=G(\varepsilon_Bf)\mathfrak{h}_A\eta_A$. However:
\begin{equation*}\eta_BF(\varepsilon_Bf)\mathfrak{j}_A=G(\varepsilon_Bf)\eta_{KA}\mathfrak{j}_A=G(\varepsilon_Bf)\mathfrak{h}_A\eta_A\end{equation*}
using 1.\ and 3. Then we need $R(\mathfrak{j}_A)\mathfrak{f}_{KA}(\alpha)\leq R(\eta_{(KA,\gamma_A)})R(\mathfrak{h}_A)\mathfrak{g}_{KA}(\alpha)$, but we know that
\begin{equation*}R(\eta_A)R(\mathfrak{h}_A)\mathfrak{g}_{KA}(\alpha)=R(\mathfrak{j}_A)R(\eta_{KA})\mathfrak{g}_{KA}(\alpha),\end{equation*}
so the inequality follows from 2.

To show this functor is full, take a natural transformation $\eta:F'\to G'$ such that $R(\mathfrak{j}_A)\mathfrak{f}_{KA}(\alpha)\leq R(\eta_{(KA,\gamma_A)})R(\mathfrak{h}_A)\mathfrak{g}_{KA}(\alpha)$ for any $\alpha\in P(KA)$ satisfying $\alpha\leq P(\gamma_A)\mathfrak{k}_{KA}(\alpha)$, and we prove that 1., 2.\ and 3.\ hold. Take any $f:A\to B$, so that $K(f):(KA,\gamma_A)\to(KB,\gamma_B)$ is a $\CC_K$ arrow; apply naturality to $K(f)$ so $\eta_{(KB,\gamma_B)}F'(Kf)=G'(Kf)\eta_{(KA,\gamma_A)}$. However, $F'(Kf)=F(\varepsilon_B)FK(f)\mathfrak{j}_A=F(\varepsilon_B)\mathfrak{j}_BF(f)=F(f)$, using naturality and coherence of $\mathfrak{j}$; similarly $G'(Kf)=G(f)$, so $\eta$ is a natural transformation from $F$ to $G$. To show 2., take any $\beta\in PB$, we want $\mathfrak{f}_B(\beta)\leq R(\eta_B)\mathfrak{g}_B(\beta)$, but by some computation we did above, we know that $\mathfrak{f}_B=R(\mathfrak{j}_B)\mathfrak{f}_{KB}\mathfrak{k}_B$---and similarly $\mathfrak{g}_B=R(\mathfrak{h}_B)\mathfrak{g}_{KB}\mathfrak{k}_B$---, so $\mathfrak{f}_B(\beta)=R(\mathfrak{j}_B)\mathfrak{f}_{KB}\mathfrak{k}_B(\beta)\leq R(\eta_B)R(\mathfrak{h}_B)\mathfrak{g}_{KB}\mathfrak{k}_B(\beta)=R(\eta_B)\mathfrak{g}_B(\beta)$, using the fact that $\mathfrak{k}_B(\beta)\in P_K(B)$ and definition of $\eta$. Finally, observe that by definition of $\gamma$, we have that $\gamma_A:(KA,\gamma_A)\to(K^2A,\gamma_{K^2A})$ is a $\CC_K$-arrow, so apply naturality of $\eta$ with respect to $\gamma_A$ to obtain $G'(\gamma_A)\eta_{(KA,\gamma_A)}=\eta_{(K^2A,\gamma_{K^2A})}F'(\gamma_A)$, i.e.\ $G(\varepsilon_{KA}\gamma_A)\mathfrak{h}_A\eta_A=\eta_{KA}F(\varepsilon_{KA}\gamma_A)\mathfrak{j}_A$, hence $\mathfrak{h}_A\eta_A=\eta_{KA}\mathfrak{j}_A$, so that also 3.\ holds.
Again, faithfulness follows by definition, and it is essentially surjective because of the properties of adjunction. It is clear that the quasi-inverse is actually an inverse.\end{proof}
\end{proposition}

\section{A comonad on the indexed poset $P$}\label{sect:comonad}
We now have all the ingredients to start our discussion concerning how to translate into the language of doctrines, seen as a generalization of the doctrine of well-formed formulae, the process of adding a constant of some fixed sort and adding a sentence to a theory.
We give an intuition of these two translations.

We start with the generalization of the concept of adding a constant. Constants in a first-order language are $0$-ary terms, terms are arrows in the category of contexts, the empty list of variables is the terminal object of the category of contexts, thus a \emph{constant} for a doctrine is just an arrow with the terminal object in the base category of the doctrine as domain object; the codomain object of the constant is the \emph{sort} of that constant. So, given an object $X$ in the base category of a doctrine $P$, \emph{adding a constant of sort $X$} means, roughly speaking, defining from $P$ a new doctrine $P_X$ endowed with a constant of sort $X$ for $P_X$, in a universal way. This construction and its universal property are formally stated in \cref{coroll:const}.

We then generalize the concept of adding an axiom. In first-order logic, we can add any sentence---that is, a formula with no free variables---to a theory, so that the sentence becomes a new axiom in the extended theory. Since the empty list of variables is the terminal object in the category of contexts, a \emph{sentence} in a doctrine is just an element in the fiber over the terminal object in the base category, and an \emph{axiom} in a doctrine is just a sentence that is the true constant in the fiber over the terminal object. So, given a sentence $\varphi$ in a doctrine $P$, \emph{adding the axiom $\varphi$} means, roughly speaking, defining from $P$ a new doctrine $P_\varphi$ in which $\varphi$ is the top element of the fiber over the terminal object---hence ``$\varphi$ is true''---, in a universal way. Observe that to do so, it is reasonable that the doctrine has at least the true constant, so we will work with primary doctrines---and not just with any doctrine, as we could do with the construction that adds a constant instead. This construction and its universal property are formally stated in \cref{coroll:ax}.

As hinted in the introduction, although these two may seem like separate processes, we show that they can be computed simultaneously.
\begin{notation} From now on, $P:\CC\op\to\Pos$ is a fixed primary doctrine, unless otherwise specified.
\end{notation}

Fix an object $X$ in the base category \CC, and an element $\varphi\in P(X)$.
\begin{center}
\begin{tikzcd}

\CC^{\text{op}}\arrow[rr,"(X\times\blank)^{\text{op}}"] \arrow[dr,"P"' ,""{name=L}]&&\CC^{\text{op}}\arrow[dl,"P" ,""'{name=R}]\\
&\Pos\arrow[rightarrow,"\mathfrak{f}","\cdot"', from=L, to=R, bend left=10]

\end{tikzcd}
\end{center}

Consider the product functor $X\times\blank:\CC\to\CC$ sending $(A\xrightarrow{f}B)$ to $(X\times A\xrightarrow{\id{X}\times f}X\times B)$, and define each component of the natural transformation $\mathfrak{f}:P\xrightarrow{\cdot} P\circ(X\times\blank)^{\text{op}}$ as follows:
\begin{equation*}\mathfrak{f}_A:P(A)\to P(X\times A)\text{, }\alpha\mapsto P(\pr1)(\varphi)\land P(\pr2)(\alpha)\end{equation*}
where $\pr1$ and $\pr2$ are the projection from $X\times A$ to $X$ and $A$ respectively.

Note that $\mathfrak{f}$ is monotone and it is a natural transformation.
Then, $(X\times\blank,\mathfrak{f})$ is a 1-cell between $P$ and itself in the category $\QDot$.

We now prove that the 1-arrow $(X\times\blank,\mathfrak{f})$ is part of a 2-comonad on $P$; to do this, we have to find two 2-arrows $\varepsilon:(X\times\blank,\mathfrak{f})\xrightarrow{\cdot}\id{P}$ and $\gamma:(X\times\blank,\mathfrak{f})\xrightarrow{\cdot}(X\times\blank,\mathfrak{f})^2$ satisfying the proper diagrams. We adapt the comonad on the functor $X\times \blank$ (also known as the \emph{reader comonad}) to indexed posets.

Define $\varepsilon_A:X\times A\to A$ to be the second projection $\varepsilon_A\coloneqq \pr2$, which is clearly natural, and is indeed a 2-arrow since $\mathfrak{f}_A(\alpha)=P(\pr1)(\varphi)\land P(\pr2)(\alpha)\leq P(\pr2)(\alpha)= P(\varepsilon_A)(\id{A}(\alpha))$ for any $\alpha\in P(A)$. Then, define $\gamma_A\coloneqq\Delta_X\times \id{A}:X\times A\to X\times X\times A$, which is again natural; it is a 2-arrow if and only if $\mathfrak{f}_A(\alpha)\leq P(\Delta\times \id{})(\mathfrak{f}_{X\times A}\mathfrak{f}_A(\alpha))$, however
\begin{align*}P(\Delta\times \id{})(\mathfrak{f}_{X\times A}\mathfrak{f}_A(\alpha))&=P(\ple{\pr1,\pr1,\pr2})\big(P(\pr1)\varphi\land P(\ple{\pr2,\pr3})(P(\pr1)\varphi\land P(\pr2)\alpha)\big)\\
&=P(\ple{\pr1,\pr1,\pr2})\big(P(\pr1)\varphi\land P(\pr2)\varphi\land P(\pr3)\alpha\big)\\
&=P(\pr1)\varphi\land P(\pr2)\alpha=\mathfrak{f}_A(\alpha).\end{align*}
Finally, since the following diagrams commute:
\begin{center}
\begin{tikzcd}

&X\times\blank\arrow[d,"\gamma"]\arrow[dl,equal]\arrow[dr,equal]&&X\times\blank\arrow[d,"\gamma"]\arrow[r,"\gamma"]&(X\times\blank)^2\arrow[d,"\id{X}\times\gamma"]\\
X\times\blank&(X\times\blank)^2\arrow[l,"\varepsilon_{X\times\blank}"]\arrow[r,"\id{X}\times\varepsilon"']&X\times\blank&(X\times\blank)^2\arrow[r,"\gamma_{X\times\blank}"]&(X\times\blank)^3

\end{tikzcd}
\end{center}
we proved the following
\begin{proposition}\label{prop:com_idx}
With the notation defined above, $(P,(X\times \blank,\mathfrak f), \gamma,\varepsilon)$ is a comonad in $\QDot$.
\end{proposition}
\begin{remark} We are interested in finding a distributive law between two different comonads, both of the form seen in Proposition \ref{prop:com_idx} on the same primary doctrine $P$ seen as an indexed poset: for two objects $X,Y$ and two elements $\varphi\in P(X), \psi\in P(Y)$, the first comonad relies on the 1-cell $(X\times\blank, \mathfrak{f})$, where $\mathfrak f=P(\pr1)(\varphi)\land P(\pr2)(\blank)$, while the second one relies on the 1-cell $(Y\times\blank,\mathfrak{g})$, where $\mathfrak g=P(\pr1)(\psi)\land P(\pr2)(\blank)$.

Dualizing distributive laws for monads in \cite{distrlaw}, recall that in general, for two given comonads $(P,(K,\mathfrak{f}),\gamma,\epsilon)$ and $(P,(C,\mathfrak{g}),\lambda,\delta)$ in $\QDot$ on the same indexed poset $P:\CC\op\to\Pos$, a distributive law between two comonads is a 2-cell $\ell:(K,\mathfrak{f})\circ(C,\mathfrak{g})\to(C,\mathfrak{g})\circ(K,\mathfrak{f})$ such that $((K,\mathfrak f),\ell)$ is a lax morphism of comonads and $((C,\mathfrak g),\ell)$ is an oplax morphism of comonads. In details, in $\QDot$, we ask for $\ell:KC\xrightarrow{\cdot} CK$ to be a natural transformation, such that $\mathfrak{f}_{CA}\mathfrak{g}_A(\alpha)\leq P(\ell_A)\mathfrak{g}_{KA}\mathfrak{f}_A(\alpha)$ for any object $A$ in $\CC$ and $\alpha\in P(A)$ and such that the following diagram commute:
\[\begin{tikzcd}
	KC && CK & KC && CK \\
	{KC^2} & CKC & {C^2K} & {K^2C} & KCK & {CK^2} \\
	KC && K & KC && C \\
	& CK &&& CK
	\arrow["\ell", from=1-1, to=1-3]
	\arrow["{K(\lambda)}"', from=1-1, to=2-1]
	\arrow["{\lambda_K}", from=1-3, to=2-3]
	\arrow["\ell", from=1-4, to=1-6]
	\arrow["{C(\gamma)}", from=1-6, to=2-6]
	\arrow["{K(\delta)}", from=3-1, to=3-3]
	\arrow["\ell"', from=3-1, to=4-2]
	\arrow["{\gamma_C}"', from=1-4, to=2-4]
	\arrow["{\ell_C}"', from=2-1, to=2-2]
	\arrow["{C(\ell)}"', from=2-2, to=2-3]
	\arrow["{K(\ell)}"', from=2-4, to=2-5]
	\arrow["{\delta_K}"', from=4-2, to=3-3]
	\arrow["{\epsilon_C}", from=3-4, to=3-6]
	\arrow["\ell"', from=3-4, to=4-5]
	\arrow["{C(\epsilon)}"', from=4-5, to=3-6]
	\arrow["{\ell_K}"', from=2-5, to=2-6]
\end{tikzcd}.\]
In our case, define for each object $A$ in \CC,
\begin{equation*}\ell_A\coloneqq\ple{ \pr2,\pr1,\pr3}:X\times Y\times A\to Y\times X\times A.\end{equation*}
This is trivially a natural transformation. Moreover, recall that the comultiplication is given by $\Delta\times\id{}$, and the counit by the projection on the second component. With these definitions, the diagrams clearly commute:
\[\begin{tikzcd}
	{X\times Y\times A} && {Y\times X\times A} \\
	{X\times Y\times Y\times A} & {Y\times X\times Y\times A} & {Y\times Y\times X\times A}
	\arrow["{\ple{\pr2,\pr1,\pr3}}", from=1-1, to=1-3]
	\arrow["{\ple{\pr1,\pr2,\pr2,\pr3}}"', from=1-1, to=2-1]
	\arrow["{\ple{\pr2,\pr1,\pr3,\pr4}}"', from=2-1, to=2-2]
	\arrow["{\ple{\pr1,\pr1,\pr2,\pr3}}", from=1-3, to=2-3]
	\arrow["{\ple{\pr1,\pr3,\pr2,\pr4}}"', from=2-2, to=2-3]
\end{tikzcd}\]
\[\begin{tikzcd}
	{X\times Y\times A} && {Y\times X\times A} \\
	{X\times X\times Y\times A} & {X\times Y\times X\times A} & {Y\times X\times X\times A}
	\arrow["{\ple{\pr2,\pr1,\pr3}}", from=1-1, to=1-3]
	\arrow["{\ple{\pr1,\pr1,\pr2,\pr3}}"', from=1-1, to=2-1]
	\arrow["{\ple{\pr1,\pr3,\pr2,\pr4}}"', from=2-1, to=2-2]
	\arrow["{\ple{\pr1,\pr2,\pr2,\pr3}}", from=1-3, to=2-3]
	\arrow["{\ple{\pr2,\pr1,\pr3,\pr4}}"', from=2-2, to=2-3]
\end{tikzcd}\]
\[\begin{tikzcd}
	{X\times Y\times A} & {X\times A} & {X\times Y\times A} & {Y\times A} \\
	{Y\times X\times A} && {Y\times X\times A}
	\arrow["{\ple{\pr2,\pr1,\pr3}}"', from=1-1, to=2-1]
	\arrow["{\ple{\pr1,\pr3}}", from=1-1, to=1-2]
	\arrow["{\ple{\pr2,\pr3}}"', from=2-1, to=1-2]
	\arrow["{\ple{\pr2,\pr3}}", from=1-3, to=1-4]
	\arrow["{\ple{\pr2,\pr1,\pr3}}"', from=1-3, to=2-3]
	\arrow["{\ple{\pr1,\pr3}}"', from=2-3, to=1-4]
\end{tikzcd}.\]
Now we need to prove that $\ell$ is indeed a 2-cell, i.e.\
\[\begin{tikzcd}
	{P(A)} & {P(X\times A)} & {P(Y\times X\times A)} \\
	{P(Y\times A)} && {P(X\times Y\times A)}
	\arrow["{\mathfrak{g}_A}"', from=1-1, to=2-1]
	\arrow["{\mathfrak{f}_A}", from=1-1, to=1-2]
	\arrow["{\mathfrak{g}_{KA}}", from=1-2, to=1-3]
	\arrow["{P(\ell_A)}", from=1-3, to=2-3]
	\arrow["{\mathfrak{f}_{CA}}"', from=2-1, to=2-3]
	\arrow["\leq"{description}, draw=none, from=1-3, to=2-1]
\end{tikzcd}.\]
Compute:
\begin{align*}\mathfrak{f}_{CA}\mathfrak{g}_A(\alpha)&=\mathfrak{f}_{CA}\big(P(\pr1)\psi\land P(\pr2)\alpha\big)=P(\pr1)\varphi\land P(\pr2)\psi\land P(\pr3)\alpha\\
&=P(\ell_A)\big(P(\pr1)\psi\land P(\pr2)\varphi\land P(\pr3)\alpha\big)\\
&=P(\ell_A)\mathfrak{g}_{KA}\big(P(\pr1)\varphi\land P(\pr2)\alpha\big)= P(\ell_A)\mathfrak{g}_{KA}\mathfrak{f}_A(\alpha).\end{align*}
In this particular case, $\ell$ is actually an isomorphism.

By looking at the commutative triangles, we observe that $\ell$ is unique:
\begin{align*}
\pr1\ell_A=\pr1\ple{\pr1,\pr3}\ell_A=\pr1\ple{\pr2,\pr3}=\pr2;\\
\pr2\ell_A=\pr1\ple{\pr2,\pr3}\ell_A=\pr1\ple{\pr1,\pr3}=\pr1;\\
\pr3\ell_A=\pr2\ple{\pr2,\pr3}\ell_A=\pr2\ple{\pr2,\pr3}=\pr3.
\end{align*}
We conclude that the distributive law $\ell$ induces a composite comonad, having the 1-cell computed as $(X\times\blank,\mathfrak{f})\circ(Y\times\blank,\mathfrak {g})=(X\times Y\times\blank,\mathfrak{f}\circ\mathfrak{g})$, where
\begin{equation*}(\mathfrak{f}\circ\mathfrak{g})_A=\mathfrak{f}_{Y\times A}\mathfrak{g}_A:P(A)\to P(X\times Y\times A), \quad\alpha\mapsto P(\pr1)\varphi\land P(\pr2)\psi\land P(\pr3)\alpha.\end{equation*}
Moreover, the composite comonad induced by the distributive law is again of the form seen in Proposition \ref{prop:com_idx}, defined with respect to the object $X\times Y$ and the element $P(\pr1)\varphi\land P(\pr2)\psi$ in $P(X\times Y)$.\end{remark}

\section{The Kleisli construction for the comonad $(X\times\blank,\mathfrak f)$}\label{sect:kleisli}

We now study the Kleisli construction of the comonad $(P,(X\times \blank,\mathfrak f), \gamma,\varepsilon)$, applying the results shown in Section \ref{sub:kleisli} to this particular case. To do so, we first have to compute the Eilenberg--Moore construction. The Eilenberg--Moore category $\CC^{X\times\blank}$ has as objects pairs $(A,a)$, where $a:A\to X\times A$ is an arrow in $\CC$ such that the following diagram commutes:
\begin{center}
\begin{tikzcd}

&A\arrow[r,"a"]\arrow[d,"a"]&X\times A\arrow[d,"\gamma_A=\Delta\times \id{}"]\\
A\arrow[ur,equal]&X\times A\arrow[l,"\varepsilon_A=\pr2"]\arrow[r,"\id{}\times a"']&X\times X\times A 

\end{tikzcd}
\end{center}
so that the second component of $a=\ple{ a_1,a_2}$ must be the identity $a_2=\id{A}$, while the first one can be any map $a_1:A\to X$. Moreover, an arrow $f:(A,a)\to(B,b)$ in $\CC^{X\times\blank}$ is an arrow $f:A\to B$ in $\CC$ such that the diagram commutes:
\begin{center}
\begin{tikzcd}

A\arrow[r,"f"]\arrow[d,"{\ple{ a_1,\id{A}}}"']&B\arrow[d,"{\ple{ b_1,\id{B}}}"]\\
X\times A\arrow[r,"\id{}\times f"]&X\times B

\end{tikzcd}
\end{center}
i.e.\ we ask for the \CC-arrow $f$ to satisfy $b_1 f=a_1$. From now on, we will write $\CC^X$ instead of $\CC^{X\times\blank}$. By looking at the description of objects and arrows of the category $\CC^X$, it is easy to observe $\CC^X$ is isomorphic to the slice category $\CC/X$: they are both categories of coalgebras of the reader comonad $X\times\blank$.

As shown in Proposition \ref{prop:eil_moo}, the induced indexed poset $P^{(X,\varphi)}:{\CC^X}\op\to\Pos$ is defined as follows:
\begin{center}For
\begin{tikzcd}
(B,b)\\
(A,a)\arrow[u,"f"]
\end{tikzcd}
the reindexing is
\begin{tikzcd}
\{\beta\in P(B)\mid\beta\leq P(b)(\mathfrak{f}_B(\beta))\}\arrow[d,"P^{(X,\varphi)}(f)=P(f)_{\mid P^{(X,\varphi)}(B,b)}"]\\
\{\alpha\in P(A)\mid\alpha\leq P(a)(\mathfrak{f}_A(\alpha))\}
\end{tikzcd}
\end{center}
with the order of the subsets given by $P(B)$ and $P(A)$ respectively. Since by definition $\mathfrak{f}_A(\alpha)=P(\pr1)(\varphi)\land P(\pr2)(\alpha)$, we can write
\begin{align*}P^{(X,\varphi)}(A,a)&=\{\alpha\in P(A)\mid\alpha\leq P(\ple{ a_1,\id{}})(P(\pr1)(\varphi)\land P(\pr2)(\alpha))\}\\
&=\{\alpha\in P(A)\mid\alpha\leq P(a_1)(\varphi)\}=P(A)_{\downarrow P(a_1)(\varphi)}.\end{align*}

\begin{example}
To have an intuition of the Eilenberg-Moore construction, we compute the indexed poset $P^{(X,\varphi)}$ when $P$ is the the syntactic doctrine $\LT^\mathcal{L}_{\mathcal{T}}$ in \cref{ex:doctr}\ref{i:fbf}, $X$ is a finite list of variables $\vec{x}$ and $\varphi$ is a formula $\varphi(\vec{x})$ in the context $\vec{x}$. Objects of the base category are lists of terms $\vec{t}(\vec{z})$ substitutable to $\vec{x}$, for some finite list of variables $\vec{z}$. An arrow in the base category between $\vec{t}(\vec{z})$ and $\vec{s}(\vec{y})$ is a list term $\vec{u}(\vec{z})$ substitutable to $\vec{y}$ such that $\vec{s}([\vec{u}(\vec{z})/\vec{y}])=\vec{t}(\vec{z})$. For a list $\vec{t}(\vec{z})$, the corresponding fiber is $\{\alpha(\vec{z})\in \LT^\mathcal{L}_{\mathcal{T}}(\vec{z})\mid \alpha(\vec{z})\vdash_{\mathcal{T}}\varphi([\vec{t}(\vec{z})/\vec{x}])\}$. Finally, the reindexing along $\vec{u}(\vec{z}):\vec{t}(\vec{z})\to \vec{s}(\vec{y})$ maps $\beta(\vec{y})\in \LT^\mathcal{L}_{\mathcal{T}}(\vec{y})$ such that $\beta(\vec{y})\vdash_{\mathcal{T}}\varphi([\vec{s}(\vec{y})/\vec{x}])$ to $\beta([\vec{u}(\vec{z})/\vec{y}])\in \LT^\mathcal{L}_{\mathcal{T}}(\vec{z})$.
\end{example}

Now, consider the Kleisli category $\CC_{X\times\blank}$, i.e.\ the full subcategory of $\CC^X$ whose objects are the co-free algebras. From now on, we will write $\CC_X$ instead on $\CC_{X\times\blank}$. Observe that an arrow $f=\ple{ f_1,f_2}:X\times A\to X\times B$ has to satisfy $\pr1 f=\pr1$, so $f_1$ must be the first projection $\pr1$ and the map $f$ is uniquely determined by its second component $f_2:X\times A\to B$. For this reason, from now on we will use the equivalent description of $\CC_X$, that has as objects the same as \CC, and as map $g:A\rightsquigarrow B$ is a \CC-arrow $g:X\times A\to B$---see Remark \ref{rmk:squiggly} for more details; moreover, the composition between two arrows $g:A\rightsquigarrow B$ and $h:B\rightsquigarrow C$ is the arrow $h\ple{ \pr1,g}:A\rightsquigarrow C$. A new indexed poset $P_{(X,\varphi)}$ is trivially induced on the Kleisli category by simply taking the restriction of $P^{(X,\varphi)}$ on ${\CC_X}\op$, so that $P_{(X,\varphi)}:{\CC_X}\op\to\Pos$ is defined as follows:
\begin{center}for
\begin{tikzcd}
B\\
A\arrow[u,"g", squiggly]
\end{tikzcd}
the reindexing is
\begin{tikzcd}
P(X\times B)_{\downarrow P(\pr1)(\varphi)}\arrow[d,"{P_{(X,\varphi)}(g)=P(\ple{ \pr1,g })_{\mid P_{(X,\varphi)}(B)}}"]\\
P(X\times A)_{\downarrow P(\pr1)(\varphi)}
\end{tikzcd}.
\end{center}

The Kleisli construction also comes with a 1-arrow between indexed posets $(F_X,\mathfrak{f}):P\to P_{(X,\varphi)}$, where $F_X$ is the co-free functor that sends $f:A\to B$ to its precomposition with the second projection $f\pr2:A\rightsquigarrow B$, and where each component of $\mathfrak f$ is defined as $\mathfrak{f}_A=P(\pr1)(\varphi)\land P(\pr2)(\blank):P(A)\to P_{(X,\varphi)}(A)$.

\begin{remark}\label{r:axiom-quotient}
It is easy to see that for a fixed pair $X\in\CC, \varphi\in P(X)$, fiber $P(X\times A)_{\downarrow P(\pr1)(\varphi)}$ of the doctrine $P_{(X,\varphi)}$ over an object $A$ is isomorphic to the quotient of $P(X\times A)$ where $[\alpha]\leq[\beta]$ if and only if $P(\pr1)(\varphi)\land \alpha\leq\beta$ in $P(X\times A)$. One direction of the isomorphism takes any $\alpha\in P(X\times A)$, $\alpha\leq P(\pr1)(\varphi)$, and sends it to its equivalence class; the other direction takes any equivalence class $[\alpha]$ for some $\alpha\in P(X\times A)$ and sends it to the conjunction $P(\pr1)(\varphi)\land \alpha$ of any representative with $P(\pr1)(\varphi)$.
\end{remark}

Roughly speaking, if we think of $P$ as the syntactic doctrine $\LT^\mathcal{L}_{\mathcal{T}}$ in \cref{ex:doctr}\ref{i:fbf}, we can interpret the intexed poset $P_{(X,\varphi)}$ as follows. Suppose we want to add a new constant $c$ to a first-order language $\mathcal{L}$ and the formula $\varphi([c/x])$ to the theory $\mathcal{T}$---thus we have a new axiom. We obtain a new language $\mathcal{L}'=\mathcal{L}\cup\{c\}$ and a new theory in the extended language $\mathcal{T}'=\mathcal{T}\cup\{\varphi([c/x])\}$: the indexed poset $P_{(X,\varphi)}$ should be thought as the syntactic doctrine $\LT^\mathcal{L'}_{\mathcal{T'}}$. Since adding a constant to a language clearly doesn't change the lists of finite variables, the objects of the base category of $P_{(X,\varphi)}$ are the same as the ones in the base category of $P$. A term $u(\vec{z})$ in the language $\mathcal{L}'$ is a term $u'(\,(x);\vec{z}\,)$ in the language $\mathcal{L}$ with the constant $c$ substituted to the variable $x$, thus a term $u(\vec{z}): \vec{z}\to (y)$ in the extended language is just an old term $u'(\,(x);\vec{z}\,): (\,(x);\vec{z}\,)\to (y)$ evaluated in $[c/x]$. This is the intuition behind the fact that an arrow $A\rightsquigarrow B$ in the base category $\CC_X$ of $P_{(X,\varphi)}$ is just an old arrow $X\times A \to B$. In particular, we are able to pick the constant $c$ added to the language as the $0$-ary term $c: ()\to (x)$ in the language $\mathcal{L'}$, that is interpreted as $\pr1: \tmn\rightsquigarrow X$.
Similarly, a formula $\alpha(\vec{z})$ in the language $\mathcal{L'}$ is a formula $\alpha'((x);\vec{z})$ in the language $\mathcal{L}$ with the constant $c$ substituted to the variable $x$. This is why an element in $P_{(X,\varphi)}(A)$ is in particular an element in $P(X\times A)$. To conclude, let $\alpha(\vec{z}), \beta(\vec{z})$ be two formulae in the language $\mathcal{L'}$ and let $\alpha'((x);\vec{z}),\beta'((x);\vec{z})$ be formulae in the language $\mathcal{L}$ such that $\alpha(\vec{z})=\alpha'([c/x];\vec{z})$ and $\beta(\vec{z})=\beta'([c/x];\vec{z})$. Then
\begin{align*}\alpha(\vec{z})\vdash_{\mathcal{T}'}\beta(\vec{z}) &\text{\quad if and only if\quad}\varphi([c/x])\land \alpha(\vec{z})\vdash_{\mathcal{T}}\beta(\vec{z})\\
 &\text{\quad if and only if\quad}\varphi([c/x]) \land\alpha'([c/x];\vec{z})\vdash_{\mathcal{T}}\beta'([c/x];\vec{z})\\
 &\text{\quad if and only if\quad} \varphi(x)\land \alpha'((x);\vec{z})\vdash_{\mathcal{T}}\beta'((x);\vec{z}).\end{align*}
This is the intuition behind the presentation of the fibers of $P_{(X,\varphi)}$ described in \cref{r:axiom-quotient} above.
These ideas will be formally stated and proved in \cref{ex:const,ex:ax}, but for now they are convenient intuitions to keep in mind while reading the rest of the paper.

\section{The doctrine $P_{(X,\varphi)}$ and its inherited properties}\label{sect:pres_prop}
We now want to study if some properties of $P$ can be translated to $P_{(X,\varphi)}$, and when so, if they are preserved by the 1-arrow $(F_X,\mathfrak{f})$.
First of all, we check that $\CC_X$ has finite products, preserved by $F_X$, so that both $P$ and $P_{(X,\varphi)}$ are doctrines and $(F_X,\mathfrak f)$ is a 1-arrow in $\Dott$. Then we study a few other properties of $\CC_X$ that can be inherited from \CC; then we will take a look to various properties of $P$.
\begin{center}
\begin{tikzcd}
\CC^{\text{op}}\arrow[rr,"{F_X}^{\text{op}}"] \arrow[dr,"P"' ,""{name=L}]&&{\CC_X}^{\text{op}}\arrow[dl,"P_{(X,\varphi)}" ,""'{name=R}]\\
&\Pos\arrow[rightarrow,"\mathfrak{f}","\cdot"', from=L, to=R, bend left=10]
\end{tikzcd}
\end{center}
We begin by collecting some elementary results regarding the category $\CC_X$. The key part is that $\CC_X$ has products, so that $P_{(X,\varphi)}$ is a doctrine. The part of the statement about products follows from the dual of Proposition 2.2 in \cite{kleisli}. We then provide the proof of the other properties, as we could not find precise references.
\begin{proposition}
Let $\CC$ be a category with finite products, and $\CC_X$ be the Kleisli category of the comonad $(\CC, X\times\blank, \Delta_X\times \id{},\pr2)$. Then the category $\CC_X$ has finite products, and the co-free functor $F_X:\CC\to\CC_X$ preserves finite products.
Moreover:
\begin{enumerate}
\item If the category $\CC$ is closed, then $\CC_X$ is closed and the co-free functor $F_X:\CC\to\CC_X$ preserves the exponential.
\item Suppose that $\CC$ has initial object $I$. The category $\CC_X$ has initial object and the co-free functor $F_X:\CC\to\CC_X$ preserves the initial object if and only if the endofunctor $X\times\blank:\CC\to\CC$ preserves the initial object.
\item Suppose that $\CC$ has binary coproducts. The category $\CC_X$ has binary coproducts and the co-free functor $F_X:\CC\to\CC_X$ preserves binary coproducts if and only if the endofunctor $X\times\blank:\CC\to\CC$ preserves binary coproducts.
\end{enumerate}
\begin{proof}\begin{enumerate}
\item Suppose that for any object $Y$ there is a natural bijection
\begin{equation*}\Hom_\CC(\blank\times Y,\bblank)\cong\Hom_\CC(\blank,(\bblank)^Y).\end{equation*}
We sum up the naturality in the two components with the following diagrams:
\begin{center}
\begin{tikzcd}[row sep=0pt]
A\times Y\arrow[r,"f\times \id{}"]&B\times Y\arrow[r,"\widetilde{h}"]&Z &&B\times Y\arrow[r,"k"]&Z\arrow[r,"g"]&S\\
{}\arrow[rr, no head]&&{}&&{}\arrow[rr,no head]&&{}\\
A\arrow[r,"f"]&B\arrow[r,"h"]&Z^Y &&B\arrow[r,"\widehat{k}"]&Z^Y\arrow[r,"g^Y"]&S^Y\\
\end{tikzcd}
\end{center}
Consider the functor $\blank\times Y:\CC_X\to\CC_X$, that maps $f:A\rightsquigarrow B$, i.e.\ $f:X\times A\to B$, to the arrow $f\times \id{}:A\times Y\rightsquigarrow B\times Y$, i.e.\ $f\times \id{}:X\times A\times Y\to B\times Y$. Such functor is a left adjoint, since for each object $A$, there exists an object $A^Y$---which we will prove to be the exponential in \CC---and a $\CC_X$-arrow $\varepsilon_A:A^Y\times Y\rightsquigarrow A$, i.e.\ $\varepsilon_A:X\times A^Y\times Y\to A$ such that, for any object $B$ and arrow $f:B\times Y\rightsquigarrow A$, there exists a unique $\widehat{f}:B\rightsquigarrow A^Y$---which we will prove to be the same hatted arrow in \CC---such that
\begin{center}
\begin{tikzcd}
B\times Y\arrow[r,"\widehat{f}\times \id{}",squiggly]\arrow[rr,"f",squiggly,bend right]&A^Y\times Y\arrow[r,"\varepsilon_A",squiggly]&A
\end{tikzcd}
\end{center}
Define $\varepsilon_A\coloneqq \widetilde{\pr2}:X\times A^Y\times Y\to A$, the \CC-map corresponding to $\pr2:X\times A^Y\to A^Y$, so that $\varepsilon_A$ is indeed a $\CC_X$-map $A^Y\times Y\rightsquigarrow A$. We only have to check that the composition of $\CC_X$-arrows above equals to $f$, i.e.\ in \CC
\begin{equation*}\widetilde{\pr2}\circ\ple{ \pr1,\widehat{f}\times \id{}}=\widetilde{\pr2}\circ(\ple{ \pr1,\widehat{f}}\times \id{})=\widetilde{\pr2\circ\ple{ \pr1,\widehat{f}}}=\widetilde{\widehat{f}}=f\end{equation*}
At last, to prove the uniqueness of $\widehat{f}$, suppose $f'$ such that $\widetilde{\pr2}\circ\ple{ q_1,f'\times \id{}}=f$, but the left-hand side is equal to $\widetilde{f'}$, so $\widehat{\widetilde{f'}}=\widehat{f}$, i.e.\ $f'=\widehat{f}$.

To conclude, take a \CC-arrow $g:A\times Y\to B$, and its corresponding map $\widehat{g}:A\to B^Y$, we want to prove that $\widehat{F_X(g)}=F_X(\widehat{g})$.
\begin{equation*}F_X(g):X\times A\times Y\xrightarrow{\ple{ \pr2,\pr3}}A\times Y\xrightarrow{g}B\end{equation*}
\begin{equation*}F_X(\widehat{g}):X\times A\xrightarrow{\pr2}A\xrightarrow{\widehat{g}}B^Y\end{equation*}
By naturality we have:
\begin{center}
\begin{tikzcd}[row sep=0pt]

X\times A\times Y\arrow[r,"\pr2\times \id{}"]&A\times Y\arrow[r,"g"]&B\\
{}\arrow[rr, no head]&&{}\\
X\times A\arrow[r,"\pr2"]&B\arrow[r,"\widehat{g}"]&B^Y\\
\end{tikzcd}
\end{center}
i.e.\ $\widetilde{\widehat{g}\pr2}=g(\pr2\times \id{})=g\ple{ \pr2,\pr3}$, so that $\widehat{g}\pr2=\widehat{g\ple{ \pr2,\pr3}}$.

\item We show that $F_XI=I$ is initial in $\CC_X$: consider any object $A$, we look for a unique arrow $I\rightsquigarrow A$, i.e.\ a unique arrow $X\times I\to A$, but $X\times I=I$ by assumption.
Conversely, suppose that $\CC_X$ has initial object, preserved by $F_X$. Since $U_X:\CC_X\to\CC$, that acts $\big(g:C\rightsquigarrow D\big) \mapsto \big(\ple{\pr1,g}:X\times C\to X\times D\big)$, is a left adjoint, it preserves all colimits, and in particular $X\times I=I$, as claimed.

\item At first we suppose that the functor $X\times\blank$ preserves binary coproducts. Consider any pair of objects $A$, $B$, their coproduct diagram in \CC, its image in $\CC_X$ through $F_X$ and a pair of arrows $\alpha:X\times A\to V$ and $\beta:X\times B\to V$:
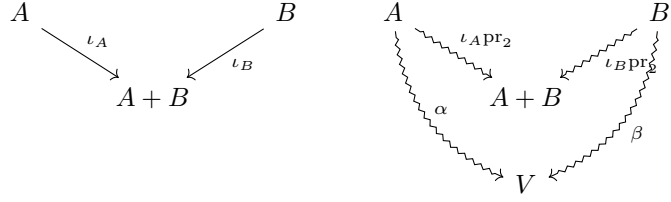
\begin{center}
\begin{tikzcd}
A\arrow[rd,"\iota_A"]&&B\arrow[ld,"\iota_B"]&A\arrow[rd,"\iota_A\pr2",squiggly]\arrow[rdd,"\alpha",squiggly, bend right]&&B\arrow[ld,"\iota_B\pr2",squiggly]\arrow[ldd,"\beta",squiggly, bend left]\\
&A+ B&&&A+ B\\
&&&&V
\end{tikzcd}
\end{center}
Define the map $A+ B\rightsquigarrow V$ to be the composition $\binom{\alpha}{\beta}\psi:X\times(A+ B)\xrightarrow[\sim]{\psi}(X\times A)+(X\times B)\xrightarrow{\binom{\alpha}{\beta}}V$, where $\psi$ is the inverse of the canonical arrow below:
\begin{center}
\begin{tikzcd}
X\times A\arrow[rd,"\iota_{X\times A}"]&&X\times B\arrow[dl,"\iota_{X\times B}"]\\
&(X\times A)+(X\times B)\arrow[d,"{\ple{\binom{\pr1}{\pr1},\binom{\iota_A \pr2}{\iota_B\pr2}}=\binom{\ple{ \pr1,\iota_A\pr2}}{\ple{ \pr1,\iota_B\pr2}}} "]\\
&X\times(A+ B)\arrow[u,bend left,"\psi"]\arrow[dl,"\pr1"]\arrow[dr,"\pr2"]\\
X&&(A+ B)
\end{tikzcd}
\end{center}
In particular,
\begin{equation*}\binom{\iota_{X\times A}}{\iota_{X\times B}}=\id{(X\times A)+(X\times B)}=\psi\binom{\ple{ \pr1,\iota_A\pr2}}{\ple{ \pr1,\iota_B\pr2}}=\binom{\psi\ple{ \pr1,\iota_A\pr2}}{\psi\ple{ \pr1,\iota_B\pr2}}.\end{equation*}
So the composition $X\times A\xrightarrow{\ple{ \pr1,\iota_A\pr2}} X\times(A+ B)\xrightarrow{\binom{\alpha}{\beta}\psi}V$ is equal to $\binom{\alpha}{\beta}\psi\ple{ \pr1,\iota_A\pr2}=\binom{\alpha}{\beta}\iota_{X\times A}=\alpha$; similarly, for $\beta$. Hence, $\binom{\alpha}{\beta}\psi$ makes the diagram commute, and it is clearly unique, so $\CC_X$ has coproducts, preserved by $F_X$ by construction.

Conversely, suppose that $\CC_X$ has coproducts, preserved by $F_X$, our claim is that in \CC, $X\times\blank$ distribute over $+$. So in $\CC_X$ take $A,B$ and their coproduct
\begin{center}
\begin{tikzcd}
A\arrow[r,squiggly,"\iota_A\pr2"]&(A+ B)&B.\arrow[l,squiggly,"\iota_B\pr2"]
\end{tikzcd}
\end{center}
Recall that $U_X:\CC_X\to\CC$, that maps $g:C\rightsquigarrow D$ to $\ple{\pr1,g}:X\times C\to X\times D$ is a left adjoint, so it preserves all colimits, and in particular $X\times(A+ B)=(X\times A)+(X\times B)$, as claimed.\qedhere\end{enumerate}\end{proof}\end{proposition}

\begin{remark}
Is it important to observe that the hypothesis about $F_X$ preserving coproducts is necessary for the equivalence in 3. of the proposition above. Indeed, suppose $\CC$ to be a bounded lattice $(R,\land,\lor)$, and fix $x\in R$; so $\CC_X=\overline{R}\coloneqq(|R|,\sqsubset)$ where $a\sqsubset a'$ if and only if $x\land a\leq a'$. The poset $\overline{R}$ has coproducts: indeed, $a\overline{\lor}b\coloneqq(x\land a)\lor(x\land b)$. Clearly $a\sqsubset a\overline{\lor}b$ and $b\sqsubset a\overline{\lor}b$; moreover, take $a,b\sqsubset y$, i.e.\ $x\land a\leq y$ and $x\land b\leq y$, then $a\overline{\lor}b\sqsubset y$ if and only if
\begin{equation*}x\land({{(x\land a)}\lor{(x\land b)}})\leq y,\end{equation*}
which holds, so that $\overline{R}$ has indeed coproducts. However, if $x\land\blank$ does not distribute over $\lor$, coproducts are not preserved---$x\land(a\lor b)\neq(x\land a)\lor(x\land b)$.
\end{remark}

We now study the structural properties of the fibers of $P$ that are inherited by $P_{(X,\varphi)}$ and preserved by the morphism $(F_X,\mathfrak f)$.

\subsubsection*{Finite meets}
\begin{proposition}\label{prop:fin_meet}
Let $P:\CC\op\to\Pos$ be a primary doctrine and $P_{(X,\varphi)}$ be the Kleisli object of the comonad $(P,(X\times\blank,\mathfrak f),\gamma,\varepsilon)$ defined by the pair $X\in\CC$ and $\varphi\in P(X)$. Then $P_{(X,\varphi)}$ is a primary doctrine, and $(F_X,\mathfrak f)$ is a primary homomorphism.
\begin{proof}
Recall that, by assumption, for any object $A$ of \CC, the poset $P(A)$ has finite meets. We want to check that $P_{(X,\varphi)}(A)= P(X\times A)_{\downarrow P(\pr1)(\varphi)}$ has finite meets too: for any two elements $\alpha,\beta\in P_{(X,\varphi)}(A)$, their meet is $\alpha\land\beta$ computed in $P(X\times A)$. Naturality of the operation is trivial. Moreover, the poset $P_{(X,\varphi)}(A)$ has a top element, which is $\ter_A\coloneqq P(\pr1)(\varphi)$, and $\ter$ is again natural.

Take $\mathfrak{f}_A:P(A)\to P_{(X,\varphi)}(A)= P(X\times A)_{\downarrow P(\pr1)(\varphi)}$. For any $\alpha,\beta\in P(A)$, one has $\mathfrak{f}_A(\alpha\land\beta)=P(\pr1)(\varphi)\land P(\pr2)(\alpha\land\beta)=P(\pr1)(\varphi)\land P(\pr2)(\alpha)\land P(\pr2)(\beta)=\mathfrak{f}_A(\alpha)\land\mathfrak{f}_A(\beta)$. Moreover, $\mathfrak{f}_A(\top_A)=P(\pr1)(\varphi)\land P(\pr2)(\top_A)=P(\pr1)(\varphi)\land \top_{X\times A}=P(\pr1)(\varphi)=\ter_A$.\end{proof}
\end{proposition}

\subsubsection*{Elementarity}
\begin{proposition}
Let $P:\CC\op\to\Pos$ be a primary doctrine and $P_{(X,\varphi)}$ be the Kleisli object of the comonad $(P,(X\times\blank,\mathfrak f),\gamma,\varepsilon)$ defined by the pair $X\in\CC$ and $\varphi\in P(X)$. If $P$ is an elementary doctrine, then $P_{(X,\varphi)}$ is an elementary doctrine, and $(F_X,\mathfrak f)$ is an elementary homomorphism.
\begin{proof}
We already proved that, since $P$ is a primary doctrine, $P_{(X,\varphi)}$ is a primary doctrine too. So now take an object $A$ and define
\begin{equation*}\delta^{X,\varphi}_A:=P(\pr1)(\varphi)\land P(\ple{\pr2,\pr3})(\delta_A)\in P_{(X,\varphi)}(A\times A)=P(X\times A\times A)_{\downarrow P(\pr1)(\varphi)},\end{equation*}
where $\delta_A\in P(A\times A)$ is the fibered equality of $P$ on $A$.
\begin{enumerate}
\item We show that in $P_{(X,\varphi)}(A)=P(X\times A)_{\downarrow P(\pr1)(\varphi)}$ we have $P(\pr1)(\varphi)\leq P(\ple{\pr1,\pr2,\pr2})(\delta^{X,\varphi}_A)$, so we compute:
\begin{align*}P(\ple{\pr1,\pr2,\pr2})(\delta^{X,\varphi}_A)&=P(\pr1)(\varphi)\land P(\ple{\pr2,\pr2})(\delta_A)=P(\pr1)(\varphi)\land P(\pr2)P(\Delta_A)(\delta_A)\\
&=P(\pr1)(\varphi)\land P(\pr2)(\top_A)=P(\pr1)(\varphi)\end{align*}
using Definition \ref{def:elem}.1.
\item For any $\alpha\in P_{(X,\varphi)}(A)$, we show that in $P_{(X,\varphi)}(A\times A)=P(X\times A\times A)_{\downarrow P(\pr1)(\varphi)}$ we have $P(\ple{\pr1,\pr2})(\alpha)\land P(\pr1)(\varphi)\land P(\ple{\pr2,\pr3})(\delta_A)\leq P(\ple{\pr1,\pr3})(\alpha)$.
To see this, observe that by property 2.\ with respect to the object $X\times A$ in Definition \ref{def:elem}, in the fiber $P(X\times A\times X\times A)$ we have
\begin{equation*}P(\ple{\pr1,\pr2})(\alpha)\land \delta_{X\times A}\leq P(\ple{\pr3,\pr4})(\alpha).\end{equation*}
Then apply $P(\ple{\pr1,\pr2,\pr1,\pr3})$ to both sides of the inequality to get that in $P(X\times A\times A)$
\begin{equation*}P(\ple{\pr1,\pr2})(\alpha)\land P(\ple{\pr1,\pr2,\pr1,\pr3})(\delta_{X\times A})\leq P(\ple{\pr1,\pr3})(\alpha).\end{equation*}
However, using properties 3.\ and 1.\ in Definition \ref{def:elem} we compute
\begin{align*}P(\ple{\pr1,\pr2,\pr1,\pr3})(\delta_{X\times A})&\geq P(\ple{\pr1,\pr2,\pr1,\pr3})\big(P(\ple{\pr1,\pr3})(\delta_{X})\land P(\ple{\pr2,\pr4})(\delta_{A})\big)\\&=P(\ple{\pr1,\pr1})(\delta_{X})\land P(\ple{\pr2,\pr3})(\delta_{A})\\&=P(\pr1)P(\Delta_X)(\delta_{X})\land P(\ple{\pr2,\pr3})(\delta_{A})\\&=P(\ple{\pr2,\pr3})(\delta_A).\end{align*}
In conclusion we get:
\begin{align*}&P(\ple{\pr1,\pr2})(\alpha)\land P(\pr1)(\varphi)\land P(\ple{\pr2,\pr3})(\delta_A)\\&\leq P(\ple{\pr1,\pr2})(\alpha)\land P(\ple{\pr2,\pr3})(\delta_A)\\&\leq P(\ple{\pr1,\pr2})(\alpha)\land P(\ple{\pr1,\pr2,\pr1,\pr3})(\delta_{X\times A})\\&\leq P(\ple{\pr1,\pr3})(\alpha)\end{align*}
as claimed.
\item For any pair of objects $A,B$, we show that in $P_{(X,\varphi)}(A\times B\times A\times B)=P(X\times A\times B\times A\times B)_{\downarrow P(\pr1)(\varphi)}$ we have $P(\ple{\pr1,\pr2,\pr4})(\delta^{X,\varphi}_A)\land P(\ple{\pr1,\pr3,\pr5})(\delta^{X,\varphi}_B)\leq \delta^{X,\varphi}_{A\times B}$, i.e.\
\begin{multline*}P(\pr1)(\varphi)\land P(\ple{\pr2,\pr4})(\delta_A)\land P(\pr1)(\varphi)\land P(\ple{\pr3,\pr5})(\delta_B)\\
\leq P(\pr1)(\varphi)\land P(\ple{\pr2,\pr3,\pr4,\pr5})(\delta_{A\times B}),\end{multline*}
but this holds by applying $P(\ple{\pr2,\pr3,\pr4,\pr5})$ to both sides of the inequality 3.\ in Definition \ref{def:elem}.
\end{enumerate}
This proves that $P_{(X,\varphi)}$ is an elementary doctrine. To conclude, we show that $(F_X,\mathfrak f)$ is an elementary morphism, but this follows from Proposition \ref{prop:fin_meet} and from the definition of the fibered equality on $A$ as $\delta^{X,\varphi}_A=\mathfrak{f}_{A\times A}(\delta_A)$. 
\end{proof}
\end{proposition}

\subsubsection*{Existential quantifier}
\begin{proposition}
Let $P:\CC\op\to\Pos$ be a primary doctrine and $P_{(X,\varphi)}$ be the Kleisli object of the comonad $(P,(X\times\blank,\mathfrak f),\gamma,\varepsilon)$ defined by the pair $X\in\CC$ and $\varphi\in P(X)$. If $P$ is an existential doctrine, then $P_{(X,\varphi)}$ is an existential doctrine, and $(F_X,\mathfrak f)$ is an existential homomorphism.
\begin{proof}
We already proved that, since $P$ is primary, $P_{(X,\varphi)}$ is primary as well. So now take objects $B,C$ and consider $P(\ple{ \pr1,\pr2}):P(X\times C)_{\downarrow P(\pr1)(\varphi)}\to P(X\times C\times B)_{\downarrow P(\pr1)(\varphi)}$, and define its left adjoint ${\exists_{X,\varphi}}^B_C$ to be the restriction of $\exists^B_{X\times C}$. Such restriction is well defined: take $\beta\leq P(\pr1)(\varphi)$ in $P(X\times C\times B)$, one has $\exists^B_{X\times C}(\beta)\leq P(\pr1)(\varphi)$ if and only if $\beta\leq P(\ple{ \pr1,\pr2})P(\pr1)(\varphi)=P(\pr1)(\varphi)$, which is true by assumption. Beck-Chevalley condition and Frobenius reciprocity for ${\exists_{X,\varphi}}^B_C$ come easy from the same properties of $\exists^B_{X\times C}$.

Consider the following diagram:
\begin{center}
\begin{tikzcd}
P(C\times B)\arrow[r,"\mathfrak{f}_{C\times B}"]\arrow[d,"\exists^B_C"{name=b}, shift left=7pt]&P(X\times C\times B)_{\downarrow P(\pr1)(\varphi)}\arrow[d,"{\exists_{X,\varphi}}^B_C=\exists^B_{X\times C}"{name=c}, shift left=7pt]\\
P(C)\arrow[r,"\mathfrak{f}_C"]\arrow[u,"P(\pr1)"{name=a}, shift left=7pt]&P(X\times C)_{\downarrow P(\pr1)(\varphi)}\arrow[u,"{P(\ple{ \pr1,\pr2})}"{name=e}, shift left=7pt]
\arrow[phantom, from=a, to=b, "\vdash"]
\arrow[phantom, from=c, to=e, "\vdash"]
\end{tikzcd}.
\end{center}
We want to prove that the square with arrows pointing down and right is commutative. To do this, recall that the following diagram is commutative because of the Beck-Chevalley condition:
\begin{center}
\begin{tikzcd}
P(C\times B)\arrow[d,"{P(\pr2\times \id{B})}"]\arrow[r,"\exists^B_C"]&P(C)\arrow[d,"{P(\pr2)}"]\\
P(X\times C\times B)\arrow[r,"\exists^B_{X\times C}"]&P(X\times C)
\end{tikzcd}.
\end{center}
So now take $\beta\in P(C\times B)$:
\begin{align*}\mathfrak{f}_C\exists^B_C(\beta)&=P(\pr1)(\varphi)\land P(\pr2)(\exists^B_C(\beta))=P(\pr1)(\varphi)\land\exists^B_{X\times C}{P(\pr2\times \id{})}(\beta)\\
&=\exists^B_{X\times C}(P(\ple{ \pr2,\pr3})(\beta)\land {P(\ple{ \pr1,\pr2})P(\pr1)}(\varphi))=\exists^B_{X\times C}\mathfrak{f}_{X\times C}(\beta),\end{align*}
which proves our claim.\end{proof}\end{proposition}
\subsubsection*{Universal quantifier}
\begin{proposition}
Let $P:\CC\op\to\Pos$ be a primary doctrine and $P_{(X,\varphi)}$ be the Kleisli object of the comonad $(P,(X\times\blank,\mathfrak f),\gamma,\varepsilon)$ defined by the pair $X\in\CC$ and $\varphi\in P(X)$. If $P$ is a universal doctrine, then $P_{(X,\varphi)}$ is a universal doctrine, and $(F_X,\mathfrak f)$ is a universal homomorphism.

Additionally, if the universal quantifier of $P$ satisfies the Frobenius reciprocity, then also the universal quantifier of $P_{(X,\varphi)}$ does.
\begin{proof}
Take a pair of objects $B,C$ and consider
\begin{equation*}P(\ple{ \pr1,\pr2}):P(X\times C)_{\downarrow P(\pr1)(\varphi)}\to P(X\times C\times B)_{\downarrow P(\pr1)(\varphi)}.\end{equation*}
Define its right adjoint ${\forall_{X,\varphi}}^B_C(\blank)\coloneqq\forall^B_{X\times C}(\blank)\land P(\pr1)(\varphi)$. To check that this yields indeed an adjunction, take $\gamma\leq P(\pr1)(\varphi)$ in $P(X\times C)$ and $\beta\leq P(\pr1)(\varphi)$ in $P(X\times C\times B)$, we want to prove that $P(\ple{ \pr1,\pr2})(\gamma)\leq\beta$ if and only if $\gamma\leq{\forall_{X,\varphi}}^B_C(\beta)=\forall^B_{X\times C}(\beta)\land P(\pr1)(\varphi)$. First of all, suppose $P(\ple{ \pr1,\pr2})(\gamma)\leq\beta$, then it follows from the adjunction $P(\ple{ \pr1,\pr2})\dashv\forall^B_{X\times C}(\beta)$ that $\gamma\leq\forall^B_{X\times C}(\beta)$; combining this with the assumption on $\gamma$, the inequality $\gamma\leq{\forall_{X,\varphi}}^B_C(\beta)$ holds. Conversely, suppose $\gamma\leq\forall^B_{X\times C}(\beta)\land P(\pr1)(\varphi)\leq\forall^B_{X\times C}(\beta)$, then the claim holds again because of the adjunction.

We now want to prove the Beck-Chevalley condition for ${\forall_{X,\varphi}}$. Take an arrow $f_2:C'\rightsquigarrow C$, i.e.\ $f_2:X\times C'\to C$, and write $f\coloneqq\ple{ \pr1,f_2}:X\times C'\to X\times C$. Recall the Beck-Chevalley condition for the universal quantifier of $P$:
\begin{center}
\begin{tikzcd}
P(X\times C\times B)\arrow[d,"{P(f\times \id{B})}"]\arrow[r,"\forall^B_{X\times C}"]&P(X\times C)\arrow[d,"{P(f)}"]\\
P(X\times C'\times B)\arrow[r,"\forall^B_{X\times C'}"]&P(X\times C')
\end{tikzcd}
\end{center}
and use it to prove the condition for ${\forall_{X,\varphi}}$:
\begin{center}
\begin{tikzcd}
P(X\times C\times B)_{\downarrow P(\pr1)(\varphi)}\arrow[d,"{P(f\times \id{B})}"]\arrow[r,"{\forall_{X,\varphi}}^B_C"]&P(X\times C)_{\downarrow P(\pr1)(\varphi)}\arrow[d,"{P(f)}"]\\
P(X\times C'\times B)_{\downarrow P(\pr1)(\varphi)}\arrow[r,"{\forall_{X,\varphi}}^B_{C'}"]&P(X\times C')_{\downarrow P(\pr1)(\varphi)}
\end{tikzcd}.
\end{center}
So take $\beta\in P(X\times C\times B)$:
\begin{align*}P(f){\forall_{X,\varphi}}^B_C(\beta)&=P(f)(\forall^B_{X\times C}(\beta)\land P(\pr1)(\varphi))=P(f)\forall^B_{X\times C}(\beta)\land P(\pr1)(\varphi)\\
&=\forall^B_{X\times C'}P(f\times \id{})(\beta)\land P(\pr1)(\varphi)={\forall_{X,\varphi}}^B_{C'}P(f\times \id{})(\beta).\end{align*}
If we ask in addition that the doctrine $P$ satisfies Frobenius reciprocity for the adjunction $P(\ple{ \pr1,\pr2})\dashv\forall^B_{X\times C}(\beta)$, then also the doctrine $P_{(X,\varphi)}$ satisfies Frobenius for the adjunction $P(\ple{ \pr1,\pr2})\dashv{\forall_{X,\varphi}}^B_{C}$: for any $\gamma\leq P(\pr1)(\varphi)$ and $\beta\leq P(\pr1)(\varphi)$,
\begin{align*}P(\ple{ \pr1,\pr2})(\gamma\land{\forall_{X,\varphi}}^B_{C}(\beta))&=P(\ple{ \pr1,\pr2})(\gamma\land{\forall}^B_{X\times C}(\beta)\land {P(\pr1)(\varphi)})\\
&=P(\ple{ \pr1,\pr2})(\gamma)\land\beta,\end{align*}
using $\gamma\land P(\pr1)\varphi=\gamma$ and Frobenius reciprocity.

To conclude, consider the following diagram:
\begin{center}
\begin{tikzcd}
P(C\times B)\arrow[r,"\mathfrak{f}_{C\times B}"]\arrow[d,"\forall^B_C"{name=b}, shift left=7pt]&P(X\times C\times B)_{\downarrow P(\pr1)(\varphi)}\arrow[d,"{\forall_{X,\varphi}}^B_C"{name=c}, shift left=7pt]\\
P(C)\arrow[r,"\mathfrak{f}_C"]\arrow[u,"P(\pr1)"{name=a}, shift left=7pt]&P(X\times C)_{\downarrow P(s_1)(\varphi)}\arrow[u,"{P(\ple{ \pr1,\pr2})}"{name=e}, shift left=7pt]
\arrow[phantom, from=a, to=b, "\dashv"]
\arrow[phantom, from=c, to=e, "\dashv"]
\end{tikzcd}
\end{center}
We want to prove that the square with arrows pointing down and right is commutative. To do this, recall that the following diagram is commutative because of the Beck-Chevalley condition:
\begin{center}
\begin{tikzcd}
P(C\times B)\arrow[d,"{P(\pr2\times \id{B})}"]\arrow[r,"\forall^B_C"]&P(C)\arrow[d,"{P(\pr2)}"]\\
P(X\times C\times B)\arrow[r,"\forall^B_{X\times C}"]&P(X\times C)
\end{tikzcd}
\end{center}
So now take $\beta\in P(C\times B)$:
\begin{align*}{\forall_{X,\varphi}}^B_C\mathfrak{f}_{C\times B}(\beta)&={\forall_{X,\varphi}}^B_C(P(\pr1)(\varphi)\land P(\ple{ \pr2,\pr3})(\beta))\\
&=\forall^B_{X\times C}(P(\pr1)(\varphi)\land P(\ple{ \pr2,\pr3})(\beta))\land P(\pr1)(\varphi).\end{align*}
On the other hand
\begin{equation*}\mathfrak{f}_C\forall^B_C(\beta)=P(\pr1)(\varphi)\land P(\pr2)(\forall^B_C(\beta)).\end{equation*}
To prove ${\forall_{X,\varphi}}^B_C\mathfrak{f}_{C\times B}(\beta)\leq\mathfrak{f}_C\forall^B_C(\beta)$, note that if holds if and only if
\begin{equation*}{\forall_{X,\varphi}}^B_C\mathfrak{f}_{C\times B}(\beta)\leq P(\pr2)\forall^B_C(\beta)\end{equation*}
but $P(\pr2)\forall^B_C=\forall^B_{X\times C}P(\pr2\times \id{})$ and moreover $P(\pr1)(\varphi)\land P({\ple{ \pr2,\pr3}})(\beta)=P(\pr1)(\varphi)\land P({\pr2\times \id{}})(\beta)\leq P(\pr2\times \id{})(\beta)$, so applying $\forall^B_{X\times C}$ to both sides of the last inequality the claim follows.

Conversely, we have $\mathfrak{f}_C\forall^B_C(\beta)\leq {\forall_{X,\varphi}}^B_C\mathfrak{f}_{C\times B}(\beta)$, if and only if $P(\ple{\pr1,\pr2})\mathfrak{f}_C\forall^B_C(\beta)\leq \mathfrak{f}_{C\times B}(\beta)$ but $P(\ple{\pr1,\pr2})=P_{(X,\varphi)}F_X(\pr1)$, so equivalently $\mathfrak{f}_{C\times B}P(\pr1)\forall^B_C(\beta)\leq \mathfrak{f}_{C\times B}(\beta)$. This proves the claim, by applying $\mathfrak{f}_{C\times B}$ to $P(\pr1)\forall^B_C(\beta)\leq\beta$---which is the counit of the adjunction.\end{proof}\end{proposition}
\subsubsection*{Implication}
\begin{proposition}
Let $P:\CC\op\to\Pos$ be a primary doctrine and $P_{(X,\varphi)}$ be the Kleisli object of the comonad $(P,(X\times\blank,\mathfrak f),\gamma,\varepsilon)$ defined by the pair $X\in\CC$ and $\varphi\in P(X)$. If $P$ is an implicational doctrine, then $P_{(X,\varphi)}$ is an implicational doctrine, and $(F_X,\mathfrak f)$ is an implicational homomorphism.
\begin{proof}
Since we already know that $P_{(X,\varphi)}$ is primary, we check that $P_{(X,\varphi)}(A)= P(X\times A)_{\downarrow P(\pr1)(\varphi)}$ is cartesian closed too: for any $\beta,\gamma\in P_{(X,\varphi)}(A)$, define $\beta\Rightarrow\gamma\coloneqq(\beta\to\gamma)\land P(\pr1)(\varphi)$. 

This is indeed a natural transformation: take $f_2:A\rightsquigarrow B$, i.e.\ $f_2:X\times A\to B$, and write for convenience $f\coloneqq\ple{ \pr1,f_2}:X\times A\to X\times B$.
\begin{center}
\begin{tikzcd}
{P(X\times B)_{\downarrow P(\pr1)(\varphi)}}\op\times P(X\times B)_{\downarrow P(\pr1)(\varphi)}\arrow[d,"P(f)\times P(f)"]\arrow[r,"\Rightarrow"]&P(X\times B)_{\downarrow P(\pr1)(\varphi)}\arrow[d,"{P(f)}"]\\
{P(X\times A)_{\downarrow P(\pr1)(\varphi)}}\op\times P(X\times A)_{\downarrow P(\pr1)(\varphi)}\arrow[r,"\Rightarrow"]&P(X\times A)_{\downarrow P(\pr1)(\varphi)}
\end{tikzcd}
\end{center}
So, take a pair $\alpha,\alpha'\in P(X\times B)_{\downarrow P(\pr1)(\varphi)}$: on the one hand it is sent to $P(f)(\alpha\Rightarrow\alpha')=P(f)((\alpha\to\alpha')\land P(\pr1)(\varphi))=(P(f)(\alpha)\to P(f)(\alpha'))\land P(\pr1)(\varphi)$; on the other hand to $P(f)(\alpha)\Rightarrow P(f)(\alpha')=(P(f)(\alpha)\to P(f)(\alpha'))\land P(\pr1)(\varphi)$, so that $\Rightarrow$ is indeed a natural transformation.

Now, to check that $P_{(X,\varphi)}(A)$ endowed with this operation is cartesian closed, take three elements $\alpha,\beta,\gamma\in P(X\times A)_{\downarrow P(\pr1)(\varphi)}$, and we prove that $\alpha\land\beta\leq\gamma$ if and only if $\alpha\leq\beta\Rightarrow\gamma$. So, suppose $\alpha\land\beta\leq\gamma$, then from $\alpha\leq\beta\to\gamma$ combined with the assumption on $\alpha$ we obtain $\alpha\leq(\beta\to\gamma)\land P(\pr1)(\gamma)$. Conversely, from $\alpha\leq(\beta\to\gamma)\land P(\pr1)(\gamma)\leq\beta\to\gamma$, it follows that $\alpha\land\beta\leq\gamma$.

Take $\mathfrak{f}_A:P(A)\to P_{(X,\varphi)}(A)=P(X\times A)_{\downarrow P(\pr1)(\varphi)}$. For any $\alpha,\beta\in P(A)$, one has on the one side $\mathfrak{f}_A(\alpha\to\beta)=P(\pr1)(\varphi)\land P(\pr2)(\alpha\to\beta)$, and on the other hand $\mathfrak{f}_A(\alpha)\Rightarrow\mathfrak{f}_A(\beta)=(P(\pr1)(\varphi)\land P(\pr2)(\alpha))\Rightarrow(P(\pr1)(\varphi)\land P(\pr2)(\beta))$. So now we prove that in any cartesian closed poset,
\begin{equation*}((x\land a)\to(x\land b))\land x=x\land(a\to b)\end{equation*}

First of all, $x\land(a\to b)\leq((x\land a)\to(x\land b))\land x$ if and only if $x\land(a\to b)\leq(x\land a)\to(x\land b)$ if and only if $x\land{(a\to b)}\land x\land {a}\leq x\land b$.

Conversely, $((x\land a)\to(x\land b))\land x\leq x\land(a\to b)$ if and only if $((x\land a)\to(x\land b))\land x\leq a\to b$ if and only if $((x\land a)\to(x\land b))\land x\land a\leq b$, but $((x\land a)\to(x\land b))\land x\land a\leq x\land b\leq b$.\end{proof}\end{proposition}

\subsubsection*{Finite joins}
\begin{proposition}
Let $P:\CC\op\to\Pos$ be a primary doctrine and $P_{(X,\varphi)}$ be the Kleisli object of the comonad $(P,(X\times\blank,\mathfrak f),\gamma,\varepsilon)$ defined by the pair $X\in\CC$ and $\varphi\in P(X)$. If $P$ is bounded, then $P_{(X,\varphi)}$ is bounded, and $(F_X,\mathfrak f)$ preserves the bottom element.
\begin{proof} The poset $P_{(X,\varphi)}(A)$ has bottom element, which is $\ini_A\coloneqq \bot_{X\times A}$, and $\ini$ is natural.

Take $\mathfrak{f}_A:P(A)\to P_{(X,\varphi)}(A)= P(X\times A)_{\downarrow P(\pr1)(\varphi)}$. Compute $\mathfrak{f}_A(\bot_A)=P(\pr1)(\varphi)\land {P(\pr2)(\bot_A)}=\bot_{X\times A}=\ini_A$, so the bottom element is preserved.\end{proof}
\end{proposition}
\begin{proposition}\label{prop:joins}
Let $P:\CC\op\to\Pos$ be a primary doctrine and $P_{(X,\varphi)}$ be the Kleisli object of the comonad $(P,(X\times\blank,\mathfrak f),\gamma,\varepsilon)$ defined by the pair $X\in\CC$ and $\varphi\in P(X)$. If $P$ has binary joins, then $P_{(X,\varphi)}$ has binary joins. If each fiber of $P$ is a distributive lattice, then $(F_X,\mathfrak f)$ preserves binary joins.
\begin{proof} To check that $P_{(X,\varphi)}(A)= P(X\times A)_{\downarrow P(\pr1)(\varphi)}$ has binary joins, take any two elements $\alpha,\beta\in P_{(X,\varphi)}(A)$, and compute their join $\alpha\lor\beta$ in $P(X\times A)$. Naturality of the operation is trivial.

Now, for any two elements $\alpha,\beta\in P(A)$, one has $\mathfrak{f}_A(\alpha\lor\beta)=P(\pr1)(\varphi)\land P(\pr2)(\alpha\lor\beta)=P(\pr1)(\varphi)\land(P(\pr2)(\alpha)\lor P(\pr2)(\beta))$

On the other hand, $\mathfrak{f}_A(\alpha)\lor\mathfrak{f}_A(\beta)=(P(\pr1)(\varphi)\land P(\pr2)(\alpha))\lor(P(\pr1)(\varphi)\land P(\pr2)(\beta))$; in general this is not equal to $\mathfrak{f}_A(\alpha\lor\beta)$, computed above. 
However, the equality holds if we ask for $P(\pr1)(\varphi)\land(\blank)$ to preserve joins, e.g.\ whenever the lattice is distributive.\end{proof}\end{proposition}
\begin{corollary}
Let $P:\CC\op\to\Pos$ be a primary doctrine and $P_{(X,\varphi)}$ be the Kleisli object of the comonad $(P,(X\times\blank,\mathfrak f),\gamma,\varepsilon)$ defined by the pair $X\in\CC$ and $\varphi\in P(X)$. If $P$ is a Heyting doctrine, then $P_{(X,\varphi)}$ is a Heyting doctrine, and $(F_X,\mathfrak f)$ is a Heyting homomorphism.
\begin{proof}Finite meets, finite joins and implication are preserved by the construction.\end{proof}
\end{corollary}

\subsubsection*{Booleanness}
\begin{proposition}\label{prop:bool}
Let $P:\CC\op\to\Pos$ be a primary doctrine and $P_{(X,\varphi)}$ be the Kleisli object of the comonad $(P,(X\times\blank,\mathfrak f),\gamma,\varepsilon)$ defined by the pair $X\in\CC$ and $\varphi\in P(X)$. If $P$ is a Boolean doctrine, then $P_{(X,\varphi)}$ is a Boolean doctrine, and $(F_X,\mathfrak f)$ is a Boolean homomorphism.
\begin{proof} We want to check that $P_{(X,\varphi)}(A)= P(X\times A)_{\downarrow P(\pr1)(\varphi)}$ is a boolean algebra: for any element $\alpha\in P_{(X,\varphi)}(A)$, define $\rceil\alpha\coloneqq\alpha\Rightarrow\bot=(\alpha\to\bot)\land P(\pr1)(\varphi)=\lnot\alpha\land P(\pr1)(\varphi)$. Since we already know that $P_{(X,\varphi)}(A)$ is a Heyting algebra, we only have to prove that $\rceil\rceil\alpha=\alpha$:
\begin{multline*}\rceil\rceil\alpha=\rceil(\lnot\alpha\land P(\pr1)(\varphi))=\lnot(\lnot\alpha\land P(\pr1)(\varphi))\land P(\pr1)(\varphi)\\
=(\lnot\lnot\alpha\lor\lnot P(\pr1)(\varphi))\land P(\pr1)(\varphi)=(\alpha\lor P(\pr1)(\lnot\varphi))\land P(\pr1)(\varphi)\\
=(\alpha\land {P(\pr1)(\varphi)})\lor {P(\pr1)({\lnot\varphi\land\varphi})}=\alpha\lor\bot=\alpha.\end{multline*}
To conclude, $(F_X,\mathfrak f)$ is Boolean since the Heyting structure is preserved by $\mathfrak{f}$.\end{proof}\end{proposition}
\subsubsection*{Variations on negation}
There are other ways to introduce negation in the context of inf-semilattices. Here we describe two examples and check that properties are again preserved.
\begin{definition}
A primary doctrine $P:\CC\op\to\Pos$ is \emph{$*$-autonomous} if for every object $A$, the poset $P(A)$ is $*$-autonomous, that is: $P(A)$ is cartesian, endowed with operation $\lnot$ such that $\lnot\lnot a=a$ for every $a\in P(A)$ and such that $a\land b\leq\lnot c$ if and only if $a\leq\lnot(b\land c)$. Moreover the operation $\lnot:P\op\to P$ yields a natural transformation.

A primary doctrine homomorphism between two $*$-autonomous doctrines is \emph{$*$-autonomous} if it preserves the negation.
\end{definition}

\begin{proposition}
Let $P:\CC\op\to\Pos$ be a primary doctrine and $P_{(X,\varphi)}$ be the Kleisli object of the comonad $(P,(X\times\blank,\mathfrak f),\gamma,\varepsilon)$ defined by the pair $X\in\CC$ and $\varphi\in P(X)$. If $P$ is a $*$-autonomous doctrine, then $P_{(X,\varphi)}$ is a $*$-autonomous doctrine, and $(F_X,\mathfrak f)$ is a $*$-autonomous homomorphism.
\begin{proof} For any $\alpha\in P_{(X,\varphi)}(A)=P(X\times A)_{\downarrow P(\pr1)(\varphi)}$, define $\rceil\alpha\coloneqq\lnot\alpha\land P(\pr1)(\varphi)$. The operation extended on each fiber $\rceil:{P_{(X,\varphi)}}\op\to{P_{(X,\varphi)}}$ is trivially a natural transformation. We call $x=P(\pr1)(\varphi)$ for simplicity, and we prove $\rceil\rceil\alpha=\alpha$.
\begin{multline*}\alpha\leq\rceil\rceil\alpha\text{ if and only if }\alpha\leq\rceil(\lnot\alpha\land\ x)\text{ if and only if }\alpha\leq\lnot(\lnot\alpha\land\ x)\land x\\
\text{ if and only if }{\alpha}\leq\lnot({\lnot\alpha}\land{x})\text{ if and only if }\alpha\land x\leq\lnot\lnot\alpha=\alpha.\end{multline*}
Conversely,
\begin{equation*}\rceil\rceil\alpha\leq\alpha\text{ if and only if }{\lnot(\lnot\alpha\land\ x)}\land {x}\leq\lnot{\lnot\alpha}\text{ if and only if }\lnot(\lnot\alpha\land\ x)\leq\lnot(x\land\lnot\alpha).\end{equation*}
Now, to prove the equivalence, take $\alpha,\beta,\gamma\leq x$, then
\begin{multline*}\alpha\land\beta\leq\rceil\gamma\text{ if and only if }\alpha\land\beta\leq\lnot\gamma\land x\text{ if and only if }\alpha\land\beta\leq\lnot\gamma\\
\text{ if and only if }\alpha\leq\lnot(\beta\land\gamma)\text{ if and only if }\alpha\leq\lnot(\beta\land\gamma)\land x=\rceil(\beta\land\gamma).\end{multline*}
To conclude, we prove that $\mathfrak{f}_A:P(A)\to P_{(X,\varphi)}(A)$ preserves the negation. On the one hand $\mathfrak{f}_A(\lnot\alpha)=P(\pr1)(\varphi)\land P(\pr2)(\lnot\alpha)=P(\pr1)(\varphi)\land \lnot P(\pr2)(\alpha)$ and, on the other hand $\rceil\mathfrak{f}_A(\alpha)=\lnot\mathfrak{f}_A(\alpha)\land P(\pr1)(\varphi)=\lnot(P(\pr1)(\varphi)\land P(\pr2)(\alpha))\land P(\pr1)(\varphi)$. To see this, it is enough to check that in a $*$-autonomous inf-semilattice we have $x\land\lnot a=\lnot(x\land a)\land x$, for any $a,x$.
First of all,
\begin{equation*}x\land\lnot a\leq\lnot(x\land a)\land x\text{ if and only if }x\land\lnot a\leq\lnot(x\land a)\text{ if and only if }x\land\lnot a\land x\leq\lnot a;\end{equation*}
conversely,
\begin{equation*}\lnot(x\land a)\land x\leq x\land\lnot a\text{ if and only if }\lnot(x\land a)\land x\leq \lnot a\text{ if and only if }\lnot(x\land a)\leq\lnot(x\land a).\qedhere\end{equation*}
\end{proof}\end{proposition}
\begin{definition}
A primary doctrine $P:\CC\op\to\Pos$ \emph{has pseudo-complements} if for every object $A$, the poset $P(A)$ has pseudo-complements, that is: the poset $P(A)$ is cartesian, endowed with an operation $\lnot$ and a bottom element $\bot$, where $\lnot a=\text{max}\{b\mid a\land b=\bot\}$; moreover the operations $\bot:\mathbf{1}\to P$, $\lnot:P\op\to P$ yield natural transformations .

A primary doctrine homomorphism between two doctrines with pseudo-complements \emph{preserves pseudo-complements} if it is bounded and preserves the negation.
\end{definition}
\begin{proposition}
Let $P:\CC\op\to\Pos$ be a primary doctrine and $P_{(X,\varphi)}$ be the Kleisli object of the comonad $(P,(X\times\blank,\mathfrak f),\gamma,\varepsilon)$ defined by the pair $X\in\CC$ and $\varphi\in P(X)$. If $P$ has pseudo-complements, then $P_{(X,\varphi)}$ has pseudo-complements, and $(F_X,\mathfrak f)$ preserves pseudo-complements.
\begin{proof}
For any $\alpha\in P_{(X,\varphi)}(A)=P(X\times A)_{\downarrow P(\pr1)(\varphi)}$, define $\rceil\alpha\coloneqq\lnot\alpha\land P(\pr1)(\varphi)$; clearly $\rceil:{P_{(X,\varphi)}}\op\to{P_{(X,\varphi)}}$ is a natural transformation. First of all we observe that $\alpha\land\rceil\alpha=\bot$, since $\alpha\land\lnot\alpha\land P(\pr1)(\varphi)=\bot$; then, suppose $\beta\leq P(\pr1)(\varphi)$ such that $\alpha\land\beta=\bot$, but from $\alpha\land\beta=\bot$, $\beta\leq\lnot\alpha$ follows in $P(X\times A)$, so $\beta\leq\lnot\alpha\land P(\pr1)(\varphi)=\rceil\alpha$, hence $\rceil\alpha=\text{max}\{\beta\mid \alpha\land\beta=\bot\}$.

To conclude, we prove that $\mathfrak{f}_A:P(A)\to P_{(X,\varphi)}(A)$ preserves the negation. Take $\mathfrak{f}_A(\lnot\alpha)$ and $\rceil\mathfrak{f}_A(\alpha)$ computed as above---in the $*$-autonomous case---, so again we check that $x\land\lnot a=\lnot(x\land a)\land x$ for any $a,x$ in a pseudo-complemented poset. First of all,
\begin{equation*}x\land\lnot a\leq\lnot(x\land a)\land x\text{ iff }x\land\lnot a\leq\lnot(x\land a)\text{, so it is sufficient }(x\land\lnot a)\land(x\land a)=\bot.\end{equation*}
Conversely,
\begin{equation*}\lnot(x\land a)\land x\leq x\land\lnot a\text{ iff }\lnot(x\land a)\land x\leq \lnot a\text{, so it is sufficient }(\lnot(x\land a)\land x)\land a=\bot.\eqno{\qedhere}\end{equation*}
\end{proof}\end{proposition}
\subsubsection*{Weak Power Objects}
Recall from Definition 4.9 in \cite{pasq} that a doctrine $P$ has weak power objects if for every object $A$ in the base category \CC, there exists an object $\pwo{A}$ and an element $\boldsymbol{\in}_A\in P(A\times\pwo{A})$ such that for any object $B$ and $\phi\in P(A\times B)$ there exists an arrow $\{\phi\}:B\to\pwo{A}$ such that $\phi=P(\id{A}\times\{\phi\})(\boldsymbol{\in}_A)$.
\begin{proposition}
Let $P:\CC\op\to\Pos$ be a primary doctrine and $P_{(X,\varphi)}$ be the Kleisli object of the comonad $(P,(X\times\blank,\mathfrak f),\gamma,\varepsilon)$ defined by the pair $X\in\CC$ and $\varphi\in P(X)$. If $P$ has weak power objects, then $P_{(X,\varphi)}$ has weak power objects.
\begin{proof}
Since $\CC_X$ has the same objects as $\CC$, for any object $A$ consider $\pwo{A}$ and the element $\mathfrak{f}_{A\times\pwo{A}}(\boldsymbol{\in}_A)\in P_{(X,\varphi)}(A\times\pwo{A})$, i.e.\ $P(\pr1)\varphi\land P(\ple{ \pr2,\pr3})(\boldsymbol{\in}_A)\in P(X\times A\times\pwo{A})_{\downarrow P(\pr1)\varphi}$. We want to prove that $(\pwo{A},\mathfrak{f}_{A\times\pwo{A}}(\boldsymbol{\in}_A))$ is a weak power object of $A$ in the doctrine $P_{(X,\varphi)}$.

To see this, we take any object $C$ and $\psi\in P_{(X,\varphi)}(A\times C)=P(X\times A\times C)_{\downarrow P(\pr1)\varphi}$ and look for an arrow $[\psi]:C\rightsquigarrow\pwo{A}$---i.e.\ a \CC-arrow $X\times C\to \pwo{A}$---such that
\begin{equation*}\psi=P_{(X,\varphi)}(\id{A}\times_{\CC_X}[\psi])\mathfrak{f}_{A\times\pwo{A}}(\boldsymbol{\in}_A).\end{equation*}
Here, the product $\id{A}\times_{\CC_X}[\psi]$ is computed in $\CC_X$, hence it is actually the \CC-arrow
\begin{equation*}\ple{ \pr2,[\psi]\ple{ \pr1,\pr3}}:X\times A\times C\to A\times\pwo{A}.\end{equation*}

Using the fact that $P$ has weak power object, we can take the object $X\times C$ and the element $P(\ple{ \pr2,\pr1,\pr3})\psi\in P(A\times X\times C)$ and we know that there exists
\begin{equation*} [\psi]\coloneqq\{P(\ple{ \pr2,\pr1,\pr3})\psi\}:X\times C\to \pwo{A}\end{equation*}
such that
\begin{equation*}P(\ple{ \pr2,\pr1,\pr3})\psi=P\left(\id{A}\times\{P(\ple{ \pr2,\pr1,\pr3})\psi\}\right)(\boldsymbol{\in}_A)=P\left(\ple{ \pr1,[\psi]\ple{ \pr2,\pr3}}\right)(\boldsymbol{\in}_A).\end{equation*}

Now compute
\begin{align*}P_{(X,\varphi)}(\id{A}\times_{\CC_X}[\psi])\mathfrak{f}_{A\times\pwo{A}}(\boldsymbol{\in}_A)&=P(\ple{ \pr1,\pr2,[\psi]\ple{ \pr1,\pr3}})\big(P(\pr1)\varphi\land P(\ple{ \pr2,\pr3})(\boldsymbol{\in}_A)\big)\\
&=P(\pr1)\varphi\land P\left(\ple{ \pr2,[\psi]\ple{ \pr1,\pr3}}\right)(\boldsymbol{\in}_A)\end{align*}
which is equal to $\psi$ if and only if $P(\ple{ \pr2,\pr1,\pr3})\psi=P(\pr2)\varphi\land P\left(\ple{\pr1,[\psi]\ple{ \pr2,\pr3}}\right)(\boldsymbol{\in}_A)$, but this is true following from the definition of $[\psi]$ and the fact that $\psi\leq P(\pr1)\varphi$.\end{proof}\end{proposition}

\section{Universal properties of $P_{(X,\varphi)}$}\label{sect:univ_prop}
Consider the following diagram for a primary doctrine $P$. In particular, $P_{(X,\varphi)}$ is primary too.
\begin{center}
\begin{tikzcd}
\CC^{\text{op}}\arrow[rr,"{F_X}^{\text{op}}"] \arrow[dr,"P"' ,""{name=L}]&&{\CC_X}^{\text{op}}\arrow[dl,"P_{(X,\varphi)}" ,""'{name=R}]\\
&\Pos\arrow[rightarrow,"\mathfrak{f}","\cdot"', from=L, to=R, bend left=10]
\end{tikzcd}
\end{center}
We can interpret this 1-arrow as follows: we are adding a constant of sort $X$ to the theory $P$, and making this constant verify $\varphi$.
Indeed, take $\mathfrak{f}_X(\varphi)\in P_{(X,\varphi)}(X)$, which is the interpretation of $\varphi$ in $P_{(X,\varphi)}$, and consider the constant $\pr1:\tmn\rightsquigarrow X$ in $\CC_X$.
\begin{notation}
When there is no confusion, the terminal object of a given category will be simply called $\tmn$. Otherwise, a subscript will specify the category in which we are computing the terminal object.
\end{notation}
This map is the $\CC$-arrow $\pr1:X\times \tmn\to X$, which is a direction of the canonical isomorphism $X\times\tmn\cong X$, whose inverse is given by $\ple{\id{X},!_X}:X\to X\times\tmn$. This induces an isomorphism also between the corresponding fibers $P(X\times \tmn)\cong P(X)$, so that $P_{(X,\varphi)}(\tmn)=P(X\times\tmn)_{\downarrow P(\pr1)\varphi}\cong P(X)_{\downarrow\varphi}$. From now on we will write $\id{X}:\tmn\rightsquigarrow X$ instead of $\pr1:\tmn\rightsquigarrow X$, and $P(X)_{\downarrow\varphi}$ instead of $P(X\times\tmn)_{\downarrow P(\pr1)\varphi}$. With this notation, we compute the reindexing of $\mathfrak{f}_X(\varphi)$ along the constant $\id{X}:\tmn\rightsquigarrow X$, and we show that it is the top element in $P_{(X,\varphi)}(\tmn)=P(X)_{\downarrow\varphi}$. Indeed:
\begin{equation*}P_{(X,\varphi)}(\id{X})\mathfrak{f}_X(\varphi)=P(\Delta_X)\mathfrak{f}_X(\varphi)=P(\Delta_X)(P(\pr1)\varphi\land P(\pr2)\varphi)=\varphi\land\varphi=\varphi,\end{equation*}
and $\varphi$ is the top element of $P_{(X,\varphi)}(\tmn)=P(X)_{\downarrow\varphi}$, which means, the interpretation of $\varphi$ evaluated in the new constant is true.

\begin{theorem}\label{thm:univ_prop}
Let $P:\CC\op\to\Pos$ be a primary doctrine. Given an object $X$ in the base category \CC and an element $\varphi\in P(X)$, the 1-arrow $(F_X,\mathfrak f):P\to P_{(X,\varphi)}$ and the $\CC_X$-arrow $\id{X}:\tmn_{\CC_X}\rightsquigarrow X$ are such that $\top\leq P_{(X,\varphi)}(\id{X})\mathfrak{f}_X(\varphi)$ in $P_{(X,\varphi)}(\tmn_{\CC_X})$, and they are universal with respect to this property, i.e.\ for any primary 1-arrow $(G,\mathfrak{g}):P\to R$, where $R:\ct{D}\op\to\Pos$ is a primary doctrine, and any $\ct{D}$-arrow $c:\tmn_\ct{D}\to G(X)$ such that $\top\leq R(c)\mathfrak{g}_X(\varphi)$ in $R(\tmn_\ct{D})$ there exists a unique up to a unique natural isomorphism primary 1-arrow $(G',\mathfrak{g}'):P_{(X,\varphi)}\to R$ such that $(G',\mathfrak{g}')\circ(F_X,\mathfrak{f})=(G,\mathfrak{g})$ and $G'(\id{X})=c$.
\begin{proof}
Consider the diagram
\begin{equation}\label{eq:cmd}
\begin{tikzcd}
(P,(X\times \blank,\mathfrak{f}),\gamma, \varepsilon)\arrow[dr, bend right=15pt,"{((G,\mathfrak{g}),\mathfrak{j})}"']\arrow[rr, "{((F_X,\mathfrak{f}),\gamma)}"]&&(P_{(X,\varphi)},(\id{},\id{}),\id{}, \id{})\arrow[dl,dashed, bend left=15pt,"{((G',\mathfrak{g}'),\id{})}"]\\
&(R,(\id{},\id{}),\id{}, \id{})
\end{tikzcd}
\end{equation}
describing the universal property for the Kleisli construction for the comonad we are studying on $P$---see diagram \eqref{eq:uni_prop} in Proposition \ref{prop:kleisli}.

So, in order to construct $(G',\mathfrak{g}'):P_{(X,\varphi)}\to R$, we must define $\mathfrak{j}$ in such a way that $((G,\mathfrak{g}),\mathfrak{j})$ is an arrow in $\Cmd^*(\QDot)$ as in \eqref{eq:cmd}, i.e.\ a natural transformation $\mathfrak{j}:G\xrightarrow{\cdot} G(X\times\blank)$, such that $\mathfrak{g}_A\leq R(\mathfrak{j}_A)\mathfrak{g}_{X\times A}\mathfrak{f}_A$ and satisfying the coherence diagrams.
Knowing that $G$ preserves products, we need to define for every object $A$ an arrow $\mathfrak{j}_A:GA\to GX\times GA$ take $\mathfrak{j}_A\coloneqq\ple{ c\cdot!_{GA}, \id{GA}}$, where $!_{GA}:GA\to\mathbf {t}_\ct{D}$ is the unique arrow from $GA$ to the teminal object.

This is a natural transformation:
\begin{center}
\begin{tikzcd}
A\arrow[d,"f"]&GA\arrow[r,"\mathfrak{j}_A"]\arrow[d,"G(f)"]&GX\times GA\arrow[d,"{G(\id{}\times f)=\id{}\times G(f)}"]\\
B&GB\arrow[r,"\mathfrak{j}_B"]&GX\times GB
\end{tikzcd}.
\end{center}
Indeed, for any $f:A\to B$, we have:
\begin{equation*}(\id{}\times G(f))\ple{ c\cdot!_{GA}, \id{GA}}=\ple{ c\cdot!_{GA}, G(f)}=\ple{ c\cdot!_{GB}, \id{GB}} G(f).\end{equation*}
Moreover, for any $\alpha\in P(A)$, we have
\begin{align*}R(\mathfrak{j}_A)\mathfrak{g}_{X\times A}\mathfrak{f}_A(\alpha)&=R(\mathfrak{j}_A)\mathfrak{g}_{X\times A}(P(\pr1)(\varphi)\land P(\pr2)(\alpha))\\
&=R(\mathfrak{j}_A)\big(R{G(\pr1)}\mathfrak{g}_X(\varphi)\land R{G(\pr2)}\mathfrak{g}_A(\alpha)\big)=R(c\cdot !_{GA})\mathfrak{g}_X(\varphi)\land \mathfrak{g}_A(\alpha),\end{align*}
using naturality of $\mathfrak g$ and the fact that $G$ preserves products.

So now observe that $\mathfrak{g}_A(\alpha)\leq R(c\cdot !_{GA})\mathfrak{g}_X(\varphi)\land \mathfrak{g}_A(\alpha)$ if and only if $\mathfrak{g}_A(\alpha)\leq R(c\cdot !_{GA})\mathfrak{g}_X(\varphi)$, but by assumption $R(c)\mathfrak{g}_X(\varphi)=\top$ in $R(\mathbf{t}_\ct{D})$, so $R(!_{GA})R(c)\mathfrak{g}_X(\varphi)=R(c\cdot!_{GA})\mathfrak{g}_X(\varphi)=\top$ in $R(GA)$, hence the inequality holds.

To conclude, we prove that the coherence diagrams commute:
\begin{center}
\begin{tikzcd}
G\arrow[r,"\mathfrak{j}"]\arrow[d,"\mathfrak{j}"]&GX\times G\blank\arrow[d,"G(\gamma)"]&G\arrow[r,"\mathfrak{j}"]\arrow[rd,"\id{}",bend right]&GX\times G\blank\arrow[d,"G(\varepsilon)"]\\
GX\times G\blank\arrow[r,"\mathfrak{j}_{X\times\blank}"]&GX\times GX\times G\blank&&G
\end{tikzcd}.
\end{center}
The first diagram commutes since
\begin{align*}G(\Delta_X\times \id{A})\mathfrak{j}_A&=(\Delta_{GX}\times \id{GA})\ple{ c\cdot !_{GA},\id{GA}}=\ple{ c\cdot !_{GA},c\cdot !_{GA},\id{GA}}\\
&=\ple{ c\cdot !_{G(X\times A)},\id{G(X\times A)}}\ple{ c\cdot !_{GA},\id{GA}},\end{align*}
while the second one commutes since
\begin{equation*}G(\varepsilon_A)\mathfrak{j}_A=G(\pr2)\mathfrak{j}_A=\pr2'\mathfrak{j}_A=\id{GA}.\end{equation*}
So now we know from the universal property that there exists a unique $(G',\mathfrak{g}'):P_{(X,\varphi)}\to R$ such that $((G,\mathfrak{g}),\mathfrak j)=((G',\mathfrak{g}'),\id{})((F_K,\mathfrak f),\gamma)$. In particular there is an arrow $(G',\mathfrak{g}'):P_{(X,\varphi)}\to R$ such that $(G',\mathfrak{g}')\circ(F_K,\mathfrak{f})=(G,\mathfrak{g})$. Moreover, if we translate the universal property in Proposition \ref{prop:kleisli} to the notation used here, we observe that $G'(g:A\rightsquigarrow B)=G(g)\mathfrak{j}_A$, where $g:X\times A\to B$ is an arrow in \CC. In particular, taking $A=\tmn_{\CC_X}$ and $g=\id{X}:X\to X$, we obtain $G'(\id{X})=\mathfrak{j}_{\tmn_\CC}=c$. Here we use the fact that $G$ preserves the terminal object, and that the product of an object with the terminal object is the object itself.

We use the definition of $G'$ on arrows to prove that $G'$ preserves products; consider the following three diagrams: on the left there is the \CC-diagram that mirrors a product diagram in $\CC_X$---in the middle---, while on the right there is the image of such product through $G'$.
\begin{center}
\begin{tikzcd}[column sep=5pt]
&X\times A\times B\arrow[dl, bend right, "\pr2"]\arrow[dr, bend left, "\pr3"]&&&A\times B\arrow[dl,squiggly, bend right, "\pr2"]\arrow[dr,squiggly, bend left, "\pr3"]&&&GA\times GB\arrow[dl, bend right, "G(\pr2)\mathfrak{j}_{A\times B}"]\arrow[dr, bend left, "G(\pr3)\mathfrak{j}_{A\times B}"]&\\
A&&B&A&&B&GA&&GB
\end{tikzcd}
\end{center}
However, since $G$ preserves products, $G(\pr2)$ and $G(\pr3)$ are respectively the second and third projections from $GX\times GA\times GB$, and these precomposed with $\mathfrak{j}_{A\times B}$ are precisely the first and second projections from $GA\times GB$, as claimed.

To show that $\mathfrak{g}'$ preserves infima and top element, recall from Proposition \ref{prop:kleisli} that $\mathfrak{g}'_A$ is the restriction of $R(\mathfrak{j}_A)\mathfrak{g}_{X\times A}$. Then notice that for any $\alpha, \beta\in P_{(X,\varphi)}(A)$ we have
\begin{equation*}\mathfrak{g}'_A(\alpha\land\beta)=R(\mathfrak{j}_A)\mathfrak{g}_{X\times A}(\alpha\land\beta)=R(\mathfrak{j}_A)\big(\mathfrak{g}_{X\times A}(\alpha)\land\mathfrak{g}_{X\times A}(\beta)\big)=\mathfrak{g}'_A(\alpha)\land\mathfrak{g}'_A(\beta),\end{equation*}
since by assumption $\mathfrak{g}$ respect the structures, while for the top element:
\begin{align*}\mathfrak{g}'_A(P(\pr1)\varphi)&=R(\mathfrak{j}_A)\mathfrak{g}_{X\times A}P(\pr1)\varphi=R(\ple{ c \cdot !_{GA},\id{GA}})R(\pr1)\mathfrak{g}_X\varphi\\
&=R(!_{GA})R(c)\mathfrak{g}_X\varphi=\top_{GA}.\end{align*}
So $(G',\mathfrak{g}')$ is indeed a primary 1-arrow.

Finally, suppose that $(\overline{G}, \overline{\mathfrak{g}}):P_{(X,\varphi)}\to R$ is another primary 1-arrow such that $(\overline{G}, \overline{\mathfrak{g}})\circ(F_X,\mathfrak{f})=(G,\mathfrak{g})$ and $\overline{G}(\id{X})=c$.

Then we can compute the composition $((G,\mathfrak{g}),\overline{\mathfrak j})=((\overline{G},\overline{\mathfrak{g}}),\id{})((F_K,\mathfrak f),\gamma)$, where we define $\overline{\mathfrak{j}}_A:=(\id{}\circ\gamma)_A=\overline{G}(\varepsilon_{X\times A}\gamma_A)=\overline{G}(\id{X\times A})$---see Remark \ref{rmk:squiggly} and Proposition \ref{prop:kleisli}. We claim that $\overline {\mathfrak{j}}=\mathfrak{j}$, so that by uniqueness given by the universal property, the equality $(G',\mathfrak{g}')=(\overline{G},\overline{\mathfrak{g}})$ follows.

In our notation, we have to think of $\id{X\times A}$ as the map $A\rightsquigarrow X\times A$, and it is uniquely defined by its two components: the first one is $\pr1:A\rightsquigarrow X$, the second one is $\pr2:A\rightsquigarrow A$. Observe that $\pr2$ is the identity of $A$ in $\CC_X$, while $\pr1:A\rightsquigarrow X$ is the composition of the unique arrow $A\rightsquigarrow\tmn_\CC$ and the constant $\id{X}:\tmn_{\CC_X}\rightsquigarrow X$. Since $\overline G$ preserves products, $\overline{\mathfrak{j}}_A=\overline G(\id{X\times A}):GA\to GX\times GA$ must be the identity of $GA$ on the second component; in particular the second components of $\overline{\mathfrak{j}}_A$ and $\mathfrak{j}_A$ are the same. Concerning the first component we have $\overline{G}(\pr1)=\overline{G}(\id{X})\overline{G}(!_{X\times A})=c\cdot !_{GA}$, i.e.\ also the first component of $\overline{\mathfrak{j}}_A$ coincides with the first component of $\mathfrak{j}_A$, hence the two maps coincide as claimed.\end{proof}
\end{theorem}
\begin{theorem}\label{thm:prop}
Let $P$, $P_{(X,\varphi)}$, $R$, $(G,\mathfrak{g}):P\to R$ be the doctrines and a morphism with the same assumption of Theorem \ref{thm:univ_prop}. Then
\begin{enumerate}[label=(\roman*)]
\item if $P,R$ and $(G,\mathfrak{g})$ are elementary, then $\mathfrak{g}'$ preserves the fibered equality;
\item if $P,R$ and $(G,\mathfrak{g})$ are existential, then $\mathfrak{g}'$ preserves the existential quantifier;
\item if $P,R$ and $(G,\mathfrak{g})$ are universal, then $\mathfrak{g}'$ preserves the universal quantifier;
\item if $P,R$ and $(G,\mathfrak{g})$ are implicational, then $\mathfrak{g}'$ preserves the implication;
\item if $P,R$ are bounded, with top and bottom elements preserved by $\mathfrak{g}$, then $\mathfrak{g}'$ preserves them;
\item if $P,R$ have binary joins, preserved by $\mathfrak{g}$, then $\mathfrak{g}'$ preserves binary joins;
\item if $P,R$ and $(G,\mathfrak{g})$ are respectively Heyting or Boolean, then $\mathfrak{g}'$ preserves the corresponding structure.
\end{enumerate}
\begin{proof}
\begin{enumerate}[label=(\roman*)]
\item We need to check that $\mathfrak{g'}_{A\times A}(\delta^{X,\varphi})_A=\delta_{GA}$.
\begin{align*}\mathfrak{g'}_{A\times A}(\delta^{X,\varphi})_A&=R(\mathfrak{j}_{A\times A})\mathfrak{g}_{X\times A\times A}\big(P(\pr1)(\varphi)\land P(\ple{\pr2,\pr3})(\delta_A)\big)\\
&=R(\mathfrak{j}_{A\times A})R(\pr1)\mathfrak{g}_X(\varphi)\land R(\ple{c!_{GA\times GA},\pr2,\pr3})R(\ple{\pr2,\pr3})\mathfrak{g}_{A\times A}(\delta_A)\\
&=R(!_{GA\times GA})R(c)\mathfrak{g}_X(\varphi)\land \mathfrak{g}_{A\times A}(\delta_A)=\top_{GA\times GA}\land\delta_{GA}=\delta_{GA}.\end{align*}
\item We prove that $\exists^{GB}_{GC}\mathfrak{g'}_{C\times B}=\mathfrak{g}'_{C}{\exists_{X,\varphi}}^B_C$:
\begin{align*}\exists^{GB}_{GC}\mathfrak{g'}_{C\times B}&=\exists^{GB}_{GC}R({\mathfrak{j}_{C\times B}})\mathfrak{g}_{X\times C\times B}=R(\mathfrak{j}_C)\exists^{GB}_{GX\times GC}\mathfrak{g}_{X\times C\times B}\\
&=R(\mathfrak{j}_C)\mathfrak{g}_{X\times C}\exists^B_{X\times C}=\mathfrak{g}'_{C}{\exists_{X,\varphi}}^B_C.\end{align*}
\item The proof is similar to the one above, with a little alteration:
\begin{align*}\forall^{GB}_{GC}\mathfrak{g'}_{C\times B}&=\forall^{GB}_{GC}R(\mathfrak{j}_{C\times B})\mathfrak{g}_{X\times C\times B}\\
&=R(\mathfrak{j}_C)\forall^{GB}_{GX\times GC}\mathfrak{g}_{X\times C\times B}=R(\mathfrak{j}_C)\mathfrak{g}_{X\times C}\forall^B_{X\times C}\end{align*}
and also
\begin{align*}\mathfrak{g}'_{C}{\forall_{X,\varphi}}^B_C&=R(\mathfrak{j}_C)\mathfrak{g}_{X\times C}\big(\forall^B_{X\times C}(\blank)\land P(\pr1)\varphi\big)\\
&=R(\mathfrak{j}_C)\mathfrak{g}_{X\times C}\forall^B_{X\times C}(\blank)\land{R(\mathfrak{j}_C)\mathfrak{g}_{X\times C}P(\pr1)\varphi}=R(\mathfrak{j}_C)\mathfrak{g}_{X\times C}\forall^B_{X\times C}.\end{align*}
Note that ${R(\mathfrak{j}_C)\mathfrak{g}_{X\times C}P(\pr1)\varphi}=\top_{GC}$ since $\mathfrak{g}'$ preserves the top element---see Theorem \ref{thm:univ_prop}.
\item Take $\alpha, \beta\in P_{(X,\varphi)}(A)$:
\begin{align*}\mathfrak{g}'_A(\alpha\Rightarrow\beta)&=R(\mathfrak{j}_A)\mathfrak{g}_{X\times A}((\alpha\rightarrow\beta)\land P(\pr1)\varphi)\\
&=\big(R(\mathfrak{j}_A)\mathfrak{g}_{X\times A}(\alpha)\rightarrow R(\mathfrak{j}_A)\mathfrak{g}_{X\times A}(\beta)\big)\land {R(\mathfrak{j}_A)\mathfrak{g}_{X\times A}P(\pr1)\varphi}=\mathfrak{g}'_A(\alpha)\rightarrow\mathfrak{g}'_A(\beta).\end{align*}
\item Compute $\mathfrak{g}'_A(\ini_A)=R(\mathfrak{j}_A)\mathfrak{g}_{X\times A}(\bot_A)=\bot_{GA}$.
\item Take $\alpha, \beta\in P_{(X,\varphi)}(A)$, then
\begin{equation*}\mathfrak{g}'_A(\alpha\lor\beta)=R(\mathfrak{j}_A)\mathfrak{g}_{X\times A}(\alpha\lor\beta)=\mathfrak{g}'_A(\alpha)\lor\mathfrak{g}'_A(\beta).\end{equation*}
Observe that in order to prove this point it was not necessary to ask for the condition that $P(\pr1)\varphi\land(\blank)$ preserves finite joins in $P(X\times A)$, which was necessary for $\mathfrak{f}$ to preserve finite joins in Proposition \ref{prop:joins}: it is enough to ask $\mathfrak{g}$ to preserve them.
\item It follows trivially combining the previous properties.\qedhere
\end{enumerate}\end{proof}
\end{theorem}
A stronger result for Theorem \ref{thm:univ_prop} holds. Let again $P:\CC\op\to\Pos$ be a primary doctrine; fix on object $X$ in the base category, and an element $\varphi\in P(X)$. For any other primary doctrine $R:\ct{D}\op\to\Pos$, define the category $\PD_{(X,\varphi)}(P,R)$ whose objects are pairs of the kind $\big((G,\mathfrak g), c:\tmn_\ct{D}\to GX\big)$, where $(G,\mathfrak g)\in\PD(P,R)$, such that $\top\leq R(c)\mathfrak{g}_X(\varphi)$ in $R(\tmn_\ct{D})$, and whose arrows are 2-arrows preserving the constant, meaning $\theta:\big((G,\mathfrak g),c\big)\to\big((H,\mathfrak h),d\big)$ is a 2-arrow $\theta:(G,\mathfrak g)\to(H,\mathfrak h)$ in $\PD$ such that $c\theta_X=d$. There is an obvious functor induced by precomposition with $(F_X,\mathfrak f)$, from $\PD(P_{(X,\varphi)},R)$ to $\PD_{(X,\varphi)}(P,R)$: it maps any $\xi:(K,\mathfrak k)\to(K',\mathfrak k ')$ into $\xi_{F_X}:\big((K,\mathfrak k)\circ(F_X,\mathfrak f),K(\id{X})\big)\to\big((K',\mathfrak k')\circ(F_X,\mathfrak f),K'(\id{X})\big)$. This is well defined on objects since $R(K(\id{X}))(\mathfrak{k}_{F_XX}\mathfrak{f}_X\varphi)=\mathfrak{k}_\tmn P_{(X,\varphi)}(\id{X})\mathfrak{f}_X\varphi=\top$, and well defined on arrows since $K(\id{X})\xi_X=K'(\id{X})$ by naturality.
\begin{theorem}\label{thm:prec}
Let $P$ be a primary doctrine. Given an object $X$ in the base category and an element $\varphi\in P(X)$, the functor $\PD(P_{(X,\varphi)},R)\to\PD_{(X,\varphi)}(P,R)$ induced by precomposition with $(F_X,\mathfrak f)$ is an equivalence of categories for any primary doctrine $R$.
\begin{proof}
The functor is essentially surjective following from Theorem \ref{thm:univ_prop} and faithfulness is trivial since $F_X$ is the identity on objects. To show that the functor is full, take any 2-arrow $\theta:\big((K,\mathfrak k)\circ(F_X,\mathfrak f),K(\id{X})\big)\to\big((K',\mathfrak k')\circ(F_X,\mathfrak f),K'(\id{X})\big)$ and prove that $\theta:(K,\mathfrak k)\to(K',\mathfrak{k}')$ is in $\PD(P_{(X,\varphi)},R)$. First of all, we check that it is a natural transformation $K\to K'$: take any $f:A\rightsquigarrow B$ in $\CC_X$ and break it as the composition of $\id{X\times A}:A\rightsquigarrow X\times A$ and $f\ple{\pr2,\pr3}=F_X(f):X\times A\rightsquigarrow B$; moreover observe that $\id{X\times A}$ has as first projection the composition of the unique arrow $A\rightsquigarrow\tmn$ and the constant $\id{X}:\tmn\rightsquigarrow X$, and as second projection the identity $\pr2:A\rightsquigarrow A$---see the end of the proof of Theorem \ref{thm:univ_prop}. So the naturality diagram becomes:
\[\begin{tikzcd}
	A && KA &&&& {K'A} \\
	{X\times A} && {KX\times KA} &&&& {K'X\times K'A} \\
	B && KB &&&& {K'B}
	\arrow["{\theta_A}", from=1-3, to=1-7]
	\arrow["{\ple{K(\id{X})!_{KA},\id{KA}}}", from=1-3, to=2-3]
	\arrow["{KF_X(f)}", from=2-3, to=3-3]
	\arrow["{\theta_B}"', from=3-3, to=3-7]
	\arrow["{\ple{K'(\id{X})!_{K'A},\id{K'A}}}"', from=1-7, to=2-7]
	\arrow["{K'F_X(f)}"', from=2-7, to=3-7]
	\arrow["{\id{X\times A}}", squiggly, from=1-1, to=2-1]
	\arrow["{F_X(f)}", squiggly, from=2-1, to=3-1]
	\arrow["f"', curve={height=24pt}, squiggly, from=1-1, to=3-1]
	\arrow["{\theta_X\times\theta_A}"{description}, from=2-3, to=2-7]
	\arrow["{K(f)}"', curve={height=40pt}, from=1-3, to=3-3]
	\arrow["{K'(f)}", curve={height=-40pt}, from=1-7, to=3-7]
\end{tikzcd}.\]
The lower square commutes since $\theta:KF_X\xrightarrow{\cdot} K'F_X$ by assumption, while the upper square commutes since $K(\id{X})\theta_X=K'(\id{X})$. To conclude, we need for any $\CC$-object $A$ and any $\alpha\in P_{(X,\varphi)}(A)=P(X\times A)_{\downarrow P(\pr1)\varphi}$ the inequality $\mathfrak{k}_A(\alpha)\leq R(\theta_A)\mathfrak{k}'_A(\alpha)$ to hold. In particular $\alpha\in P(X\times A)$, so we can consider
\begin{equation*}\mathfrak{f}_{X\times A}\alpha=P(\pr1)(\varphi)\land P(\ple{\pr2,\pr3})(\alpha)\in P_{(X,\varphi)}(X\times A)\subseteq P(X\times X\times A);\end{equation*}
apply then naturality of $\mathfrak{k}$ to $\id{X\times A}:A\rightsquigarrow X\times A$ to observe that
\begin{multline*}R(\ple{K(\id{X})!_{KA},\id{KA}})\mathfrak{k}_{X\times A}\mathfrak{f}_{X\times A}\alpha=\mathfrak{k}_AP_{(X,\varphi)}(\id{X\times A})\mathfrak{f}_{X\times A}\alpha\\
=\mathfrak{k}_AP(\ple{\pr1,\pr1,\pr2})\big(P(\pr1)(\varphi)\land P(\ple{\pr2,\pr3})(\alpha)\big)=\mathfrak{k}_A(P(\pr1)(\varphi)\land \alpha)=\mathfrak{k}_A(\alpha)\end{multline*}
since $\alpha\leq P(\pr1)(\varphi)$. Moreover, since in particular $\theta:(K,\mathfrak k)\circ(F_X,\mathfrak f)\to(K',\mathfrak k')\circ(F_X,\mathfrak f)$, we know that $\mathfrak{k}_{X\times A}\mathfrak{f}_{X\times A}(\alpha)\leq R(\theta_{X\times A})\mathfrak{k}'_{X\times A}\mathfrak{f}_{X\times A}(\alpha)$. So we have:
\begin{align*}\mathfrak{k}_A(\alpha)&=R(\ple{K(\id{X})!_{KA},\id{KA}})\mathfrak{k}_{X\times A}\mathfrak{f}_{X\times A}\alpha\\
&\leq R(\ple{K(\id{X})!_{KA},\id{KA}})R(\theta_{X\times A})\mathfrak{k}'_{X\times A}\mathfrak{f}_{X\times A}(\alpha)\\
&=R(\theta_{A})R(\ple{K'(\id{X})!_{K'A},\id{K'A}})\mathfrak{k}'_{X\times A}\mathfrak{f}_{X\times A}(\alpha)=R(\theta_A)\mathfrak{k}'_A(\alpha).\qedhere\end{align*}
\end{proof}
\end{theorem}
The process studied above in Theorem \ref{thm:univ_prop} describes how to add a constant of sort $X$ that verifies a formula $\varphi$ in a universal way. Taking the particular case when $X=\tmn$ is the terminal object, we are not adding any constant---the functor $\tmn\times\blank:\CC\to\CC$ is essentially the identity---, and we are just requiring $\varphi\in P(\tmn)$ to be true in the new doctrine---i.e.\ we are adding the axiom $\varphi$ to the theory $P$---, in a universal way. In this case we write $P_\varphi:\CC\op\to\Pos$; for any given $\CC$-arrow $f:A\to B$, we have $P_\varphi(f):P_\varphi(B)\to P_\varphi(A)$, computed as
\begin{equation*}P(f):P(B)_{\downarrow P(!_B)\varphi}\to P(A)_{\downarrow P(!_A)\varphi}.\end{equation*}
The 1-arrow $P\to P_\varphi$ becomes $(\id{\CC}, \mathfrak{f})$, where $\mathfrak{f}_A:P(A)\to P_\varphi(A)$ maps an element $\alpha\in P(A)$ to $P(!_A)\varphi\land \alpha\in P_\varphi(A)$. All the additional properties of $P$ described in Section \ref{sect:pres_prop} are clearly recovered by $P_\varphi$.
\begin{corollary}\label{coroll:ax}
Let $P:\CC\op\to\Pos$ be a primary doctrine. Given an element $\varphi\in P(\tmn)$, the 1-arrow $(\id{\CC},\mathfrak f):P\to P_\varphi$ is such that $\top\leq \mathfrak{f}_\tmn(\varphi)$ in $P_\varphi(\tmn)$, and it is universal with respect to this property, i.e.\ for any primary 1-arrow $(G,\mathfrak{g}):P\to R$, where $R:\ct{D}\op\to\Pos$ is a primary doctrine, such that $\top\leq \mathfrak{g}_\tmn(\varphi)$ in $R(\tmn_\ct{D})$ there exists a unique up to a unique natural isomorphism primary 1-arrow $(G',\mathfrak{g}'):P_{\varphi}\to R$ such that $(G',\mathfrak{g}')\circ(\id{\CC},\mathfrak f)=(G,\mathfrak{g})$.
\end{corollary}

The category corresponding to $\PD_{(X,\varphi)}(P,R)$ in Theorem \ref{thm:prec} for some primary doctrine $R:\ct{D}\op\to\Pos$ in this case is called $\PD_\varphi(P,R)$: objects are morphisms $(G,\mathfrak{g})\in\PD(P,R)$ such that $\top\leq\mathfrak{g}_\tmn(\varphi)$ in $R(\tmn_\ct{D})$ and arrows are 2-arrows of $\PD$. In particular $\PD_\varphi(P,R)$ is a full subcategory of $\PD(P,R)$. Precomposition with $(\id{\CC},\mathfrak f)$ is a functor from $\PD(P_\varphi,R)$ to $\PD(P,R)$, and has image in $\PD_\varphi(P,R)$: given $(K,\mathfrak k):P_\varphi\to R$, the composition $(K,\mathfrak k)(\id{\CC},\mathfrak f)$ is such that $\mathfrak {(kf)}_\tmn(\varphi)=\mathfrak{k}_\tmn\mathfrak f _\tmn(\varphi)=\mathfrak{k}_\tmn(\varphi)=\top$ in $R(\tmn_\ct{D})$, since $\varphi$ is the top element in $P_\varphi(\tmn)$.
\begin{corollary}
Let $P$ be a primary doctrine. Given an element $\varphi\in P(\tmn)$, precomposition with $(\id{\CC},\mathfrak f)$
\begin{equation*}\blank\circ(\id{\CC},\mathfrak f):\PD(P_\varphi,R)\to\PD_\varphi(P,R)\end{equation*}
is an equivalence of categories for any primary doctrine $R$.
\end{corollary}
Similarly, we can take the particular case when $\varphi=\top_X\in P(X)$ is the top element, so we are not making any formula true---the natural transformation $\mathfrak{f}_A:P(A)\to P(X\times A)$ represent the inclusion of formulae of sort $A$ in the formulae of the same sort but in a language with a new constant---, and we are just adding a constant $\id{X}:\tmn\rightsquigarrow X$, in a universal way. In this case we write $P_X:\CC_X\op\to\Pos$; for any given $\CC_X$-arrow $f:A\rightsquigarrow B$, we have $P_X(f):P_X(B)\to P_X(A)$ computed as
\begin{equation*}P(\ple{\pr1,f}):P(X\times B)\to P(X\times A).\end{equation*}
The 1-arrow $P\to P_X$ becomes $(F_X, \mathfrak{f})$, where $\mathfrak{f}_A:P(A)\to P_X(A)$ maps an element $\alpha\in P(A)$ to $P(\pr2)(\alpha)\in P_X(A)$. All the additional properties of $P$ described in Section \ref{sect:pres_prop} are clearly recovered by $P_X$. Observe that for this construction, the assumption that the starting doctrine $P$ is primary is not needed: the reader can go through all proofs removing the primary assumption on the starting doctrine when $\varphi=\top_X$ and every computation still works.
\begin{corollary}\label{coroll:const}
Let $P:\CC\op\to\Pos$ be a doctrine. Given an object $X$ in the base category, the 1-arrow $(F_X,\mathfrak f):P\to P_X$ and the $\CC_X$-arrow $\id{X}:\tmn_{\CC_X}\rightsquigarrow X$ are universal, i.e.\ for any 1-arrow $(G,\mathfrak{g}):P\to R$, where $R:\ct{D}\op\to\Pos$ is a doctrine, and any $\ct{D}$-arrow $c:\tmn_\ct{D}\to G(X)$ there exists a unique up to a unique natural isomorphism 1-arrow $(G',\mathfrak{g}'):P_X\to R$ such that $(G',\mathfrak{g}')\circ(F_X,\mathfrak f)=(G,\mathfrak{g})$ and $G'(\id{X})=c$.
\end{corollary}
The category corresponding to $\PD_{(X,\varphi)}(P,R)$ in Theorem \ref{thm:prec} can be defined for any doctrine $R:\ct{D}\op\to\Pos$, and in this case is called $\Dott_X(P,R)$: objects are pairs $\big((G,\mathfrak{g}),c:\tmn\to GX\big)$ where $(G,\mathfrak{g})\in\Dott(P,R)$ and arrows are 2-arrows of $\Dott$ preserving the constant. Precomposition with $(F_X,\mathfrak f)$ induces a functor from $\Dott(P_X,R)$ to $\Dott_X(P,R)$: it maps any $\xi:(K,\mathfrak k)\to(K',\mathfrak k ')$ into $\xi_{F_X}:\big((K,\mathfrak k)\circ(F_X,\mathfrak f),K(\id{X})\big)\to\big((K',\mathfrak k')\circ(F_X,\mathfrak f),K'(\id{X})\big)$. This is well defined on arrows since $K(\id{X})\xi_X=K'(\id{X})$ by naturality.
\begin{corollary}
Let $P$ be a doctrine. Given an object $X$ in the base category, the functor $\ple{\blank\circ(F_X,\mathfrak f),\blank(\id{X})}:\Dott(P_X,R)\to\Dott_X(P,R)$ induced by precomposition with $(F_X,\mathfrak f)$ is an equivalence of categories for any doctrine $R$.
\end{corollary}
\begin{remark}
We showed how to obtain from the universal 1-arrow $(F_X,\mathfrak f_{(X,\varphi)}):P\to P_{(X,\varphi)}$ for fixed object $X$ and element $\varphi\in P(X)$, both universal 1-arrows $(F_X,\mathfrak f_X):P\to P_X$ in Corollary \ref{coroll:const} for a fixed object $X$ and $(\id{\CC},\mathfrak f_\varphi):P\to P_\varphi$ in Corollary \ref{coroll:ax} for a fixed element $\varphi$ in $P(\tmn)$ as particular cases. Note that we wrote some subscripts to avoid confusion between the constructions. We now show that we can recover the first 1-arrow from the other two. To do so, take a primary doctrine $P:\CC\op\to\Pos$, fix an object $X$ and an element $\varphi\in P(X)$. Apply the construction that adds a constant to obtain $(F_X,\mathfrak{f}_X):P\to P_X$. Now consider the primary doctrine $P_X:\CC_X\op\to\Pos$ and the element $\varphi$ in the fiber over the terminal object $P_X(\tmn)=P(X)$. Apply the construction that adds an axiom to obtain $(\id{\CC_X},\mathfrak{f}_\varphi):P_X\to(P_X)_\varphi$.
\[\begin{tikzcd}
	\CC\op && {\CC_X\op} && {\CC_X\op} \\
	\\
	&& \Pos
	\arrow[""{name=0, anchor=center, inner sep=0}, "{P_X}"{description, pos=0.3}, from=1-3, to=3-3]
	\arrow[""{name=1, anchor=center, inner sep=0}, "P"', from=1-1, to=3-3]
	\arrow["{F_X\op}", from=1-1, to=1-3]
	\arrow["{\id{\CC_X}\op}", from=1-3, to=1-5]
	\arrow[""{name=2, anchor=center, inner sep=0}, "{(P_X)_\varphi}", from=1-5, to=3-3]
	\arrow["{\mathfrak f_X}"', curve={height=-6pt}, shorten <=8pt, shorten >=8pt, from=1, to=0]
	\arrow["{\mathfrak f_\varphi}"', curve={height=-6pt}, shorten <=8pt, shorten >=8pt, from=0, to=2]
\end{tikzcd}\]
Compute for each object $A$, the poset $(P_X)_\varphi(A)=P_X(A)_{\downarrow P_X(!_{X\times A})\varphi}$, where $!_{X\times A}:A\rightsquigarrow\tmn$ is the unique $\CC_X$-arrow from $A$ to $\tmn$. The reindexing along this arrow is $P_X(!_{X\times A}):P_X(\tmn)\to P_X(A)$, that maps $\varphi$ to $P(\pr1)\varphi$, so $(P_X)_\varphi(A)=P(X\times A)_{\downarrow P(\pr1)\varphi}$ which is exactly how the fibers of $P_{(X,\varphi)}$ are computed. Then compute reindexing in $(P_X)_\varphi$: given $f:A\rightsquigarrow B$, we know that $(P_X)_\varphi(f)$ is defined as the restriction of $P_X(f)$, that is $P(\ple{\pr1,f})$, which is how reindexing are computed in $P_{(X,\varphi)}$. So the functor $(P_X)_\varphi$ is $P_{(X,\varphi)}$. Moreover, observe that the composition of the 1-arrows is $(\id{\CC_X},\mathfrak{f}_\varphi)(F_X,\mathfrak{f}_X)=(F_X,\mathfrak{f}_{(X,\varphi)})$.
\end{remark}

\begin{remark}
Using the description of the fibers of $P_{(X,\varphi)}$ in \cref{r:axiom-quotient}, it is worth noticing that the two constructions of adding a constant and adding an axiom have essentially the same shape on the base category and on the fibers. On the one hand, in the new base category, homomorphism between two objects $A$ and $B$ are old homomorphism between $X\times A$ and $B$; more precisely $\Hom_{\CC_X}(A,B)=\Hom_\CC(X\times A,B)$. On the other hand, looking at a poset with finite meets as a category with finite products, we have for any pair $\alpha, \beta\in P(X\times A)$ that $\Hom([\alpha],[\beta])=\Hom(P(\pr1)(\varphi)\land\alpha,\beta)$. Roughly speaking, in each fiber poset, adding a constant of sort $P(\pr1)(\varphi)$ means to force the existence of a (necessarily unique) arrow from the terminal object to $P(\pr1)(\varphi)$, i.e.\ forcing $P(\pr1)(\varphi)$ to become true, hence in our construction we are actually adding an axiom.
\end{remark}

In the following, we apply separately the two constructions to a doctrine of well-formed formulae in some language $\mathcal{L}$ and theory $\mathcal{T}$. At first, we apply the construction that adds a constant to a doctrine, and show that there is an isomorphism between this doctrine and the doctrine of well-formed formulae in the language with a new constant symbol. Then we apply the construction that adds an axiom to a doctrine, and show that there is an isomorphism between this doctrine and the doctrine of well-formed formulae where the theory has a new axiom.
\begin{example}\label{ex:const}
Let $\mathcal{L}$ be a first-order language and $\mathcal{T}$ be a theory. Consider the doctrine of well-formed formulae $\LT^\mathcal{L}_{\mathcal{T}}:\ctx_\mathcal{L}\op\to\Pos$ and the fixed object $(x)$ in the base category $\ctx_\mathcal{L}$. On the one hand, consider the 1-arrow $(F_{(x)},\mathfrak f):\LT^\mathcal{L}_{\mathcal{T}}\to(\LT^\mathcal{L}_{\mathcal{T}})_{(x)}$, where $(\LT^\mathcal{L}_{\mathcal{T}})_{(x)}:(\ctx_\mathcal{L})_{(x)}\op\to\Pos$. Arrows in $(\ctx_\mathcal{L})_{(x)}$ are of the form $\vec t(\,(x);\vec z\,):\vec z\rightsquigarrow\vec y$, and the fibers $(\LT^\mathcal{L}_{\mathcal{T}})_{(x)}(\vec z)$ are $\LT^\mathcal{L}_{\mathcal{T}}(\,(x);\vec z\,)$ for any list of variables $\vec z$. On the other hand consider the doctrine $\LT^{\mathcal{L}\cup\{c\}}_{\mathcal{T}}:\ctx_{\mathcal{L}\cup\{c\}}\op\to\Pos$, where $c$ is a constant symbol not appearing in $\mathcal{L}$. There is a trivial 1-arrow $(E,\mathfrak e):\LT^\mathcal{L}_{\mathcal{T}}\to\LT^{\mathcal{L}\cup\{c\}}_{\mathcal{T}}$: the functor $E$ is defined by the inclusion of terms in the extended language, the natural transformation $\mathfrak e$ is defined by the inclusion of formulae. The universal property of $(F_{(x)},\mathfrak f)$ defines a unique $(E',\mathfrak e'):(\LT^\mathcal{L}_{\mathcal{T}})_{(x)}\to \LT^{\mathcal{L}\cup\{c\}}_{\mathcal{T}}$ such that $(E',\mathfrak e')(F_{(x)},\mathfrak f)=(E,\mathfrak e)$ and such that $E'\big(\id{(x)}:()\rightsquigarrow(x)\big)=\big(c:()\to(x)\big)$.
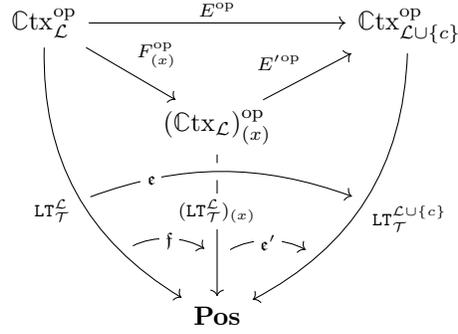
\begin{center}
\begin{tikzcd}
\ctx_\mathcal{L}\op\arrow[dr,"F_{(x)}\op"]\arrow[ddddr,"\LT^\mathcal{L}_{\mathcal{T}}"'{name=f}, bend right,""'{name=a, near end}]\arrow[rr,"E\op"]&&\ctx_{\mathcal{L}\cup\{c\}}\op\arrow[ddddl,"\LT^{\mathcal{L}\cup\{c\}}_{\mathcal{T}}"{name=g}, bend left,""{name=e, near end}]\\
&(\ctx_\mathcal{L})_{(x)}\op\arrow[ddd,"(\LT^\mathcal{L}_{\mathcal{T}})_{(x)}"{description, pos=0.4}, ""{name=b, near end},""'{name=c, near end}]\arrow[ur, "{E'}\op"]\\ \\ \\
&\Pos
\arrow[from=a, to=b, bend left, shorten=2mm, "\mathfrak{f}" description]
\arrow["\mathfrak{e}'"{description}, from=c, to=e, bend left, shorten=2mm]
\arrow[from=f, to=g, bend left=20, shorten=2mm, "\mathfrak{e}"{description, near start}, crossing over]
\end{tikzcd}
\end{center}
The functor $E'$ maps an arrow $\vec t(\,(x);\vec z\,):\vec z\rightsquigarrow\vec y$ to the term $\vec t(\,[c/x];\vec z\,):\vec z\to\vec y$ in $\ctx_{\mathcal{L}\cup\{c\}}$. For a given pair $\vec t(\,(x);\vec z\,), \vec s(\,(x);\vec z\,)$ such that $\vec t(\,[c/x];\vec z\,)=\vec s(\,[c/x];\vec z\,)$, substitute again $[x/c]$ and get $\vec t=\vec s$, so $E'$ is faithful. Then, for a given a term $u(\vec z)$ in the language ${\mathcal{L}\cup\{c\}}$, we can consider $c$ as a variable and substitute each occurrence of $c$ with $x$, to obtain a term $u'(\,(x);\vec z\,)$ obviously written in the language $\mathcal L$: in particular $E'(u')=u'(\,[c/x];\vec z\,)=u(\vec z)$, so $E'$ is full. Moreover, since $E'$ is the identity on objects, $E'$ is an isomorphism.

Concerning formulae, a component of the natural transformation $\mathfrak{e}'_{(\vec z)}$ sends a formula $\alpha(\,(x);\vec z\,)$ in $(\LT^\mathcal{L}_{\mathcal{T}})_{(x)}(\vec z)=\LT^\mathcal{L}_{\mathcal{T}}(\,(x);\vec z\,)$ to the formula $\alpha(\,[c/x];\vec z\,)\in\LT^{\mathcal{L}\cup\{c\}}_{\mathcal{T}}(\vec z)$. A similar argument to the one that showed fullness of the functor $E'$ proves that $\mathfrak{e}'$ is a natural isomorphism.

To conclude, we can say that the doctrine $(\LT^\mathcal{L}_{\mathcal{T}})_{(x)}$ is again a doctrine of well-formed formulae.
\end{example}
\begin{example}\label{ex:ax}
Let $\mathcal{L}$ be a first-order language and $\mathcal{T}$ be a theory. Consider the doctrine of well-formed formulae $\LT^\mathcal{L}_{\mathcal{T}}:\ctx_\mathcal{L}\op\to\Pos$ and the fixed $\mathcal{L}$-sentence $\varphi\in\LT^\mathcal{L}_{\mathcal{T}}()$. On the one hand, consider the 1-arrow $(\id{\ctx_\mathcal{L}},\mathfrak f):\LT^\mathcal{L}_{\mathcal{T}}\to(\LT^\mathcal{L}_{\mathcal{T}})_{\varphi}$, where $(\LT^\mathcal{L}_{\mathcal{T}})_{\varphi}:\ctx_\mathcal{L}\op\to\Pos$. Its fibers $(\LT^\mathcal{L}_{\mathcal{T}})_{\varphi}(\vec z)$ are by definition $\LT^\mathcal{L}_{\mathcal{T}}(\vec z)_{\downarrow\varphi}$ for any list of variables $\vec z$. On the other hand consider the doctrine $\LT^\mathcal{L}_{\mathcal{T}\cup\{\varphi\}}:\ctx_{\mathcal{L}}\op\to\Pos$. There is an obvious 1-arrow $(\id{\ctx_\mathcal{L}},\mathfrak e):\LT^\mathcal{L}_{\mathcal{T}}\to\LT^{\mathcal{L}}_{\mathcal{T}\cup\{\varphi\}}$: the natural transformation $\mathfrak e$ is defined by the quotient of formulae with respect to the extended theory, meaning that for each component $\vec x$ it maps any $\mathcal{T}$-provable sequent $\alpha(\vec x)\vdash_{\mathcal T}\beta(\vec x)$ into the $\mathcal{T}\cup\{\varphi\}$-provable sequent $\alpha(\vec x)\vdash_{\mathcal{T}\cup\{\varphi\}}\beta(\vec x)$. To use the universal property of $(\id{\ctx_\mathcal{L}},\mathfrak f)$, we need to check that $\mathfrak{e}_{()}$ maps $\varphi\in\LT^\mathcal{L}_{\mathcal{T}}()$ to the top element of $\LT^\mathcal{L}_{\mathcal{T}\cup\{\varphi\}}()$. However this is true since clearly $\top\vdash_{\mathcal{T}\cup\{\varphi\}}\varphi$. Consequently there exists a unique $(E',\mathfrak e'):(\LT^\mathcal{L}_{\mathcal{T}})_{\varphi}\to \LT^{\mathcal{L}}_{\mathcal{T}\cup\{\varphi\}}$ such that $(E',\mathfrak e')(\id{\ctx_\mathcal{L}},\mathfrak f)=(\id{\ctx_\mathcal{L}},\mathfrak e)$.
\begin{center}
\begin{tikzcd}
\ctx_\mathcal{L}\op\arrow[dr,"\id{\ctx_\mathcal{L}}\op"]\arrow[ddddr,"\LT^\mathcal{L}_{\mathcal{T}}"'{name=f}, bend right,""'{name=a, near end}]\arrow[rr,"\id{\ctx_\mathcal{L}}\op"]&&\ctx_\mathcal{L}\op\arrow[ddddl,"\LT^{\mathcal{L}}_{\mathcal{T}\cup\{\varphi\}}"{name=g}, bend left,""{name=e, near end}]\\
&\ctx_\mathcal{L}\op\arrow[ddd,"(\LT^\mathcal{L}_{\mathcal{T}})_{\varphi}"{description, pos=0.4}, ""{name=b, near end},""'{name=c, near end}]\arrow[ur, "\id{\ctx_\mathcal{L}}\op"]\\ \\ \\
&\Pos
\arrow[from=a, to=b, bend left, shorten=2mm, "\mathfrak{f}" description]
\arrow["\mathfrak{e}'"{description}, from=c, to=e, bend left, shorten=2mm]
\arrow[from=f, to=g, bend left=20, shorten=2mm, "\mathfrak{e}"{description, near start}, crossing over]
\end{tikzcd}
\end{center}
The functor $E'$ is the identity.

A component of the natural transformation $\mathfrak{e}'_{\vec{x}}$ sends a formula $\alpha(\vec x)$ in $(\LT^\mathcal{L}_{\mathcal{T}})_{\varphi}(\vec x)=\LT^\mathcal{L}_{\mathcal{T}}(\vec x)_{\downarrow\varphi}$ to the formula $\alpha(\vec x)\in\LT^{\mathcal{L}}_{\mathcal{T}\cup\{\varphi\}}(\vec x)$. Define the inverse function: it maps $\beta(\vec x)$ to $\beta(\vec x)\land \varphi$. This is well defined and monotone, since if we take $\alpha(\vec x)\vdash_{\mathcal{T}\cup\{\varphi\}}\beta(\vec x)$, it easily follows that $\alpha(\vec x)\land\varphi\vdash_{\mathcal{T}}\beta(\vec x)\land\varphi$. On the one hand, take $\alpha(\vec x)$ such that $\alpha(\vec x)\vdash_{\mathcal{T}}\varphi$, apply $\mathfrak{e}'_{\vec{x}}$ to get $\alpha(\vec x)$, and then send it to $\alpha(\vec x)\land\varphi$, and observe that $\alpha(\vec x)\land\varphi\dashv\vdash_{\mathcal{T}}\alpha(\vec x)$ using the initial assumption on $\alpha(\vec x)$. Conversely, take $\beta(\vec x)\in\LT^{\mathcal{L}}_{\mathcal{T}\cup\{\varphi\}}(\vec x)$, send it to $\beta(\vec x)\land\varphi$, and then apply $\mathfrak{e}'_{\vec{x}}$ to get $\beta(\vec x)\land\varphi\in \LT^{\mathcal{L}}_{\mathcal{T}\cup\{\varphi\}}(\vec x)$. Observe that $\beta(\vec x)\land\varphi\dashv\vdash_{\mathcal{T}\cup\{\varphi\}}\beta(\vec x)$. So $\mathfrak{e}'$ is indeed a natural isomorphism.

To conclude, we can say that the doctrine $(\LT^\mathcal{L}_{\mathcal{T}})_{\varphi}$ is again a doctrine of well-formed formulae.
\end{example}

Going towards a conclusion, in the following example we prove that if we apply our construction to the subsets doctrine $\mathscr{P}$ trying to force the existence of an arrow in the empty set, the resulting doctrine collapses to a trivial one, meaning that all of its fibers are singletons.
\begin{example}\label{ex:collapse}
Consider the doctrine $\mathscr{P}:\Set\op\to\Pos$ and do the construction with respect to the object $\emptyset$.
\begin{center}
\begin{tikzcd}
\Set^{\text{op}}\arrow[rr,"{F_\emptyset}^{\text{op}}"] \arrow[dr,"\mathscr P"' ,""{name=L}]&&{\Set_\emptyset}^{\text{op}}\arrow[dl,"\mathscr{P}_{\emptyset}" ,""'{name=R}]\\
&\Pos\arrow[rightarrow,"\mathfrak{f}","\cdot"', from=L, to=R, bend left=10]
\end{tikzcd}
\end{center}
In the category $\Set_\emptyset$ all objects are isomorphic: indeed, for any pair of sets $A,B$, there exists exactly one homomorphism $A\rightsquigarrow B$ corresponding to the inclusion $\emptyset\times A\to B$ of the empty set in $B$. The fiber over a set $A$ is $\mathscr{P}_\emptyset(A)=\mathscr{P}(\emptyset\times A)=\mathscr{P}(\emptyset)=\{\emptyset\}$, hence every fiber is trivial. In particular, the natural transformation $\mathfrak f$ is the constant function on every component.
\end{example}

A similar problem arises if we apply our construction again to the subset doctrine $\mathscr{P}$ trying to force the empty set $\emptyset\in\mathscr{P}(\{*\})$ as a new axiom. While the base category $\Set$ is the same after the construction, the fibers all collapse as in \cref{ex:collapse}.
\begin{example}
Consider the doctrine $\mathscr{P}:\Set\op\to\Pos$ and do the construction with respect to the element $\emptyset\in\mathscr{P}(\{*\})$.
\begin{center}
\begin{tikzcd}
\Set^{\text{op}}\arrow[rr,"{\id{}}"] \arrow[dr,"\mathscr P"' ,""{name=L}]&&{\Set}^{\text{op}}\arrow[dl,"\mathscr{P}_{\emptyset\subseteq\{*\}}" ,""'{name=R}]\\
&\Pos\arrow[rightarrow,"\mathfrak{f}","\cdot"', from=L, to=R, bend left=10]
\end{tikzcd}
\end{center}
Note that we wrote $\mathscr{P}_{\emptyset\subseteq\{*\}}$to remark the fact that here the empty set is an element in the fiber over the terminal object of $\Set$, and not an object in the base cateogry as it was in \cref{ex:collapse}.
The fiber over a set $A$ is $\mathscr{P}_{\emptyset\subseteq\{*\}}(A)=\mathscr{P}(A)_{\downarrow\emptyset}=\{S\subseteq A\mid S\subseteq\emptyset\}=\{\emptyset\}$, hence every fiber is trivial. Also in this case, the natural transformation $\mathfrak f$ is the constant function on every component.
\end{example}

\begin{definition}
Let $P:\CC\op\to\Pos$ and $R:\ct{D}\op\to\Pos$ be two doctrines and let $(F,\mathfrak f):P\to R$ be a model of $P$ in $R$. The model $(F,\mathfrak f)$ is \emph{conservative} if every component of $\mathfrak f$ is full as a functor between posets.
\end{definition} 
The idea is that the model $(F_X,\mathfrak f):P\to P_{(X,\varphi)}$ is conservative if, for two given formulae $\alpha, \beta$ in a fiber of the doctrine $P$---in the old language---, if $\mathfrak{f}_A\alpha\leq\mathfrak{f}_A\beta$---the interpretation of $\alpha$ proves the interpretation of $\beta$ in the new theory of the extended language---, then $\alpha\leq\beta$ in $P(A)$---$\alpha$ was already proving $\beta$ in the old theory.
The examples above show that in general the morphism $(F_X,\mathfrak f):P\to P_{(X,\varphi)}$ is not conservative. However, if the doctrine $P$ is existential, we have the following characterization.

\begin{proposition}
Let $P:\CC\op\to\Pos$ be an existential doctrine and $P_{(X,\varphi)}$ be the Kleisli object of the comonad $(P,(X\times\blank,\mathfrak f),\gamma,\varepsilon)$ defined by the pair $X\in\CC$ and $\varphi\in P(X)$. The model $(F_X,\mathfrak f):P\to P_{(X,\varphi)}$ is conservative if and only if $\top_\tmn\leq\exists^X_\tmn\varphi$ in $P(\tmn)$.
\begin{proof}
First of all, suppose that for every object $A$ and for every pair $\alpha,\beta\in P(A)$, if $\mathfrak{f}_A\alpha\leq\mathfrak{f}_A\beta$ then $\alpha\leq \beta$. Consider the unit of the adjunction $\varphi\leq P(!_X)\exists^X_\tmn\varphi$ in $P(X)$. However, apply $\mathfrak f_\tmn:P(\tmn)\to P_{(X,\varphi)}(\tmn)=P(X)_{\downarrow\varphi}$ to both $\top_\tmn$ and $\exists^X_\tmn\varphi$: we get
\begin{align*}
\mathfrak f_\tmn(\top_\tmn)&=P(!_X)(\top_\tmn)\land\varphi=\varphi\\
\mathfrak f_\tmn(\exists^X_\tmn\varphi)&=P(!_X)(\exists^X_\tmn\varphi)\land\varphi.\end{align*}
So from $\mathfrak{f}_\tmn\top\leq\mathfrak{f}_\tmn\exists^X_\tmn\varphi$, by fullness we get $\top_\tmn\leq\exists^X_\tmn\varphi$.

Conversely, suppose that $\top_\tmn\leq\exists^X_\tmn\varphi$ in $P(\tmn)$ and $\mathfrak f_A(\alpha)\leq\mathfrak f_A(\beta)$ in $P_{(X,\varphi)}(A)$, for some pair $\alpha,\beta\in P(A)$, or equivalently $P(\pr1)\varphi\land P(\pr2)\alpha\leq P(\pr2)\beta$ in $P(X\times A)$. We switch the coordinates and work in $P(A\times X)$: we get $P(\pr2)\varphi\land P(\pr1)\alpha\leq P(\pr1)\beta$. Then use the adjunction $\exists^X_A\dashv P(\pr1)$ to obtain $\exists^X_A\big(P(\pr2)\varphi\land P(\pr1)\alpha\big)\leq\beta$ in $P(A)$. Using Frobenius reciprocity and Beck-Chevalley condition we then get $\exists^X_A\big(P(\pr2)\varphi\land P(\pr1)\alpha\big)=\exists^X_A P(\pr2)\varphi\land\alpha=P(!_A){\exists^X_\tmn\varphi}\land\alpha=P(!_A)\top_\tmn\land\alpha=\alpha$, where we used the assumption in the third equality. Hence $\alpha\leq\beta$ as claimed.\end{proof}\end{proposition}

\begin{example}\label{ex:costruzione-parti}
Consider the doctrine $\mathscr{P}:\Set\op\to\Pos$ and do the construction with respect to the object $X$.
\begin{center}
\begin{tikzcd}
\Set^{\text{op}}\arrow[rr,"{F_X}^{\text{op}}"] \arrow[dr,"\mathscr P"' ,""{name=L}]&&{\Set_X}^{\text{op}}\arrow[dl,"\mathscr{P}_{X}" ,""'{name=R}]\\
&\Pos\arrow[rightarrow,"\mathfrak{f}","\cdot"', from=L, to=R, bend left=10]
\end{tikzcd}
\end{center}
In the category $\Set_X$ objects are sets, and an arrow $A\rightsquigarrow B$ is a function $X\times A\to B$. The fiber over a set $A$ is $\mathscr{P}_X(A)=\mathscr{P}(X\times A)=\mathscr{P}(X\times A)$, and the reindexing along $f:A\rightsquigarrow B$ maps $S\subseteq X\times B$ to $\{(x,a)\in X\times A\mid (x,f(x,a))\in S\}$. The natural transformation $\mathfrak f$ at an object $A\in \Set$ maps $S\subseteq A$ to $X\times S\subseteq X\times A$.
\end{example}

\begin{example}\label{ex:costruzione-parti-2}
Consider the doctrine $\mathscr{P}:\Set\op\to\Pos$ and do the construction with respect to the object $X$ and the element $Y\in\mathscr{P}(X)$.
\begin{center}
\begin{tikzcd}
\Set^{\text{op}}\arrow[rr,"{F_X}^{\text{op}}"] \arrow[dr,"\mathscr P"' ,""{name=L}]&&{\Set_X}^{\text{op}}\arrow[dl,"\mathscr{P}_{(X,Y)}" ,""'{name=R}]\\
&\Pos\arrow[rightarrow,"\mathfrak{f}","\cdot"', from=L, to=R, bend left=10]
\end{tikzcd}
\end{center}
In the category $\Set_X$ is defined as in \cref{ex:costruzione-parti}. The fiber over a set $A$ is
\[\mathscr{P}_{(X,Y)}(A)=\mathscr{P}(X\times A)_{\downarrow \pr1^{-1}[Y]}=\{S\subseteq X\times A\mid S\subseteq Y\times A\}=\mathscr{P}(Y\times A),\]
and the reindexing along $f:A\rightsquigarrow B$ maps $S\subseteq Y\times B$ to $\{(x,a)\in Y\times A\mid (x,f(x,a))\in S\}$. The natural transformation $\mathfrak f$ at an object $A\in \Set$ maps $S\subseteq A$ to $Y\times S=(X\times S)\cap(Y\times A)\subseteq Y\times A$.
\end{example}

\begin{remark}
Comparing \cref{ex:costruzione-parti,ex:costruzione-parti-2} above, we can observe that for a given pair of sets $X,Y$ with $Y\subseteq X$, the construction that adds a constant of sort $Y$ to the subset doctrine as in \cref{ex:costruzione-parti}, and the construction that adds a constant of sort $X$ and then the axiom $Y$ evaluated in the new constant as in \cref{ex:costruzione-parti-2}, have the same fibers: in both cases, the fiber over a set $A$ is $\mathscr{P}(Y\times A)$. The only difference between the two constructions is in the arrows of the base categories. However, by the universal property of $\mathscr{P}\to\mathscr{P}_{(X,Y)}$, we have the following factorisation of the doctrine homomorphism $\mathscr{P}\to\mathscr{P}_Y$:

\begin{equation}\label{d:costruzione-parti}
\begin{tikzcd}
\Set\op\arrow[dr,"F_X\op"]\arrow[ddddr,"\mathscr{P}"'{name=f}, bend right,""'{name=a, near end}]\arrow[rr,"F_Y\op"]&&\Set_Y\op\arrow[ddddl,"\mathscr{P}_Y"{name=g}, bend left,""{name=e, near end}]\\
&\Set_X\op\arrow[ddd,"\mathscr{P}_{(X,Y)}"{description, pos=0.4}, ""{name=b, near end},""'{name=c, near end}]\arrow[ur, dashed]\\ \\ \\
&\Pos
\arrow[from=a, to=b, bend left, shorten=2mm, "\mathfrak{f}_{X,Y}"{pos=0.4}]
\arrow[""{description}, from=c, to=e, bend left, shorten=2mm, dashed]
\arrow[from=f, to=g, bend left=20, shorten=2mm, "\mathfrak{f}_Y"{description, near start}, crossing over]
\end{tikzcd}
\end{equation}
Indeed, consider the arrow $\iota:\tmn\rightsquigarrow X$ in $\Set_Y$ corresponding to the inclusion $\iota:Y\to X$ in $\Set$. The arrow $\iota$ is such that $\top\leq\mathscr{P}_Y(\iota)(\mathfrak{f}_Y)_X(Y)$ in $\mathscr{P}_Y(\{*\})$, since $Y$ is the top element of $\mathscr{P}_Y(\{*\})$.
\[\begin{tikzcd}[row sep=tiny]
	{\mathscr{P}(X)} & {\mathscr{P}_Y(X)=\mathscr{P}(Y\times X)} & {\mathscr{P}_Y(\{*\})=\mathscr{P}(Y)} \\
	Y & {Y\times Y} & {\{y\in Y\mid (y,y)\in Y\times Y\}=Y.}
	\arrow["{(\mathfrak{f}_Y)_X}", from=1-1, to=1-2]
	\arrow["{\ple{\id{Y},\iota}^{-1}}", from=1-2, to=1-3]
	\arrow[maps to, from=2-1, to=2-2]
	\arrow[from=2-2, to=2-3]
\end{tikzcd}\]
Thus, by \cref{thm:univ_prop}, we have the factorisation in \eqref{d:costruzione-parti}. The functor $\Set_X\to\Set_Y$ between the base categories is the identity on objects, and maps a function $g:X\times A\to B$ to its restriction $g(\iota\times\id{A}):Y\times A\to B$; the natural transformation is the identity on each fiber.
\end{remark}

\begin{example}
We provide an interpretation of the universal property in \cref{coroll:ax}. Let $\mathcal{L}=\{\cdot, e,(\blank)^{-1}\}$ be the language of groups, let $\mathcal{T}$ be the theory of groups and let $\varphi\in \LT^\mathcal{L}_\mathcal{T}()$ be the commutativity axiom $\forall x\forall y (x\cdot y = y \cdot x)$. Set-teoretic models for the theory $\mathcal{T}$, i.e.\  groups, are in a one-to-one correspondence with doctrine homomorphisms $\LT^\mathcal{L}_\mathcal{T}\to\mathscr{P}$ preserving boolean, existential and elementary structure \cite{quotcomplfoun}. Similarly, abelian groups are in a one-to-one correspondence with boolean existential elementary doctrine homomorphisms $\LT^\mathcal{L}_{\mathcal{T}\cup\{\varphi\}}\to\mathscr{P}$. As we saw in \cref{ex:ax}, the doctrine $\LT^\mathcal{L}_{\mathcal{T}\cup\{\varphi\}}$ is obtained from the doctrine $\LT^\mathcal{L}_\mathcal{T}$ by adding the axiom $\varphi$. So let $H$ be a group and let $(G,\mathfrak{g}):\LT^\mathcal{L}_\mathcal{T}\to\mathscr{P}$ be the correspondent doctrine homomorphism. The group $H$ is abelian if and only if $(G,\mathfrak{g})$ factors through the extension morphism $\LT^\mathcal{L}_\mathcal{T}\to(\LT^\mathcal{L}_\mathcal{T})_{\varphi}$. This is a way to read what is stated in \cref{coroll:ax}: roughly speaking, given a primary doctrine $P$ and an element $\varphi\in P(\tmn)$, a model $(G,\mathfrak{g}):P\to R$ satisfies $\varphi$ if and only if $\top\leq\mathfrak{g}_\tmn(\varphi)$ in $R(\tmn)$ if and only if the model $(G,\mathfrak{g})$ factors through $P\to P_\varphi$. 
\end{example}

\section*{Acknowledgement}
The author would like to acknowledge her Ph.D.\ supervisors Sandra Mantovani and Pino Rosolini: this paper is the second chapter of her Ph.D.\ thesis, defended at Universit\`a degli Studi di Milano. The author would also like to thank the anonymous referee for his remarks and suggestions, and Jacopo Emmenegger for his valuable comments.

\section*{Conflict of interest statement}
The author declares no conflicts of interest.

\vskip 1cm
\bibliography{Biblio}
\bibliographystyle{alpha}
\end{document}